\documentclass[12pt]{article}
\usepackage{amsmath}
\usepackage{amsfonts}
\usepackage{amssymb}
\usepackage{theorem} 

\title{Polish $G$-spaces and continuous logic}
\author{A. Ivanov
\thanks{The first author is supported by Polish National Science Centre grant DEC2011/01/B/ST1/01406} 
, B. Majcher-Iwanow
\thanks{The second author is partially supported by ESF Short Visit Grant no. 5419} }
\date{ } 

\newtheorem{thm}{Theorem}[section]
\newtheorem{lem}[thm]{Lemma}
\newtheorem{definicja}[thm]{Definition}
\newtheorem{cor}[thm]{Corollary}
\newtheorem{prop}[thm]{Proposition}
\newtheorem{remark}[thm]{Remark}  

\begin{document}
\maketitle
\topskip 20pt

\begin{quote}
{\bf Abstract.} 
We analyse logic actions of  Polish groups which arise  
in continuous logic. 
We extend the generalised model theory of H.Becker from 
\cite{becker2} to the case of Polish $G$-spaces when $G$ 
is an arbitrary Polish group. 
\bigskip


\end{quote}

\bigskip

\section{Introduction}

The paper is devoted to Polish group actions 
from the point of view of continuous logic. 
Let $({\bf Y},d)$ be a Polish space and $Iso({\bf Y})$ 
be the corresponding isometry group 
endowed with the pointwise convergence topology. 
Then $Iso ({\bf Y})$ is a Polish group. \parskip0pt 

It is worth noting that any Polish group $G$ can be 
realised as a closed subgroup of the isometry group 
$Iso ({\bf Y})$ of an appropriate Polish space  $({\bf Y},d)$. 
Moreover it is shown by J.Melleray in \cite{melleray} 
(Theorem 6) that $G$ can be chosen as the automorphism 
group of a continuous metric structure on $({\bf Y},d)$ 
which is approximately ultrahomogeneous.  

For any countable continuous signature $L$ the set 
${\bf Y}_L$ of all continuous metric $L$-structures on 
$({\bf Y},d)$ can be considered as a Polish 
$Iso({\bf Y})$-space. 
We call this action {\em logic} and show that it is universal 
for Borel reducibility of orbit equivalence relations 
of Polish $G$-spaces with closed $G\le Iso({\bf Y})$.  

Note that for any tuple $\bar{s}\in {\bf Y}$ the map 
$g\rightarrow d(\bar{s},g(\bar{s}))$ can be considered 
as a {\em grey} subgroup of $G$. 
Grey subsets  and subgroups of metric spaces (groups) 
were introduced in \cite{BYBM} as {\em graded} counterparts 
of subsets and subgroups. 
For any continuous sentence $\phi$ the map 
${\bf Y}_L \rightarrow [0,1]$ defined by $M\rightarrow \phi^M$ 
is a grey subset of  ${\bf Y}_L$. 

Typical notions naturally arising for logic actions  
can be applied in the general case of a Polish 
$G$-space ${\bf X}$ with $G$ as above. 
If we consider $G$ together with a family of grey 
subgroups as above, then distinguishing an appropriate family 
$\mathcal{B}$ of grey subsets of ${\bf X}$ we arrive at 
the situation very similar to the logic space ${\bf Y}_L$. 
For example we can treat elements of $\mathcal{B}$ as 
continuous formulas. 

In the case of the logic action of closed subgroups of 
$S_{\infty}$ on the space of $X_L$ of discrete structures 
on $\omega$ this approach was realised by H.Becker in 
\cite{becker} and \cite{becker2}. 
Imitating topologies generated by sets of the form 
$Mod(\phi,\bar{s})=\{ M\in X_L : M\models \phi (\bar{s})\}$, 
where $\bar{s}\in \omega$ and formulas $\phi$ are chosen 
from a countable fragment of $L_{\omega_1 \omega}$, 
he introduces the concept of a {\em nice topology}. 
Many theorems of traditional model theory can be 
generalised to topological statements concerning 
spaces with nice topologies. 
  
We generalise this approach to the general case 
of Polish $G$-actions removing the Becker's assumption 
that $G<S_{\infty}$. 
It turns out that some basic notions from \cite{becker2} 
have to be replaced by their grey counterparts. 
In particular we define {\em grey Vaught transforms} and   
{\em grey nices bases}. 
In Section 2 we prove an existence theorem  
stating that under natural circumstances appropriately 
defined nice topologies can be always found. 

When a family of grey subsets $\mathcal{B}$ 
generates such a topology as a nice basis  
(this will be described later) then treating 
elements of $\mathcal{B}$ as continuous formulas, 
we may extend theorems of continuous model theory 
to some topological statements.  
In Section 4 we obtain topological versions 
of several theorems from logic, 
for example Ryll-Nardzewski's theorem.  

The construction of nice bases arises in the 
most natural form when one considers the case 
of continuous logic actions over $\mathbb{U}$, 
the Urysohn space of diameter 1. 
Let $L$ be a countable continuous signature and 
$\mathbb{U}_L$ be the $Iso(\mathbb{U})$-space of 
all continuous $L$-structures. 
Our main theorem of Section 3 rouphly states that 
for any countable fragment $\mathcal{L}$ of 
$L_{\omega_1 \omega}$, grey subsets associated with continuous 
$\mathcal{L}$-formulas form a nice basis on $\mathbb{U}_L$.  

We may interpret this theorem that $\mathbb{U}$ 
is the continuous counterpart of $\omega$ 
in moving from the case of discrete logic 
$S_{\infty}$-actions to the case of actions 
on spaces of continuous structures.  
Some other reasons supporting this intuition 
can be found in  \cite{EFP} and in Section 2 of \cite{CL}. 

This motivates some very basic questions. 
A theorem  from \cite{becker2} (Corollary 1.13) 
states that in the case of logic actions of 
$S_{\infty}$ on the space of countable structures 
each nice topology is defined by model sets 
$Mod (\phi )$ of formulas of a countable fragment 
of $L_{\omega_1 \omega}$.  
Can this statement be extended to $\mathbb{U}_L$?   

In Section \ref{LopEsc} we give the positive answer. 
In fact this  shows that nice topologies induced 
by countable continuous fragments  
of $L_{\omega_1 \omega}$ are ubiquitous.  
This issue depends on the L\'{o}pez-Escobar theorem 
on invariant Borel subsets of logic spaces.  
Note that the classical   L\'{o}pez-Escobar theorem was 
one of the ingredients of Corollary 1.13 of \cite{becker2}. 
We use the version of the L\'{o}pez-Escobar theorem 
for continuous logic obtained by S.Coskey and M.Lupini 
in \cite{CL} 
.  
Very recently I.Ben Yaacov, A.Nies and 
T.Tsankov have proved  in \cite{BNT} 
another natural continuous version of 
the L\'{o}pez-Escobar theorem.   

Viewing the logic space ${\bf Y}_L$ 
as a Polish space one may consider Borel/algorithmic 
complexity of interesting subsets of ${\bf Y}_L$. 
This is the main concern of Section 5. 
It demonstrates some new setting arising in the 
approach of the logic space of continuous structures.  


We now give some preliminaries in detail.

\paragraph{Polish group actions.}

A {\em Polish space (group)} is a separable, completely
metrizable topological space (group).
Sometimes we extend the corresponding metric to 
tuples by 
$$
d((x_1 ,...,x_n ), (y_1 ,...,y_n ))= max (d(x_1 ,y_1 ),...,d(x_n ,y_n )). 
$$ 
If a Polish group $G$ continuously acts on a Polish space ${\bf X}$,
then we say that ${\bf X}$ is a {\em Polish $G$-space}. 
We  say that a subset of ${\bf X}$ is {\em invariant} if
it is $G$-invariant. \parskip0pt

If $B$ is a subset of ${\bf X}$ and $u$ is a non-empty open 
subset of $G$ then let 
$$
B^{\star u} = \{ x\in X:\{ g\in u:gx\in B\} \mbox{ is comeagre in }u\} ,  
$$ 
$$ 
B^{\Delta u}=\{ x\in X:\{ g\in u:gx\in B\} \mbox{ is not meagre in }u\} .
$$
These operations are called Vaught transforms. 
Their properties can be found in Section 5 of  \cite{bk}.

Let $({\bf Y},d)$ be a Polish space and $Iso({\bf Y})$ 
be the corresponding isometry group 
endowed with the pointwise convergence toplogy. 
Then $Iso ({\bf Y})$ is a Polish group. 
A compatible left-invariant metric can be obtained as follows: 
fix a countable dense set $S=\{ s_i : i\in \{ 1,2,...\} \}$ 
and then define for two isometries $\alpha$ 
and $\beta$ of ${\bf Y}$ 
$$ 
\rho_{S} (\alpha ,\beta )= \sum_{i=1}^{\infty} 2^{-i} min(1, d(\alpha (s_i ),\beta (s_i ))) .
$$ 
The metric completion of $(Iso({\bf Y},d), \rho_{S} )$ 
can be naturally considered as a semigroup of isometric 
embeddings of $({\bf Y},d)$ into itself. 
Let $In({\bf Y})$ be the semigroup of all isometric 
embeddings of this space. 
 
We will study closed subgroups of $Iso({\bf Y})$.  
We fix a dense countable set $\Upsilon \subset Iso({\bf Y})$. 
In any closed subgroups of $Iso({\bf Y},d)$ we distinguish 
the base consisting of all sets of the form 
$N_{\sigma ,q} = \{ \alpha : \rho_S (\alpha, \sigma )<q \}$, $\sigma \in \Upsilon$ 
and $q\in \mathbb{Q}$.

\paragraph{Continuous structures.} 

We now fix a countable continuous signature 
$$
L=\{ d,R_1 ,...,R_k ,..., F_1 ,..., F_l ,...\} . 
$$ 
Let us recall that a {\em metric $L$-structure} 
is a complete metric space $(M,d)$ with $d$ bounded by 1, 
along with a family of uniformly continuous operations on $M$ 
and a family of predicates $R_i$, i.e. uniformly continuous maps 
from appropriate $M^{k_i}$ to $[0,1]$.   
It is usually assumed that to a predicate symbol $R_i$ 
a continuity modulus $\gamma_i$ is assigned so that when 
$d(x_j ,x'_j ) <\gamma_i (\varepsilon )$ with $1\le j\le k_i$ 
the corresponding predicate of $M$ satisfies 
$$ 
|R_i (x_1 ,...,x_j ,...,x_{k_i}) - R_i (x_1 ,...,x'_j ,...,x_{k_i})| < \varepsilon . 
$$ 
It happens very often that $\gamma_i$ coincides with $id$. 
In this case we do not mention the appropriate modulus. 
We also fix continuity moduli for functional symbols. 

Note that each countable structure can be considered 
as a complete metric structure with the discrete $\{ 0,1\}$-metric.  

Atomic formulas are the expressions of the form $R_i (t_1 ,...,t_r )$, 
$d(t_1 ,t_2 )$, where $t_i$ are terms (built from functional $L$-symbols). 
In metric structures they can take any value from $[0,1]$.   
{\em Statements} concerning metric structures are usually 
formulated in the form 
$$
\phi = 0, 
$$ 
where $\phi$ is a {\em formula}, i.e. an expression built from 
0,1 and atomic formulas by applications of the following functions: 
$$ 
x/2  \mbox{ , } x\dot- y= max (x-y,0) \mbox{ , } min(x ,y )  \mbox{ , } max(x ,y )
\mbox{ , } |x-y| \mbox{ , } 
$$ 
$$ 
\neg (x) =1-x \mbox{ , } x\dot+ y= min(x+y, 1) \mbox{ , } sup_x \mbox{ and } inf_x . 
$$ 
Sometimes statements are called {\em conditions}; we will use both names. 
A {\em theory} is a set of statements without free variables 
(here $sup_x$ and $inf_x$ play the role of quantifiers). 

We often extend the set of formulas by the application 
of {\em truncated products} by positive rational numbers.   
This means that when $q\cdot x$ is greater than $1$, the 
truncated product of $q$ and $x$ is $1$. 
Since the context is always clear, we preserve the same notation $q\cdot x$. 
The continuous logic after this extension does not differ from 
the basic case. 
   
It is worth noting that the choice of the set of connectives 
guarantees that for any continuous relational structure $M$, 
any formula $\phi$ is a $\gamma$-uniform 
continuous function from the 
appropriate power of $M$ to $[0,1]$, where 
$\gamma (\varepsilon )$ is of the form  
$$ 
\frac{1}{n} \cdot  
min \{ \gamma ' (\varepsilon ) : \gamma ' 
\mbox{ is a continuity modulus of an } 
L\mbox{-symbol appearing in the formula} \} , 
$$ 
$$ 
\mbox{ where the number $n$ only depends on the complexity of  }\phi . 
$$ 
This follows from the fact that when $\phi_1$ and $\phi_2$ have 
continuity moduli $\gamma_1$ and $\gamma_2$ respectively, 
then the formula $f(\phi_1 ,\phi_2 )$ obtained by applying 
a binary connective $f$, has a continuity modulus 
of the form $min (\gamma_1 (\frac{1}{2} x) , \gamma_2 (\frac{1}{2} x) )$. 

It is observed in Appendix A of \cite{BYU} that instead of continuity moduli 
one can consider {\em inverse continuity moduli}. 
Slightly modifying that place in \cite{BYU} we define it as follows. 

\begin{definicja} \label{inverse} 
A continuous monotone function $\delta :[0,1]\rightarrow [0,1]$ with $\delta(0) =0$ 
 is an inverse  continuity modulus of a map $F(\bar{x}) :  {\bf X}^n \rightarrow [0,1]$ 
if for any $\bar{a}$, $\bar{b}$ from ${\bf X}^n$, 
$$ 
|F(\bar{a} )-F(\bar{b})| \le \delta (max_{i\le n} (d(a_i ,b_i )))
$$ 
\end{definicja}  

The choice of the connectives above guarantees that 
the following statement holds. 

\begin{lem} \label{ContMod} 
For any continuous relational structure $M$, 
where any $n$-ary relation has $n\cdot id$ as an inverse 
continuity modulus, any formula $\phi$ admits an inverse continuity 
modulus which is of the form $k \cdot id$, where $k$ 
depends on the complexity of $\phi$. 
\end{lem} 

The proof is easy. 
For example it is applied in the proof that when $\phi_1$ and $\phi_2$ 
have inverse continuity moduli $\delta_1$ and $\delta_2$ respectively, 
then the formula $f(\phi_1 ,\phi_2 )$ obtained by applying 
a binary connective $f$, has an inverse continuity modulus 
of the form $2 \cdot max (\delta_1 , \delta_2 )$.  
This lemma will be used below in several arguments. 

For a continuous structure $M$ defined on $({\bf Y},d)$ let 
$Aut(M)$ be the subgroup of $Iso ({\bf Y},d)$ consisting of 
all isometries preserving the values of atomic formulas. 
It is easy to see that $Aut(M)$ is a closed subgroup with 
respect to the topology on $Iso ({\bf Y})$ defined above. 

For every $c_1 ,...,c_n \in M$ and $A\subseteq M$ 
we define the $n$-type $tp(\bar{c}/A)$ of $\bar{c}$ over $A$ 
as the set of all $\bar{x}$-conditions with parameters from $A$ 
which are satisfied by $\bar{c}$ in $M$.  
Let $S_n (T_A )$ be the set of all $n$-types over $A$ 
of the expansion of the theory $T$ by constants from $A$. 
There are two natural topologies on this set. 
The {\em logic topology} is defined by the basis consisting of 
sets of types of the form $[\phi (\bar{x})<\varepsilon ]$, 
i.e. types containing some $\phi (\bar{x})\le \varepsilon'$ with 
$\varepsilon '<\varepsilon$.   
The logic topology is compact. 

The $d$-topology is defined by the metric 
$$
d(p,q)= inf \{ max_{i\le n} d(c_i ,b_i )| \mbox{ there is a model } M \mbox{ with } M\models p(\bar{c})\wedge q(\bar{b})\}. 
$$ 
By Propositions 8.7 and 8.8 of \cite{BYBHU} the $d$-topology is finer 
than the logic topology and $(S_n (T_A ),d)$ is a complete space.

Definability in continuous structures is introduced as follows. 

\begin{definicja} 
Let $A\subseteq M$. 
A predicate $P:M^n \rightarrow [0,1]$ is definable in 
$M$ over $A$ if there is a sequence $(\phi_k (x) :k\ge 1 )$ 
of $L(A)$-formulas such that predicates interpreting 
$\phi_k (x)$ in $M$ converge to $P(x)$ uniformly in $M^n$. 
\end{definicja} 

A theory $T$ is {\em separably categorical} if any 
two separable models of $T$ are isomorphic. 
By Theorem 12.10 of \cite{BYBHU} a complete theory $T$ 
is separably categorical if and only if for each $n>0$, 
every $n$-type $p$ is principal. 
The latter means that for every model $M\models T$, the 
predicate $dist(x,p(M))$ is definable over $\emptyset$. 

Another property equivalent to separable categoricity states 
that for each $n>0$, the metric space $(S_n (T),d)$ is compact.  
In particular for every $n$ and every $\varepsilon$ there is 
a finite family of principal $n$-types $p_1 ,...,p_m$ so that 
their $\varepsilon$-neighbourhoods cover $S_n(T)$. 

In the classical first order logic a countable 
structure $M$ is $\omega$-categorical if and only 
if $Aut(M)$ is an {\em oligomorphic} permutation group, 
i.e. for every $n$, $Aut(M)$ has finitely many orbits 
on $M^n$. 
In continuous logic we have the following modification.  

\begin{definicja} 
An isometric action of a group $G$ on a metric space 
$({\bf X},d)$ is said to be approximately oligomorphic 
if for every $n\ge 1$ and $\varepsilon >0$ there is 
a finite set $F\subset {\bf X}^n$ such that 
$$ 
G\cdot F = \{ g\bar{x} : g\in G \mbox{ and } \bar{x}\in F\}
$$
is $\varepsilon$-dense in $({\bf X}^n,d)$. 
\end{definicja} 

Assuming that $G$ is the automorphism group of a non-compact 
separable continuous metric structure $M$, $G$ is approximately 
oligomorphic if and only if  the structure $M$ is separably 
categorical (C. Ward Henson, see Theorem 4.25 in \cite{scho}). 
It is also known that separably categorical structures are 
{\em approximately homogeneous} in the following sense: 
if $n$-tuples $\bar{a}$ and $\bar{c}$ have the same types 
(i.e. the same values $\phi (\bar{a})=\phi (\bar{b})$ for 
all $L$-formulas $\phi$) then for every $c_{n+1}$ and 
$\varepsilon >0$ there is an tuple $b_1 ,...,b_n ,b_{n+1}$ 
of the same type with $\bar{c},c_{n+1}$, so that 
$d(a_i, b_i )\le \varepsilon$ for $i\le n$.  
In fact for any $n$-tuples $\bar{a}$ and $\bar{b}$ 
there is an automorphism $\alpha$ of $M$ such that 
$$
d(\alpha (\bar{c}),\bar{a})\le d(tp(\bar{a}),tp(\bar{c})) +\varepsilon . 
$$  
(i.e $M$ is {\em strongly $\omega$-near-homogeneous} 
in the sense of Corollary 12.11 of \cite{BYBHU}). 

The following notion is helpful when we study some concrete 
examples, for example the Urysohn space. 
A relational continuous structure $M$ is 
{\em approximately ultrahomogeneous} if for any $n$-tuples 
$(a_1 ,..,a_n )$ and $(b_1 ,...,b_n )$ with the same 
quantifier-free type (i.e. with the same values of 
predicates for corresponding subtuples) and any 
$\varepsilon >0$ there exists $g\in Aut(M)$ such that 
$$ 
max \{ d(g(a_j ),b_j ): 1\le j \le n\} \le \varepsilon . 
$$  
As we already mentioned any Polish group can be chosen 
as the automorphism group of a continuous metric 
structure which is approximately ultrahomogeneous.  

The bounded Urysohn space $\mathbb{U}$ (see Section 3) 
is {\em ultrahomogeneous} in the traditional sense: 
any partial isomorphism between two tuples extends 
to an automorphism of the structure \cite{U}.  
Note that this obviously implies that $\mathbb{U}$ is 
approximately ultrahomogeneous. 

In Section 3 we will use the continuous version of 
$L_{\omega_1 \omega}$ from \cite{BYIov}. 
We remind the reader that continuous 
$L_{\omega_1 \omega}$-formulas are defined by the standard 
procedure applied to countable conjunctions and disjunctions 
(see \cite{BYIov}). 
Each continuous infinite formula depends on finitely many free variables. 
The main demand is the existence of continuity moduli of such 
formulas.  
It is usually assumed that a continuity modulus 
$\delta_{\phi ,x}$ satisfies the equality 
$$ 
\delta_{\phi ,x} (\varepsilon) = sup \{ \delta_{\phi ,x} (\varepsilon '): 0<\varepsilon' <\varepsilon \} 
$$ 
and 
$$ 
\delta_{\bigwedge \Phi ,x}(\varepsilon )= sup \{ \delta'_{\bigwedge \Phi ,x}(\varepsilon '): 0<\varepsilon' <\varepsilon \},  
\mbox{ where } \delta'_{\bigwedge \Phi ,x}= inf \{ \delta_{\phi ,x}: \phi\in \Phi \} .  
$$

\section{The space of metric structures.} 

In the first part of this section we introduce 
logic actions of isometry groups on spaces of continuous 
structures. 
In the second part we prove that these spaces are universal 
for Borel reducibility of orbits equivalence relations. 
We finish the section by a short discussion concerning possible 
developments of this result. 

\subsection{Logic action} 

We now fix a countable continuous signature 
$$
L=\{ d,R_1 ,...,R_k ,..., F_1 ,..., F_l ,...\}
$$ 
and a Polish space $({\bf Y},d)$. 
Let $S$ be a dense countable subset of ${\bf Y}$. 
Let $seq(S)=\{ \bar{s}_i :i\in \omega\}$ be the set (and an enumeration) 
of all finite sequences (tuples) from $S$. 
Let us define the space of metric $L$-structures on $({\bf Y},d)$. 
Using the recipe as in the case of $Iso({\bf Y})$ 
we introduce a metric on the set of $L$-structures as follows. 
Enumerate all tuples of the form $(\varepsilon ,j,\bar{s})$, where  
$\varepsilon \in \{ 0,1\}$ and when $\varepsilon =0$,  $\bar{s}$ is a tuple 
from $seq(S)$ of the length of the arity of $R_j$,  and for $\varepsilon =1$, 
$\bar{s}$ is a tuple from $seq(S)$ of the length of the arity of $F_j$.  
For metric $L$-structures $M$ and $N$ let   
$$ 
\delta_{seq(S)} (M ,N )= \sum_{i=1}^{\infty} \{ 2^{-i} |R^M_j (\bar{s} )-R^N_j (\bar{s} )| 
\footnote{resp. $ 2^{-i} d(F^M_j (\bar{s} ),F^N_j (\bar{s} ))$ when $\varepsilon =1$}  
\mbox{ : }  i \mbox{ is the number of } (\varepsilon ,j,\bar{s}) \} . 
$$ 
Since the predicates and functions are uniformly continuous 
(with respect to moduli of $L$) and $S$ is dense 
in ${\bf Y}$, we see that $\delta_{seq(S)}$ is a complete metric. 
Moreover by an appropriate choice of rational values for 
$R_j (\bar{s})$ we find a countable dense subset of metric 
structures on ${\bf Y}$, i.e. the space obtained is Polish. 
We denote it by ${\bf Y}_L$. 
It is clear that $Iso ({\bf Y})$ acts on ${\bf Y}_L$ continuously. 
Thus we consider ${\bf Y}_L$ as an $Iso({\bf Y})$-space and 
call it the {\em space of the logic action} on ${\bf Y}$. 
 
It is convenient to consider the following basis 
of the topology of ${\bf Y}_L$. 
Fix a finite sublanguage $L'\subset L$, a finite subset 
$S'\subset S$, a finite tuple $q_1 ,...,q_t \in {\bf Q} \cap [0,1]$ 
and a rational $\varepsilon \in [0,1]$ with $1-\varepsilon < 1/2$.  
Consider a diagram $D$ of $L'$ on $S'$ of some inequalities of the form   
$$
d(F_j (\bar{s}) ,s' ) > \varepsilon \mbox{ , } d(F_j (\bar{s}) ,s' )< 1 - \varepsilon , 
$$
$$ 
|R_j (\bar{s}) - q_i | > \varepsilon \mbox{ , } |R_j (\bar{s}) - q_i |< 1 - \varepsilon ,\mbox{ with } \bar{s}\in seq(S'), s'\in S' . 
$$
(i.e. in the case of relations we consider 
negations of statements of the form: 
$|R_j (\bar{s}) - q_i |\le \varepsilon$ , 
$|R_j (\bar{s}) - q_i |\ge 1 - \varepsilon$).  
The set of metric $L$-structures realizing $D$ is 
an open set of the  topology of ${\bf Y}_L$ and 
the family of sets of this form is a basis of this topology.   
Compactness theorem for continuous logic (see \cite{BYU}) 
shows that the topology is compact. 
We will call it logic too. 

The following proposition is very helpful. 

\begin{prop} \label{EsLo}
For any continuous formula $\phi(\bar{v},\bar{w})$ of the 
language $L$ there is a natural number $n$ such that 
for any tuple $\bar{a}\in {\bf Y}$ 
and $\varepsilon \in [0,1]$, the subset 
$$
Mod(\phi ,\bar{a},<\varepsilon )=\{ (M,\bar{c}):M \models\phi (\bar{a},\bar{c})<\varepsilon \} 
$$
$$
\mbox{ ( or } 
Mod(\phi ,\bar{a},>\varepsilon )=\{ (M,\bar{c}):M \models\phi (\bar{a},\bar{c})>\varepsilon \} \mbox{ ) } 
$$ 
of the corresponding space ${\bf Y}_{L\bar{c}}$ of 
$\bar{c}$-expansions of $L$-structures, belongs to ${\bf \Sigma}_n$. 
\end{prop}

{\em Proof.}
The proof is by induction on the complexity of $\phi$.
Assume that $\phi$ is atomic. 
A straightforward argument shows that if, for example, 
$\phi$ is $P(\bar{v},\bar{c})$, then for every $\varepsilon$ 
the set of $L\bar{c}$-structures 
satisfying $\phi (\bar{a},\bar{c})<\varepsilon$ is open. 

If $\phi$ is of the form $\neg\psi (\bar{v},\bar{w})$, 
then $Mod(\phi ,\bar{a},< \varepsilon )$ is the set 
$Mod( \psi ,\bar{a},>1-\varepsilon )$, 
i.e. of the next ${\bf \Sigma}_n$-class with respect to  
sets of the form $Mod (\psi, \bar{a}, <\varepsilon' )$.

If $\phi$ is of the form 
$max_{i\le k}\psi_{i}(\bar{v},\bar{w})$,
$k \in\omega$, then $Mod(\phi ,\bar{a},< \varepsilon )$ is the
intersection of the sets $Mod( \psi_i ,\bar{a},<\varepsilon )$, $i\le k$.
The cases of other Boolean connectives are similar. 
\parskip0pt 

When $\phi =inf_{u}\psi (\bar{v},u,\bar{w})$,
the corresponding subset is the (countable) union 
$$
\bigcup \{ Mod (\psi ,\bar{a},s_{k+1},<\varepsilon ) : s_{k+1} \in S \} , 
$$
i.e. of the same ${\bf \Sigma}_n$-class with $\psi$. 

When $\phi =sup_{u}\psi (\bar{v},u,\bar{w})$, the subset 
$Mod(\phi ,\bar{a},< \varepsilon )$ is the (countable) union of all 
possible intersections 
$$
\bigcap \{ Mod (\psi ,\bar{a},s_{k+1},<\varepsilon' ) : s_{k+1} \in S \} 
\mbox{ , } \varepsilon' <\varepsilon \mbox{ and } \varepsilon' \in \mathbb{Q} .
$$
$\Box$ 


\subsection{Reduction} 

Let $({\bf Y},d)$ be a Polish metric space with diameter not greater than 1 and   
$G< Iso({\bf Y},d)$ be a Polish group. 
The following theorem is the main result of Section 1. 

\begin{thm} \label{reduction} 
There is a continuous relational signature $L^*$ such that 
for any Polish $G$-space ${\bf X}$ there is a Borel 1-1-map 
${\cal M}: {\bf X} \rightarrow {\bf Y}_L$ such that 
for any $x,x'\in {\bf X}$ 
structures ${\cal M}(x)$ and ${\cal M}(x')$ are isomorphic 
if and only if $x$ and $x'$ are in the same $G$-orbit. 
\end{thm} 

In other words the map ${\cal M}$ is a Borel $G$-invariant 1-1-reduction 
of the $G$-orbit equivalence relation on ${\bf X}$ to the 
$Iso ({\bf Y})$-orbit equivalence relation on the space ${\bf Y}_{L^*}$. 
Let us start with some preliminaries. 
Let $\langle {\bf X}, \tau ,d^{\tau} \rangle$
be a Polish $G$-space with a basis ${\mathcal{A}}=\{ A_l :l\in \omega\}$. 
To describe a reduction of the $G$-space ${\bf X}$ to 
an $Iso({\bf Y})$-space of countinuous structures on $({\bf Y},d)$  
we use \cite{melleray} and some standard  
ideas already applied for closed subgroups of $S_{\infty}$ 
(see Section 6.1 of \cite{hjorth}). 

We will assume that $d$ and $d^{\tau}$ have values from $[0,1]$ 
(if necessary we may replace them by $\frac{d(x,y)}{1+ d(x,y)}$). 
We fix some countable dense set $S\subset {\bf Y}$ and 
enumerate $S=\{ s_1 ,s_2 ,... \}$ and all orbits of $G$ of 
finite tuples of $S$ (i.e. of $Seq(S )$). 
For the closure of such an $n$-orbit $C$ define a predicate 
$R_{\overline{C}}$ on $({\bf Y},d)$ by 
$$ 
R_{\overline{C}}(y_1 ,...,y_n ) = d((y_1 ,...,y_n ),\overline{C}) 
\mbox{ ( i.e. } inf \{ d(\bar{y}, \bar{c}): \bar{c} \in \overline{C}\} ).   
$$ 
It is proved in Theorem 6 of \cite{melleray} that the 
continuous structure $M$ of all these predicates on ${\bf Y}$ 
is approximately ultrahomogeneous and $G$ is its automorphism group. 

Let $L$ be the language of $M$. 
For every pair of natural numbers $k>0$ and $l$ we add to $L$ 
a predicate $R_{k,l}$ of arity $k$. 
The extended language will be denoted by $L^*$. 
Then to every $x\in {\bf X}$ we assign an $L^*$-expansion of $M$ 
where the predicates $R_{k,l}(y_1 ,...,y_k )$ are interpreted as follows: 
$$
inf \{ max (d((h(y_1 ),...,h(y_k )),(s_1 ,...,s_k )), d^{\tau}(h x ,x')): 
x'\in A_l \mbox{ and } h \in G\} . 
$$
It is easily seen, that these predicates are uniformly 
continuous with respect to the continuous modulus $id$ 
for each variable.  
Let $M(x)$ denote this expansion.  
Let ${\cal M}$ denote the map $x\rightarrow M(x)$. 
The following proposition implies the theorem above. 

\begin{prop} \label{1} 
The map ${\cal M}$ is a Borel $G$-invariant 1-1-reduction 
of the $G$-orbit equivalence relation on ${\bf X}$ to the 
$Iso ({\bf Y})$-orbit equivalence relation on the space ${\bf Y}_{L^*}$ 
of all $L^*$-stuctures. 
Moreover the ${\cal M}$-preimage of any open subset of 
${\bf Y}_{L^*}$ belongs to ${\bf \Delta}_3$ and the ${\cal M}$-image 
of any open subset of ${\bf X}$ belongs to $F_{\sigma}$. 
\end{prop} 

{\em Proof.} 
To see $G$-invariantness note that the condition $gx =x'$ 
with $g\in G$, implies the property that for every $k,l\in \omega$ 
and $y_1 ,...,y_k \in {\bf Y}$
$$ 
R^{M(x')}_{k,l}(y_1 ,...,y_k ) =R^{M(x)}_{k,l}(g^{-1}(y_1 ),...,g^{-1}(y_k ))
$$
(i.e. $g$ maps $M(x)$ to $M(x')$). 
This follows from the fact that for any $h \in G$ 
$$ 
max (d((h(y_1 ),...,h(y_k )), (s_1 ,...,s_k)), d^{\tau}(h x'  ,A_l )) = 
$$
$$
max (d((hg(g^{-1}(y_1 )),...,hg(g^{-1}(y_k )),(s_1 ,...,s_k)), 
d^{\tau} (hg x ,A_l )).  
$$  

Let us check that the map $x\rightarrow M(x)$ is injective. 
Assume $x,x'\in {\bf X}$ and 
$x\not= x'$. 
Then there are basic open sets $A_l$ and $A_m$ such that 
$d(A_l ,A_m ) >0$, $A_l \cap \{ x,x' \} = \{ x\}$ 
and $A_m \cap \{ x,x' \} = \{ x' \}$. 
Since the $G$-action is continuous, there is an open set 
$V\subset G$ containing the identity  
such that $V x\subseteq A_l$ and $V x' \subseteq A_m$. 
We may think that $V$ consists of all $h\in G$ such that 
for some sufficiently small $\varepsilon$ and a natural $k$ 
$$
\sum_{i=1}^{k} 2^{-i} min(1, d(h(c_i ),c_i )) < \varepsilon .
$$ 
This obviously means that $R^{M(x)}_{k,l} (c_1 ,...,c_{k})$ 
(which is $0$) cannot be equal to  
$R^{M(x')}_{k,l}(c_1 ,...,c_k )$, i.e. $M(x)\not= M(x')$.  

We can now see that if the structures $M(x)$ 
and $M(x')$ are isomorphic then $x$ and $x'$ belong to 
the same $G$-orbit. 
Indeed, by the choice of relations in $M(x)$ such an 
isomorphism can be realized by an element (say $g$) of $G$.  
Then we see that $M(x')=M(gx)$ and thus  by the definition 
of relations $R_{k,l}$ (in particular for tuples of the form 
$c_1 ,...,c_k$), $x'=gx$. 

Let us prove the last statement of the proposition.  
Fix a finite sublanguage $L'\subset L^*$, a finite subset 
$S'\subset S$, a finite tuple $q_1 ,...,q_t \in {\bf Q} \cap [0,1]$ 
and a small rational $\varepsilon \in [0,1]$.  
Consider a diagram $D$ of $L'$ on $S'$ consisting of 
atomic statements and negations (in the standard sense)  
of atomic statements of the following form: 
$$ 
|R (\bar{s}) - q_i | \le \varepsilon \mbox{ , }  |R (\bar{s}) - q_i | \ge \varepsilon . 
$$ 
The set of metric $L^*$-structures realizing $D$ is an $F_{\sigma}$-set 
of the logic topology.  
Since each $M(x)$ belongs to the closed subset of all 
expansions of $M$ we will assume that $D$ is consistent 
with the elementary diagram of $M$. 
As a result $D$ is determined by formulas of the form 
( with or without $\neg$ in its standard meaning) 
$$ 
| R_{k,l}(s_{i_1} ,...,s_{i_k} ) - q_i| \le \varepsilon \mbox{ , } 
| R_{k,l}(s_{i_1} ,...,s_{i_k} ) - q_i| \ge \varepsilon \mbox{ , } s_{i_j} \in C' . 
$$ 
Thus the set of $x$ with $M(x)\models D$ is a finite intersection 
of sets of the form 
$$ 
\{ x \mbox{ : } inf \{ max (d((h(s_{i_1} ),...,h(s_{i_k} )), 
(s_1 ,...,s_k)), d^{\tau}(h x ,x')): x'\in A_l \mbox{ and } h \in G\}| <\varepsilon \} , 
$$ 
$$ 
\{ x \mbox{ : } inf \{ max (d((h(s_{i_1} ),...,h(s_{i_k} )), 
(s_1 ,...,s_k)), d^{\tau}(h x ,x')): x'\in A_l \mbox{ and } h \in G\}| \le \varepsilon \}  
$$
(or their complements), where the first one is open 
and thus the second one is an $F_{\sigma}$-set. 
From this we conclude that the preimage of 
an open set belongs to ${\bf \Delta}_3$. 

Now consider the case of the ${\mathcal{M}}$-image 
of the basic open set $A_l$. 
Note that for any point $a$ of $A_l$ there is some $A_m$ 
with the closure satisfying $a\in \overline{A}_m \subseteq A_l$.  
On the other hand it is easy to see that the ${\mathcal{M}}$-image 
of the closure $\overline{A}_m$ consists of all $M(x)$ such that 
for every $k$, $R_{k,m}(s_1 ,...,s_k ) = 0$. 
It is clear that this is a closed set. 
This means that the ${\cal M}$-image of $A_l$ is an $F_{\sigma}$-set. 
$\Box$ 

\bigskip 

The proof of Theorem 2.2 of the paper \cite{CL} 
contains a construction which is equivalent to 
the proof  of Theorem \ref{reduction} and is based 
on Section 2.6 of \cite{bk}. 
Theorem 2.2 of \cite{CL} complements it by 
applications of the universality of the isometry group 
of the Urysohn sphere $\mathbb{U}$ and 
a version of the L\'{o}pez-Escobar theorem. 
As a result the space ${\bf Y}$ can be replaced 
by $\mathbb{U}$ and the image of $\mathcal{M}$ 
can be taken to be $L_{\omega_1 \omega}$-axiomatisable.  

These results are slightly connected with the conjecture of G.Hjorth, 
that for any Borel equivalence relation $E$ Borel reducible to 
an orbit equivalence relation of a Polish $G$-space there is 
a Polish $G$-space ${\bf X}$ such that the orbit equivalence 
relation $E^{{\bf X}}_{G}$ is Borel and 
$E\le_{B} E^{{\bf X}}_{G}$. 
Since any Polish $G$ can be realised as an isometry group, 
we may always assume that $E$ is Borel reducible to 
the orbit equivalence relation of an {\em logic} action. 
We hope that this observation may bring some 
model-theoretic tools for the Hjorth's conjecture.

\section{Grey subsets and nice bases} 

In this section we introduce the main notions 
of the paper and develope the basic theory of 
nice topologies. 
In particular in the second part of the section 
we prove an existence theorem for nice topologies. 
The first part developes necessary techniques. 

\subsection{Grey subsets} 

Proposition \ref{EsLo} naturally fits to the notion 
of grey subsets, introduced in \cite{BYM}. 
Let us recall this notion. 

A function $\phi$ from a space ${\bf X}$ to $[-\infty ,+\infty ]$ 
is {\em upper (lower) semi-continuous} if the set 
$\phi_{<r}$ (resp. $\phi_{>r}$) is open for all $r\in \mathbb{R}$ 
(here $\phi_{<r} = \{ z\in {\bf X}: \phi (z) <r \}$, a {\em cone}). 
A {\em grey subset} of ${\bf X}$, denoted 
$\phi \sqsubseteq {\bf X}$, is a function 
${\bf X}\rightarrow [0,\infty ]$. 
It is {\em open (closed)}, $\phi \sqsubseteq_o {\bf X}$ 
(resp. $\phi \sqsubseteq_c {\bf X}$), 
if it is upper (lower) semi-continuous. 
We also write $\phi \in {\bf \Sigma}_1$ when 
$\phi \sqsubseteq_o {\bf X}$ and we write 
$\phi \in {\bf \Pi}_1$ when $\phi \sqsubseteq_c {\bf X}$.  
We will assume below that values of 
a grey subset belong to $[0,1]$. 

It is observed in \cite{BYM} if $\Phi$ is a family 
of upper semi-continuous functions, then 
$inf \Phi : x\rightarrow inf \{ f(x):f\in \Phi\}$ 
is upper semi-continuous as well. 
If additionally ${\bf X}$ admits a countable 
base then there exists a countable subfamily 
$\Phi_0 \subseteq \Phi$ such that 
$inf \Phi_0 = inf \Phi$.

When $G$ is a Polish group, then a grey subset 
$H \sqsubseteq G$ is called a {\em grey subgroup} if 
$$
H(1)=0 \mbox{ , }\forall g\in G (H(g)=H(g^{-1})) \mbox{ and } 
\forall g, g'\in G (H(gg')\le H(g)+H(g')). 
$$ 
This is equivalent to Definition 2.5 from \cite{BYM}. 
It is worth noting that by Lemma 2.6 of \cite{BYM} 
an open grey subgroup is clopen. 

By induction on $\alpha <\omega_1$ we define Borel classes 
${\bf \Sigma}_{\alpha}$, ${\bf \Pi}_{\alpha}$ with 
${\bf \Pi}_{\alpha} = \{ \neg \phi :\phi \in {\bf \Sigma}_{\alpha} \}$. 
We will say that a grey subset $\phi$ is ${\bf \Sigma}_{\alpha}$ 
if $\phi = inf \Phi$ for some countable family of 
grey subsets from $\bigcup \{ {\bf \Pi}_{\gamma}: \gamma <\alpha \}$. 
Note that in this case $\phi_{<r} \in 
{\bf \Sigma}_{\alpha}^{{\bf X}}$ for all $r\in [0,1]$. 
On the other hand if $\phi \in {\bf \Pi}_{\alpha}$, then 
for any $r>0$, $\phi_{\le r}\in {\bf \Pi}_{\alpha}^{{\bf X}}$.  
In this paper we consider merely Borel grey subsets.

By standard inductive argument one can prove that 
for any Borel $A\subseteq X$, its characteristic 
function  ${\bf O}_A$ defined by ${\bf O}_A(x)=
\begin{cases}
0 &\mbox{ if } x\in A ;\\
1 &\mbox{ if } x\not\in A
\end{cases}$ 
is a Borel grey subset. 
Moreover ${\bf O}_A$ is a 
${\bf \Sigma}_{\alpha}$(${\bf \Pi}_{\alpha}$)-grey 
subset if and only if $A\in {\bf \Sigma}_{\alpha}^{{\bf X}}$ 
(resp. $A\in {\bf \Pi}_{\alpha}^{\bf X}$).
Sometimes we will identify subsets with their 
characteristic grey subsets.
In particular $G$ may stand for both $G$ and ${\bf O}_G$.

It is clear that for every continuous structure $M$ 
(defined on ${\bf Y}$) any continuous formula $\phi (\bar{x})$ 
defines a clopen grey subset of $M^{|\bar{x}|}$. 
Moreover note that when $\phi (\bar{x},\bar{c})$ is a continuous 
formula with parameters $\bar{c}\in M$ and $\delta$ is a linear inverse 
continuous modulus for $\phi (\bar{x}, \bar{y})$ (see Definition \ref{inverse}), 
then $\phi$ is invariant with respect to the open grey subgroup 
$H_{\delta, \bar{c}} \sqsubseteq Aut(M)$ defined by  
$$
H_{\delta, \bar{c}}(g) =\delta ( max (d(c_1 ,g(c_1 )),...,d(c_n ,g(c_n )))) \mbox{, where } g\in Aut(M) ,
$$
in the following sense: 
$$
\phi(g(\bar{a}),\bar{c}) \le \phi (\bar{a},\bar{c}) + H_{\delta ,\bar{c}}(g) .    
$$ 
The fact that $H_{\delta ,\bar{c}}$ is a grey subgroup follows 
from the condition that the action is isometric and $\delta$ is linear. 
The invariantness is a consequence of  
$$
\phi (g(\bar{a}),g(\bar{c})) = \phi (\bar{a},\bar{c}) 
$$ 
together with uniform $\delta$-continuity of  $\phi (g(\bar{x}),g(\bar{y}))$.   
We arrive at the following definition. 
\begin{definicja} 
Let ${\bf X}$ be a continuous $G$-space.  
A grey subset $\phi \sqsubseteq {\bf X}$ is called  invariant with 
respect to a grey subgroup $H \sqsubseteq G$ if 
for any $g\in G$ and $x\in {\bf X}$ we have $\phi (g(x)) \le \phi (x) \dot+ H(g)$. 
\end{definicja}  

Since $H(g)=H(g^{-1})$, the inequality from the definition is 
equivalent to $\phi (x) \le \phi (g(x)) \dot+ H(g)$. 
 
In fact Proposition \ref{EsLo} says that any continuous sentence  
$\phi (\bar{c})$ defines a grey subset of ${\bf Y}_L$ which belongs 
to ${\bf \Sigma}_n$  for some $n$: 
$$
\phi (\bar{c}) \mbox{ takes } M \mbox{ to the value } \phi^{M} (\bar{c}). 
$$ 
We extend this as follows. 

\begin{lem} 
Let $\delta$ be an inverse continuity modulus for $\phi (\bar{x})$, 
which is linear. 
The grey subset defined by $\phi (\bar{c}) \sqsubseteq {\bf Y}_L$ 
is invariant with respect to the {\bf grey stabiliser}  
$H_{\delta ,\bar{c}}\sqsubseteq Iso({\bf Y})$ defined 
as above, i.e. 
$$
H_{\delta ,\bar{c}}(g) =\delta ( max (d(c_1 ,g(c_1 )),...,d(c_n ,g(c_n )))) \mbox{, where } g\in Iso({\bf Y}) . 
$$
\end{lem} 

{\em Proof.} 
Since $\phi^{g(M)}(g(\bar{c})) = \phi^{M}(\bar{c})$ 
and $\phi^{g(M)}(\bar{x})$ is a $\delta$-uniformly continuous 
map ${\bf Y}^{|\bar{c}|} \rightarrow [0,1]$, we have 
$$
\phi^{g(M)}(\bar{c}) \le \phi^{M}(\bar{c}) + H_{\delta ,\bar{c}}(g) . 
$$ 
$\Box$   

\bigskip

A grey subset $\phi \sqsubseteq {\bf X}$ is called {\em meagre} 
in an open $\psi \sqsubseteq {\bf X}$ if there is $r>0$ such that 
$\phi_{<s}$ is meagre in $\psi_{<s}$ for all $s\le r$ 
(this extends Definition 1.4(iii) of \cite{BYM}). 

We now define two Vaught transforms as follows. 
Let $G$ continuously act on ${\bf X}$. 
For any non-empty open $J\sqsubseteq G$ let 
$$
\phi^{\Delta J}(x) = inf \{r\dot+ s \mbox{ : }\{ h: \phi (h(x))<r \} 
\mbox{ is not meagre in } J_{<s} \} , 
$$  
$$ 
\phi^{*J}(x) = sup \{ r\dot- s\mbox{ : }\{ h: \phi (h(x))\le r \} 
\mbox{ is not comeagre in } J_{<s} \} .  
$$ 
Using the original topological Vaugt transforms 
one can rewrite the above definitions as follows:
$$
\phi^{\Delta J}(x) = inf \{r\dot+ s \mbox{ : } x\in (\phi_{<r})^{\Delta(J_{<s})} \} , 
$$  
$$ 
\phi^{*J}(x) = sup \{ r\dot- s\mbox{ : } x\not\in (\phi_{\le r})^{*(J_{<s})} \}. 
$$
Note that since $J$ is an open grey subset, if $\{ h: \phi (h(x))<r \}$ is 
not meagre in  $J_{<s}$, then it is not meagre in any $J_{<s'}$ with $s<s'$.

The lemma below expresses some basic relationships 
betweeen usual topological Vaught transforms
and their grey versions introduced above. 

\begin{lem} \label{grey}  
Let $A\subseteq {\bf X}$ be a Borel set  and $\phi\sqsubseteq {\bf  X}$ 
be a Borel  grey set.
Let $u\subseteq G$  be an open set and $J\sqsubseteq G$ be an open  grey set.
Then for any  $r\in (0,1]\cap{\mathbb{Q}}$ we have the following equalities: \\ 
(1) $(\phi_{<r})^{\Delta u}=(\phi^{\Delta {\bf O}_u})_{< r}$ and 
$(\phi_{\le r})^{*u}=(\phi^{*{\bf O}_u})_{\le  r}$; \\ 
(2) $A^{\Delta (J_{<r})}=({\bf O}_A^{\Delta J})_{<r}$
and $A^{* (J_{\le r})}=({\bf O}_A^{* J})_{\le r}$.
\end{lem}

{\em Proof}. (1)  
For any  $s\in (0,1]$ we have  $({\bf O}_u )_{<s} =u$ and so 
$(\phi_{<r})^{\Delta u}=(\phi_{<r})^{\Delta ({\bf O}_u )_{<s}}$ 
and $(\phi_{\le r})^{*u}=(\phi_{\le r})^{*({\bf O}_u )_{<s}}$.
Using this  we argue as follows: 
$$ 
\begin{array}{c}
x\in (\phi_{<r})^{\Delta u}\\
\Updownarrow\\
x\in (\bigcup\{\phi_{<r'}: r'<r\} )^{\Delta u}\\
\Updownarrow\\
(\exists r'<r)(x\in (\phi_{<r'})^{\Delta u})\\
\Updownarrow\\
(\exists r'<r)(\forall s>0)(x\in (\phi_{<r'})^{\Delta ({\bf O}_u )_{<s}})\\
\Updownarrow\\
(\exists r'<r)(inf\{r''\dot +s: x\in 
(\phi_{<r''})^{\Delta({\bf O}_u )_{<s}}\}\ \le  r')\\ 
\Updownarrow\\
x\in (\phi^{\Delta {\bf O}_u})_{< r} . \end{array} 
$$ 
Since the set of rational numbers is dense in $[0,1]$ 
we may assume that $\bigcup$ and $\bigcap$ above and below 
are applied to countable families. 
The next  equality can be derived in a similar way.

$$ 
\begin{array}{c}
x\in (\phi_{\le r})^{* u}\\
\Updownarrow\\
x\in (\bigcap\{\phi_{\le r'}: r'> r\} )^{* u}\\
\Updownarrow\\
(\forall r'> r)(x\in (\phi_{\le  r'})^{* u})\\
\Updownarrow\\
(\forall r'> r)(\forall s>0)(x\in (\phi_{\le r'})^{* ({\bf O}_u )_{<s}})\\
\Updownarrow\\
(\forall r'> r)(sup\{r''\dot- s: x\not \in 
(\phi_{<r''})^{*({\bf O}_u )_{<s}}\}\ <  r')\\ 
\Updownarrow\\
x\in (\phi^{* {\bf O}_u })_{\le r} . \end{array} 
$$ 

(2)  We argue as follows: 
$$\begin{array}{c}
x\in A^{\Delta (J_{<r})}\\
\Updownarrow\\
x\in \bigcup\{A^{\Delta (J_{<r'})}: r'<r\}\\
\Updownarrow\\
(\exists r'<r)(x\in A^{\Delta (J_{<r'})})\\
\Updownarrow\\
(\exists r'<r)(\forall s>0)(x\in (({\bf O}_A)_{<s})^{\Delta (J_{<r'})}\\
\Updownarrow\\
(\exists r'<r)(inf\{r''\dot +s: x\in 
({\bf O}_A)_{<s}^{\Delta(J_{<r''})}\}\ \le  r')\\ 
\Updownarrow\\
x\in ({\bf O}_A^{\Delta J})_{< r} . \end{array} 
$$ 

Again, the next  equality can be proved in a similar way.
$\Box$ 

\bigskip

\begin{lem}\label{prop} 
Let $J\sqsubseteq G$ be an open grey subset. Then: \\
(1) $\phi^{*J} = 1-(1- \phi)^{\Delta J}$, i.e.
$\phi^{*J}(x) = 1-(1- \phi)^{\Delta J}(x)$ for all $x\in {\bf X}$. \\  
(2) $\phi^{\Delta J} \le  \phi^{* J}$, i.e.
$\phi^{\Delta J}(x) \le \phi^{* J}(x)$ for all $x\in {\bf X}$. \\
(3) If $\phi$ is a grey ${\bf \Sigma}_{\alpha}$-subset, 
then $\phi^{\Delta J} $ is also ${\bf \Sigma}_{\alpha}$. 
If $\phi $ is a grey ${\bf \Pi}_{\alpha}$-subset, 
then $\phi^{*J} (x )$ is also ${\bf \Pi}_{\alpha}$. \\ 
(4) Vaught transforms of Borel grey subsets are Borel. 
\end{lem} 

{\em Proof.} First observe that for any $r,s\in [0,1]$ we have
$1-(r\dot-s)=(1-r)\dot+s$.
Now statement (1) is straightforward:  
$$\begin{array}{c}
1-(\phi)^{*J}(x) =
1- sup \{r\dot- s \mbox{ : }x\not\in (\phi_{\le r})^{*(J_{<s})} \} =  
\\
inf \{ (1-r)\dot+ s\mbox{ : }x\in {\bf X}\setminus (\phi_{\le r})^{*(J_{<s})} \} = 
inf \{(1-r)\dot+ s\mbox{ : }  x\in (\phi_{> r})^{\Delta(J_{<s})} \} = 
\\
inf \{(1-r)\dot+ s\mbox{ : } x\in ((1-\phi)_{< 1-r})^{\Delta(J_{<s})} \}=
inf \{r'\dot+ s\mbox{ : } x\in ((1-\phi)_{< r'})^{\Delta(J_{<s})} \}=\\
(1- \phi)^{\Delta J}(x).
\end{array}$$ 

To prove (2) take any $x\in {\bf X}$ and $r\in [0,1)$ 
such that $\phi^{\Delta J}(x)>r$.
Then for some $\varepsilon > 0$ we have 
$\phi^{\Delta J}(x)>r\dot+ 4\varepsilon$, and so 
$x\not\in (\phi_{\le r\dot+ 3\varepsilon})^{\Delta(J_{< \varepsilon} )}$.
Then by the properties of the original 
topological Vaught transforms we see that
$x\not\in (\phi_{\le r\dot+ 3\varepsilon})^{* (J_{< \varepsilon } )}$.
Hence $\phi^{* J}(x)\ge (r\dot+ 3\varepsilon )\dot- \varepsilon$, 
i.e. $\phi^{* J}(x)> r$.

Statement (3) follows by induction starting 
with the case when $\phi$ is open. 
In this case for any $x \in {\bf X}$ and any open $U\subseteq G$ 
the set $\{ h\in U : \phi (h(x))<r' \}$ is open too. 
Thus the set 
$$
\{ x \mbox{ : } \{ h\in U : \phi (h(x))<r' \} \mbox{ is not meagre in } J_{<s} \} 
$$
coincides with 
$$
\{ x \mbox{ : } \{ h\in U : \phi (h(x))<r' \} \cap J_{<s} \not= \emptyset \} . 
$$
Since the action of $G$ is continuous and both 
$\phi$ and $J$ are open, the latter set is open too. 
Now note that $(\phi^{\Delta J})_{<r}$ is a union of 
sets of this form (taking $U=G$ above).  

When $\phi = inf \Psi$ where $\Psi$ is a countable subfamily of 
$\bigcup_{\gamma<\alpha} {\bf \Pi}_{\gamma}$, then 
$$
\phi^{\Delta J} (x) = inf \{ r\dot+ s \mbox{ : }\{ h: inf \Psi (h(x))<r \} \mbox{ is not meagre in } J_{<s} \} . 
$$
Since $\{ h: inf \Psi (h(x))<r \}$ is not meagre 
in $J_{<s}$ if and only if one of $\{ h: \psi (h(x))<r \}$, 
$\psi \in \Psi$, is not meagre in $J_{<s}$, we see that 
$\phi^{\Delta J} = inf \{ \psi^{\Delta J} : \psi \in \Psi \}$, 
i.e. belongs to ${\bf \Sigma}_{\alpha}$.  

The case of $\phi^{*J}$ with $\phi\in {\bf \Pi}_{\alpha}$ 
follows from the $\Delta J$-case and statement (1). 

Statement (4) follows from statement (3).  $\Box$ 

\bigskip

The argument of statement (3) is also 
applied in the following lemma.

\begin{lem} \label{inv}
Let $H$ be  an open  grey subgroup of $G$.  \\ 
(1) If $\phi$ is an open grey subset then  $\phi^{\Delta H} \le \phi$; \\ 
(2) For any grey subset $\phi$ both $\phi^{*H}(x)$ 
and $\phi^{\Delta H}(x)$ are $H$-invariant:  
$$
| \phi^{*H}(h(x)) - \phi^{*H}(x)|\le H(h) \mbox{ and }$$
$$ 
| \phi^{\Delta H}(h(x)) - \phi^{\Delta H}(x)|\le H(h) . 
$$ 
Moreover if  $\phi$ is $H$-invariant, then 
$$ 
\phi^{*H}(x) = \phi (x) = \phi^{\Delta H}(x). 
$$   
\end{lem} 

{\em Proof.} 
(1) Let $\phi (x)=r$. 
Then for any $\varepsilon >0$ the set 
$\{ h: \phi (h(x))<r +\varepsilon  \}$ 
is open and intersects open $H_{<\varepsilon}$ 
(they contain the neutral element).  
Thus $\{ h: \phi (h(x))<r+\varepsilon \}$ 
is not meagre in $H_{<\varepsilon}$ and 
$\phi^{\Delta H}(x) \le r\dot+ 2\varepsilon$. 

(2)  Let  $H(h) =t$.
Since $H$ is a grey subgroup, $H(h^{-1})=t$ and 
$H(gh), H(gh^{-1})\le H(g)\dot +t$, for all $g\in G$.
Note that for any $r,s\in (0,1]$ we have 
$$ 
\{ g: \phi (g(h(x)))\le r \}=\{ g: \phi (g(x))\le r \}h^{-1} 
\mbox{ and } H_{<s}h^{-1}\subseteq H_{<s\dot+ t}.  
$$ 
Thus if the set $\{ g: \phi (g(x))\le r \}$ 
is not comeagre in $H_{<s}$, 
then the set $\{ g: \phi (g(h(x)))\le r \}$ 
is not comeagre in $H_{<s\dot+t}$.
Hence 
$$
\{r\dot- (s\dot+ t)\ :\ x\not\in (\phi_{\le r})^{*(H_{<s})}\}
\subseteq  \{r\dot- s'\ :\ h(x)\not\in (\phi_{\le r})^{*(H_{<s'})}\} .
$$ 
Since for every $u,w,v\in [0,1]$ we have 
$u\dot-(w\dot+ v)=(u\dot- w)\dot- v$,
then the latter implies  
$\phi^{*H}(x)\dot- t\le \phi^{*H}(h(x))$.
Replacing $x$ by $h(x)$ and $h$ by $h^{-1}$, 
we obtain $\phi^{*H}(h(x)) \le \phi^{*H}(x) \dot+ t$.

A similar argument works for $\phi^{\Delta H}$. 

To see the last statement let $\phi (x) =r$. 
Take any $\varepsilon >0$ and consider the set 
$\{ h: r\dot -\varepsilon < \phi (h(x))<r\dot +\varepsilon \}$. 
Since this set is comeagre in (coincides with) the open $H_{<\varepsilon}$, 
we see that
$$ 
|\phi^{*H}(x) - \phi (x))| <\varepsilon  \mbox{ and } |\phi^{\Delta H}(x) - \phi (x))| < \varepsilon .  
$$ 
$\Box$

\bigskip

If $H$ is a grey subgroup, then for every $g\in G$ we define 
the grey coset $Hg$ and the grey conjugate $H^g$ as follows:
$$
\begin{array}{l@{\ = \ }l}
Hg(h)&H(hg^{-1})\\
H^g(h)&H(ghg^{-1}).
\end{array}
$$
Observe that if $H$ is open, then $Hg$ is an open grey subset and
$H^g$ is an open grey subgroup.

\begin{lem}\label{H^g} 
Let $H$ be an open grey subgroup, $g\in G$ and $\rho=Hg$.
Let $\phi\sqsubseteq \bf{X}$ be a grey subset. 
Then both $\phi^{\Delta\rho}$ and $\phi^{*\rho}$ are $H^g$-invariant.
\end{lem}

{\em Proof.} Let  $H^g(h)=t$.
For every $f\in G$ we have   
$$ 
\rho(fh^{-1})=H(fh^{-1}g^{-1})=H((fg^{-1})(gh^{-1}g^{-1}))\le \rho(f)+H^g(h^{-1}) 
\mbox{ , i.e. } \rho(fh^{-1} )\le \rho(f)+t . 
$$ 
Hence for every $s\in ( 0,1]$ we have 
$\rho_{<s}h^{-1}\subseteq \rho_{<s\dot+t}$. 
Using this and   
$$ 
\{ f: \phi (f(h(x)))< r \}=\{ f: \phi (f(x))< r \}h^{-1} , \mbox{ where }r\in (0,1] ,
$$
we see that if the set $\{ f: \phi (f(x))< r \}$ is not meagre in $\rho_{<s}$, 
then the set $\{ f: \phi (fh(x)))< r \}$ is not meagre in $\rho_{<s\dot+t}$, 
 i.e. if $x\in (\phi_{<r})^{\Delta(\rho_{<s})}$ then
$h(x)\in (\phi_{<r})^{\Delta(\rho_{<s\dot+t})}$.
Therefore $\{r\dot+s'\ : h(x)\in (\phi_{< r})^{\Delta(\rho_{<s'})}\}
\supseteq  \{r\dot+(s\dot+t) : x\in (\phi_{< r})^{\Delta(\rho_{<s})}\}$.
This implies $\phi^{\Delta\rho}(hx)\le \phi^{\Delta\rho}(x)\dot+H^g(h)$.

In a similar way we proceed with $\phi^{*\rho}$.
$\Box$ 

\bigskip 

The following lemma 
will be used in Section 4. 
In the proof we develope the arguments 
from Lemmas \ref{prop}(3) and \ref{inv}(1). 

\begin{lem} \label{greyold} 
Let $H$ be an open grey subgroup and  
$\phi$ and $\psi$ be open grey subsets.  \\ 
If for some $\varepsilon >0$, 
$\phi_{<r}$ is contained in the closure 
of $H_{<\varepsilon}\psi_{<t}$ then 
$(\phi^{\Delta H})_{< r}$ 
is contained in the closure of 
$H_{< r} (\psi^{\Delta H})_{< t+\varepsilon})$.  
\end{lem} 

{\em Proof.}  
Take any $x\in( \phi^{\Delta H})_{<r}$.  
Find $\sigma <r$ and  
$g_1 \in \{ h: \phi (h(x))<r-\sigma \} \cap H_{<\sigma}$. 
We can choose $g_2 \in H_{<\varepsilon}$ 
and $y$ with $\psi (y) = t'<t$ and 
$g^{-1}_2 (y)$ sufficiently close to $g_1 (x)$.  
Thus for some (any) $\varepsilon' < \frac{t-t'}{2}$, 
the open set $H_{<\varepsilon}$ 
has a non-empty intersection with 
$\{ h: \psi (h g^{-1}_2 (y))<t'+\varepsilon'\}$, 
i.e 
$g^{-1}_2 (y) \in (\psi^{\Delta H})_{< t+\varepsilon}$   
and 
$g^{-1}_1 g^{-1}_2 (y)\in H_{<r}(\psi^{\Delta H})_{< t+\varepsilon}$. 
In particular $x$ belongs to the closure of 
$H_{< r}(\psi^{\Delta H})_{< t+\varepsilon}$. 
$\Box$ 

\bigskip

The following statement follows by the same proof (even with some simplifications): 
\begin{quote} 
{\em 
 if under circumstances above, $\phi_{<r}$ is contained in 
the closure  $\overline{\psi_{<t}}$ then $(\phi^{\Delta H})_{<r}$ is 
contained in $\overline{H_{<r} (\psi^{\Delta H})_{< t}}$.} 
\end{quote}

\subsection{Nice bases}

We now consider $G$ together with a distinguished 
countable family of clopen grey subsets $\mathcal{R}$. 
We will assume that a countable family of sets of the form 
$\rho_{<r}$ for $\rho \in \mathcal{R}$ and real $r$, 
forms a basis of the topology of $G$. 
In fact we usually assume it for 
$\{ \rho_{<q} : \rho\in \mathcal{R}$ and $q\in \mathbb{Q}^{+}\}$. 

We also assume that $\mathcal{R}$ consists of grey cosets, 
i.e. for such $\rho\in\mathcal{R}$ there is a grey subgroup 
$H\in \mathcal{R}$ and an element $g_0 \in G$ so that for any 
$g\in G$, $\rho (g ) = H(g g^{-1}_0 )$.

\begin{remark}\label{Gbasis}  
For every Polish group $G$ there is a countable family of open 
grey subsets $\mathcal{R}$ as above. 

{\em Indeed, fix an arbitrary  compatible right-invariant metric 
$d\le 1$ on $G$. 
Put $H(g) = d(1,g)$. 
It is easily seen that $H(x)$ is 
of an open grey subgroup such that 
$\{H_{<r}: r\in [0,1]\cap {\mathbb{Q}}^+\}$ is a basis of open 
neighbourhoods of the unity $1$.
Now take a dense countable subgroup $G_0 \subset G$ and let 
$\mathcal{R}$ consist of all  conjugates of $H$ by elements from $G_0$ 
and $G_0$-cosets of these $H^g$. 
} 

As a result we have a countable dense subgroup $G_0 <G$ so that 
$\mathcal{R}$ is closed under $G_0$-conjugacy and consists of all  
$G_0$-cosets of grey subgroups from $\mathcal{R}$.  
{\em It is easy to see that we may additionally assume that 
the set of grey subgroups from $\mathcal{R}$ is closed 
under $max$ and truncated multiplication 
by positive rational numbers  (i.e. the product is also bounded by 1).} 
\end{remark} 

We will see below that if the space $({\bf Y},d)$ is good enough 
(for example the Urysohn space), then  
the family $\mathcal{R}$ of grey subsets of $G=Iso({\bf Y})$ 
can be chosen among  grey cosets of the form 
$$ 
\rho (g) = q (d(\bar{b}, g(\bar{a}))) \mbox{ , where } q\in \mathbb{Q}^{+} 
\mbox{ and } \bar{a}\bar{b} \mbox{ is } 
$$ 
$$
\mbox{ an appropriate tuple from }S_{{\bf Y}}, 
$$ 
If the metric is bounded by 1 we mean the truncated multiplication by $q$ 
in the formula above. 

When we consider a $(G,G_0 ,\mathcal{R})$-space ${\bf X}$ we usually 
distinguish a similar family too: we choose a countable family 
$\mathcal{U}$ of open (clopen) grey subsets of ${\bf X}$ so that 
a countable family of sets of the form 
$\sigma_{<r}$ for $\sigma \in \mathcal{U}$ and real $r$, 
forms a basis of the topology of ${\bf X}$.
To formalize this let us introduce the following notion.  

\begin{definicja}   \label{gbasis} 
A family $\mathcal{U}$  of open grey subsets of ${\bf X}$  
is called a {\bf grey basis} of the topology $\tau$ if the family 
$\{\phi_{<r}:\phi\in{\mathcal{U}}, r\in {\mathbb{Q}}\cap (0,1)\}$ is a basis of $\tau$.
\end{definicja} 

By Proposition 2.C.2 of \cite{becker} there exists a unique
partition of ${\bf X}$, ${\bf X}=\bigcup\{ Y_{t}: t\in T\}$ into invariant
$G_{\delta}$-sets $Y_{t}$ such that every $G$-orbit from $Y_{t}$ is
dense in $Y_{t}$.
To construct this partition we define for any $t\in 2^{\mathbb{N}}$ the set
$$
Y_{t}=(\bigcap\{ GA_{j}:t(j)=1\})\cap
(\bigcap\{ {\bf X}\setminus GA_{j}:t(j)=0\})
$$ 
where $A_j$ are taken from the corresponding basis of sets the form $\sigma_{<r}$. 
Now let $T=\{ t\in 2^{\mathbb{N}}:Y_{t}\not=\emptyset\}$.
This partition is called {\em canonical}. \parskip0pt

\begin{remark} 
{\em 
When $A$ is open, $GA =A^{\Delta G}$. 
Thus each $Y_t$ as above is an intersection of sets of the form 
$(\sigma_{<r})^{\Delta G}$ or their complements. 
By Lemma \ref{grey} each piece of the canonical partition 
can be constructed as an intersection of sets of the form 
$(\sigma^{\Delta G})_{<r}$ (with $\sigma \in \mathcal{U}$) or their complements. 
Here we may assume that $r$ can take only countably many values. 
}
\end{remark} 

Along with the $d$-topology $\tau$ we shall consider
some special topology on ${\bf X}$.
For this purpose we apply the idea from \cite{becker2} 
of extension of our $\mathcal{U}$ to a {\em nice basis}.

\begin{definicja}  \label{NB}
Let $\mathcal{R}$ be a grey basis of $G$ consisting of cosets of 
open grey subgroups of $G$ which also belong to $\mathcal{R}$. 
Assume that the subfamily of $\mathcal{R}$ of all open grey subgroups 
is closed under $max$ and truncated multiplication by numbers from 
$\mathbb{Q}^{+}$. 

We say that a family $\mathcal{B}$ of Borel grey subsets
of the $G$-space $(\mathbf{X}, \tau )$ is a {\bf nice basis} 
with respect to $\mathcal{R}$ if: \\ 
(i) $\mathcal{B}$ is countable and generates the topology finer than $\tau$;\\ 
(ii) for all $\phi_1, \phi_2 \in \mathcal{B}$, the functions $\neg \phi_1$, 
$min(\phi_1 ,\phi_2)$, $max(\phi_1 ,\phi_2)$, $|\phi_1 - \phi_2 |$, 
$\phi_1 \dot- \phi_2$ $\phi_1 \dot+ \phi_2$ belong to $\mathcal{B}$;\\ 
(iii) for all $\phi \in \mathcal{B}$ and $q\in \mathbb{Q}^{+}$   
the truncated product $q\cdot \phi$   
belongs to $\mathcal{B}$; \\ 
(iv) for all $\phi\in \mathcal{B}$ and open grey subsets  
$\rho \in \mathcal{R}$ we have 
$\phi^{*\rho}, \phi^{\Delta \rho} \in \mathcal{B}$;\\ 
(v) for any $\phi\in \mathcal{B}$ there exists an open grey 
subgroup $H\in \mathcal{R}$ such that $\phi$ is $H$-invariant.
\end{definicja}

It will be usually assumed that the constant function 
$1$ belongs to $\mathcal{B}$, i.e. in particular all 
constant functions $q$, $q\in \mathbb{Q}\cap [0,1]$ 
are in $\mathcal{B}$.  

\begin{definicja} \label{nto}
A topology ${\bf t}$ on ${\bf X}$ is $\mathcal{R}$-{\bf nice} for the $G$-space
$\langle {\bf X}, \tau\rangle$ if the following conditions are satisfied.\\
(a) The topology ${\bf t}$ is Polish, ${\bf t}$ is finer than $\tau$
and the $G$-action remains continuous with respect to ${\bf t}$. \\
(b) There exists a grey basis $\mathcal{B}$ of ${\bf t}$ which is nice with respect to $\mathcal{R}$. 
\end{definicja}

\begin{thm} \label{existence} 
Let $G$ be a Polish group and  $\mathcal{R}$ be a countable grey basis satisfying 
the assumptions of Definition \ref{NB} and  the following closure property: 
\begin{quote} 
for every grey subgroup $H\in\mathcal{R}$ and every $g\in G$ 
if $Hg\in\mathcal{R}$, then $H^g\in\mathcal{R}$\footnote{the conditions of Remark \ref{Gbasis} imply these assumtions}.  
\end{quote} 
Let $\langle {\bf  X}, \tau\rangle $ be a $G$-space and $\mathcal{F}$ 
be a countable family of Borel grey subsets of ${\bf X}$ generating 
a topology finer than $\tau$ such that each $\phi\in \mathcal{F}$  
is invariant with respect to some grey subgroup $H\in \mathcal{R}$ . 
 
Then there is an $\mathcal{R}$-nice topology for $(\langle {\bf  X}, \tau\rangle , G)$ 
such that $\mathcal{F}$ consists of open grey subsets.
\end{thm}

{\em Proof.} 
First we shall  construct an increasing sequence $(\mathcal{S}_n)_{n\in\omega}$
of countable families of Borel grey subsets of $\bf{X}$ along with an 
increasing sequence  $(\mathcal{A}_n)_{n\in\omega}$ 
of countable bases of Polish topologies on $\bf{X}$.
We proceed by induction.
We put
$\mathcal{S}_0=\mathcal{F}$. 
Suppose that we have already constructed $\mathcal{S}_n$. 
Then $\{\phi_{<r}:\phi\in \mathcal{S}_n , r\in (0,1]\cap {\mathbb{Q}}^+\}$ 
is a countable family of Borel subsets of ${\bf X}$. 
Accordingly to Remark 5.1.4 in \cite{bk} it can be extended 
to a countable basis of a Polish topology on $\bf{X}$.
We define $\mathcal{A}_{n}$ to be such a basis.
Then we define $\mathcal{S}_{n+1}$  to be the least family of grey 
subsets extending  $\mathcal{S}_n\cup\{{\bf{O}}_{B}:B\in \mathcal{A}_n\}$
which is closed under conditions (ii)-(iv) of Definition \ref{NB}.
We see that $\mathcal{S}_{n+1}$ is countable and by Lemma \ref{prop} 
it consists of Borel grey subsets.

Having defined sequences $(\mathcal{S}_n)$ and $(\mathcal{A}_n)$  we  put 
$\mathcal{S}=\bigcup\{\mathcal{S}_n:n\in\omega\}$
and $\mathcal{A}=\bigcup\{\mathcal{A}_n:n\in\omega\}$.
Observe that  for every $n\in\omega$ we have
$$
(*)\ \mathcal{A}_n\subseteq \{\phi_{<r}:\phi\in \mathcal{S}_{n+1}, r\in (0,1]\cap {\mathbb{Q}}^+\}
\subseteq \mathcal{A}_{n+1}.
$$
Thus we see that  $\mathcal{A}=\{\phi_{<r}:\phi\in \mathcal{S}, r\in 
(0,1]\cap{\mathbb{Q}}^+\}.$
By Remark 5.1.4 in \cite{bk} $\mathcal{A}$ is a basis of a Polish topology
on $\bf{X}$ finer then the topology generated by $\mathcal{F}$ (say ${\bf t}'$), 
so $\mathcal{S}$ is its  grey basis.
Although the topology defined by grey  basis $\mathcal{S}$  is  Polish,
 in general it  may not preserve continuity of the action.
So we have to look for  some smaller topology. 

Obviously $\mathcal{S}$  satisfies closure conditions (ii)-(iv) of Definition \ref{NB}.
This results in the  following properties of $\mathcal{A}$.

\bigskip

{\bf Claim 1.}  $\mathcal{A}$ is a Boolean algebra of subsets of $\bf{X}$ 
closed under both (standard) Vaught transforms.\bigskip

{\em  Proof. } 
Supposing that $A,B\in \mathcal{A}$, we get  ${\bf O}_A, {\bf  O}_B\in\mathcal{S}$.
Then  by closure properties of $\mathcal{S}$ we  have $1-{\bf O}_A\in \mathcal{S}$ 
and  $ max\{{\bf O}_A, {\bf O}_B\}\in \mathcal{S}$.
Since ${ \bf X}\setminus A=(1-\bf{O}_A)_{<1}$ and  
$A\cap B=(max\{\bf{O}_A, \bf{O}_B\})_{<1}$, thus 
${\bf  X}\setminus A$ and  $A\cap B$ are elements of $\mathcal{A}$.
Next  by Lemma \ref{grey}(3),  for any $\rho\in \mathcal{R}$ and 
$r\in (0,1]\cap {\mathbb{Q}}$ we have $A^{\Delta (\rho_{<r})}=(({\bf O}_A)^{\Delta \rho})_{<r}$, 
thus $A^{\Delta (\rho_{<r})}\in \mathcal{A}$.
\bigskip

Now consider the family  
$\mathcal{B}=\{(\phi^{\Delta\rho}):\phi\in\mathcal{S}, \rho\in \mathcal{R}\}$.
We claim that it has the following properties.

\bigskip

{\bf Claim 2.} (1) If $\psi\in\mathcal{S}$ is $H$-invariant for some 
open grey subgroup $H\in\mathcal{R}$, then $\psi\in\mathcal{B}$.

(2) $\mathcal{B}$  satisfies conditions 
(ii)-(v) of Definition \ref{NB};

(3) $\mathcal{B}=\{(\phi^{*\rho}):\phi\in\mathcal{S}, \rho\in \mathcal{R}\}$.

(4) The family $\{\phi_{<r}:\phi\in\mathcal{B}, r\in (0,1]\cap {\mathbb{Q}}^+\}$ is closed under intersection.
\bigskip

{\em Proof} (1) Accordingly to Lemma \ref{inv} we have
$\psi^{\Delta H}=\psi$ and $\psi^{*H}=\psi$.

(2) First we examine condition (v). 
Take an arbitary $\phi\in\mathcal{S}$ and $\rho\in \mathcal{R}$.
By the hypothesis of the theorem, there are $g\in G$ and an open 
grey subgroup  $H\in\mathcal{R}$ such that $\rho=Hg$. 
Then $\phi^{\Delta\rho}$ is $H^g$-invariant by Lemma \ref{H^g}.

Now take any $\psi_1, \psi_2\in\mathcal{B}$ and open grey subgroup
$H_1,H_2\in\mathcal{R}$ such that $\psi_i$ is $H_i$-invariant, for $i=1,2$.
Put $H=max(H_1,H_2)$, then $\psi_1,\psi_2$ are also $H$-invariant.
Now it is easy to see that both $min(\psi_1 ,\psi_2 )$ and $max(\psi_1,\psi_2)$ 
are $H$-invariant and each of the grey subsets
$\psi_1\dot- \psi_2$,  $\psi_1\dot+ \psi_2$, $|\psi_1-\psi_2|$  
is $2\cdot H$-invariant.
By Lemma \ref{prop} $\neg \psi_1$ is $H_1$-invariant 
and $r\cdot \psi_1$ is $r\cdot H_1$-invariant.
Now the rest follows from point (1).

Since $\mathcal{B}$ is closed under both grey transforms, thus
(iv) follows from the construction of $\mathcal{B}$. 
Point (3) is an immediate consequence of (1).

(4)  Take any 
$A,  B \in  \{\phi_{<r}:\phi\in\mathcal{B}, r\in (0,1]\cap {\mathbb{Q}}^+\}$.
 There are $\phi_A,\phi_B\in \mathcal{S}$ and 
$s,t\in  (0,1]\cap {\mathbb{Q}}^+$ such that
$A=(\phi_A)_{<t}$ and $ B =( \phi_B )_{<s}$.
Then $s\cdot \phi_A, t\cdot  \phi_B \in \mathcal{B}$,  
$A=(s\cdot  \phi_A)_{<st}$,  $ B =(t\cdot   \phi_B )_{<st}$
and $A\cap  B =(max(s\cdot \phi_A, t\cdot  \phi_B ))_{<st}$.
Hence   
$A\cap  B \in \{\phi_{<r}:\phi\in\mathcal{B}, r\in (0,1]\cap {\mathbb{Q}}^+\}$.

\bigskip

{\bf Claim 3.} The  family
$\{\phi_{<r}:\phi\in\mathcal{B}, r\in (0,1]\cap {\mathbb{Q}}^+\}$ is a basis of the topology generated by  the subbasis
$\{A^{\Delta u}: A\in\mathcal{A}, u\in \{\rho_{<s}: \rho\in\mathcal{R},
s\in (0,1]\cap {\mathbb{Q}}^+\}\}$.\bigskip

{\em Proof.} 
Since the family  $\{\phi_{<r}:\phi\in\mathcal{B}, r\in (0,1]\cap {\mathbb{Q}}^+\}$ 
is closed under intersection, it suffices to show 
that  each of its elements is a union of elements of 
the family 
$\{A^{\Delta u}: A\in\mathcal{A}, u\in \{\rho_{<s}: \rho\in\mathcal{R},s\in (0,1]\cap {\mathbb{Q}}^+\}\}$ 
and vice versa.
It  follows directly from the definition of a grey
$\Delta$-transfrom, that every element from 
$\{\phi_{<r}:\phi\in\mathcal{B}, r\in (0,1]\cap {\mathbb{Q}}^+\}$ 
is a union of elements from 
$\{A^{\Delta u}: A\in\mathcal{A}, u\in \{\rho_{<s}: \rho\in\mathcal{R}, s\in (0,1]\cap {\mathbb{Q}}^+\}\}$.
On the other hand  if $A\in \mathcal{A}$, then ${\bf O}_A\in \mathcal{S}$
and $A^{\Delta (\rho_{<s})}=(({\bf O}_A)^{\Delta\rho})_{<s}$.
Hence  
$$
\{A^{\Delta u}: A\in\mathcal{A}, u\in \{\rho_{<s}: 
\rho\in\mathcal{R}, s\in (0,1]\cap {\mathbb{Q}}^+\}\}\subseteq 
\{\phi_{<r}:\phi\in\mathcal{B}, r\in (0,1]\cap {\mathbb{Q}}^+\} .
$$ 

To close the argument  we have to recall Theorem 5.2.1 from \cite{bk}.
Actually, the main part of its proof consists  of justification of  the following statement.
\begin{quote}
Le $G$ be a Polish group, $X$ be a Borel $G$-space and $\mathcal{P}$ be a countable Boolean algebra
of Borel sets generating the Borel structure of $X$ and closed under Vaught transforms.
Then the family 
$$
\{B^{\Delta u}:B\in\mathcal{P}, u \mbox{ - basic open subset of } G\}
$$
is a subbasis of a Polish topology which makes the action continuous.
\end{quote}

Now applying Claim 1 we see that 
$\{A^{\Delta u}: A\in\mathcal{A}, u\in \{\rho_{<s}: \rho\in\mathcal{R}, s\in (0,1]\cap {\mathbb{Q}}^+\}\}$  
is a subbasis of  a Polish topology ${\bf t}$ such that the $G$-action
on ${\bf  X}$ is continuous with rescpect to  ${\bf t}$.
Accordingly to Claim 3, $\mathcal{B}$ is its grey basis.
It follows from the assumptions of the theorem and Claim 2 (1) 
that ${\bf t}$ is finer than ${\bf t}'$.
Then by Claim 2 (2) ${\bf t}$ is a nice topology and $\mathcal{B}$ is its nice grey basis.
$\Box$

\bigskip 

The following observation shows that the assumptions of the theorem above are natural. 

\begin{prop} 
 Let $G$ be a Polish group and let $\mathcal{R}$ be 
a countable basis consisting of open grey cosets.
Let $\langle {\bf  X}, \tau\rangle$ be a $G$-space.
Then there is a grey basis $\mathcal{F}$ for $(\langle {\bf  X}, \tau\rangle , G)$ 
such that each $\phi\in \mathcal{F}$ is invariant with respect to some open grey 
subgroup $H\in {\mathcal{R}}$
\end{prop}

{\em  Proof.} 
Take an arbitrary grey basis  $\mathcal{U}$.
Then  $\mathcal{A}=\{\{\phi_{<r}:\phi\in{\mathcal{U}}, r\in {\mathbb{Q}}\cap (0,1)\}$ is a basis for $\tau$.
We claim that the family 
$$
\{{\bf O}_B^{\Delta H}: B\in {\mathcal{A}},  H\in \mathcal{R}\mbox{  a grey subgroup}\}
$$ 
is a grey basis for $\tau$.
To see this consider any $A\in \mathcal{A}$ and $x_0\in A$. 
By continuity of the action there are an open grey subgroup 
$H\in\mathcal{R}$, $s\in {\mathbb{Q}}\cap (0,1)$ and $B\in\mathcal{A}$ 
such that $x_0\in B\subseteq A$ and for every $x\in B$ and $h\in H_{<s}$ 
the element $hx$ belongs to $A$. 
On the one hand since $x_0\in B$ and ${\bf O}_B$ is open, 
thus by Lemma \ref{inv} (1) we have ${\bf O}_B^{\Delta H}(x_0 )=0$. 
On the other hand if $x\in ({\bf O}_B^{\Delta H})_{<s}$ then there is 
$h\in H_{<s}$ with $hx\in B$ and by $H(h^{-1})=H(h)$ we have $x=h^{-1}hx\in A$. 
Hence $x_0\in ({\bf O}_B^{\Delta H})_{<s}\subseteq A$.
$\Box$ 
\bigskip

\section{Logic action over the Urysohn space} 

The construction of nice bases naturally arises when one 
considers the case of continuous logic actions over $\mathbb{U}$, 
the Urysohn space of diameter 1. 
This is the unique Polish metric space which is universal and 
ultrahomogeneous, i.e. every isometry between finite subsets of 
$\mathbb{U}$ extends to an isometry of $\mathbb{U}$. 
The space $\mathbb{U}$ is considered in 
the continuous signature  $\langle d \rangle$. 

Let $L$ be a countable continuous signature and $\mathbb{U}_L$ 
be the $Iso(\mathbb{U})$-space of all $L$-structures defined 
as in Section 1.1. 
In the present section we will consider 
nice topologies on $\mathbb{U}_L$. 
Our main theorem states that under some natural 
circumstances for any countable fragment $\mathcal{L}$ of 
$L_{\omega_1 \omega}$, grey subsets associated with continuous 
$\mathcal{L}$-formulas form a nice basis on $\mathbb{U}_L$. 
We call it 'logical' and consider it as 
the most natural example of nice bases. 

In Subsection \ref{LopEsc}
we will show that nice topologies induced 
by countable continuous fragments  
of $L_{\omega_1 \omega}$ are ubiquitous.  
In particular Theorem \ref{Cont113} extends Corollary 1.13  
from \cite{becker2} stating that in the case of logic actions of 
$S_{\infty}$ on the space of countable structures 
each nice topology is defined by model sets $Mod (\phi )$ 
of formulas of a countable fragment of $L_{\omega_1 \omega}$.  
We use the version of the L\'{o}pez-Escobar theorem 
for continuous logic from \cite{CL}.  

It is worth noting here that it is commonly accepted that 
$\mathbb{U}$ and $Iso (\mathbb{U})$ play the same role in 
model theory for metric structures as $\omega$ 
and $S_{\infty}$ play in the first order model theory 
(see and \cite{EFP} and Section 2 of \cite{CL} ).  

In Section 3.4 we discuss a very general version of 
the propery of definability of orbits. 
This property naturally appears in 
several places of continuous model theory. 
For example it was applied in \cite{CL}.  

We also find several examples which show some 
obstacles in possible modifications of our results 
(see Section 3.3 and the final part of Section 3.1).

\subsection{Logical nice topologies}

The countable counterpart of $\mathbb{U}$ is the 
{\em rational Urysohn space of diameter 1}, $\mathbb{QU}$, 
which is both ultrahomogeneous and universal for countable 
metric spaces with rational distances and diameter $\le 1$. 
It is shown in Section 5.2 of \cite{BYBM} that there is 
an embedding of $\mathbb{QU}$ into $\mathbb{U}$ so that: \\
(i) $\mathbb{QU}$ is dense in $\mathbb{U}$; \\ 
(ii) any isometry of  $\mathbb{QU}$ extends to an isometry of  
$\mathbb{U}$ and $Iso (\mathbb{QU})$ is dense in $Iso(\mathbb{U})$;  \\ 
(iii) for any $\varepsilon>0$, any partial isometry $h$ of  
$\mathbb{QU}$ with domain $\{ a_1 ,...,a_n\}$ and any isometry 
$g$ of $\mathbb{U}$ such that $d(g(a_i ),h(a_i ))<\varepsilon$ 
for all $i$, there is an isometry $\hat{h}$ of  $\mathbb{QU}$ 
that extends $h$ and is such that for all  $x\in \mathbb{U}$, 
$d(\hat{h}(x),g(x))<\varepsilon$.

We now define a family of clopen grey subgroups of 
$Iso (\mathbb{U})$ which satisfies the conditions of 
Theorem \ref{existence}. 
Let $G_0$ be a dense countable subgroup of $Iso(\mathbb{QU})$.  
By (ii) we may consider it as a subgroup of $Iso (\mathbb{U})$. 

{\bf Family $\mathcal{R}^U$.} 
Let $\mathcal{R}_0$ be the family of all clopen 
grey subgroups of $Iso(\mathbb{U})$ of the (truncated) form 
$$ 
H_{q, \bar{s}} : g\rightarrow q \cdot d(g(\bar{s}), \bar{s}), 
\mbox{ where } \bar{s}\subset \mathbb{QU}, \mbox{ and } q\in \mathbb{Q}^{+}.   
$$ 
It is clear that $\mathcal{R}_0$ is closed under 
conjugacy by elements of $G_0$. 
Consider the closure of $\mathcal{R}_0$ under 
the function $max$ and define $\mathcal{R}^U$ 
to be the family of all $G_0$-cosets of grey 
subgroups from $max(\mathcal{R}_0 )$. 
Then $\mathcal{R}^U$ is countable and the 
family of all $(H_{q, \bar{s}})_{<l}$ where 
$H\in \mathcal{R}_0$ and $l\in \mathbb{Q}$, 
generates the topology of $Iso(\mathbb{U})$.  
Moreover it is easy to see that $G_0$ and 
$\mathcal{R}^U$ satisfy all the conditions of 
Remark \ref{Gbasis} and in particular $\mathcal{R}^U$ 
satisfies the conditions of Theorem \ref{existence}.

\begin{remark} \label{R_0} 
{\em It is worth noting that any grey subgroup 
of the form $max (H_1 ,...,H_k )$ where 
$$ 
H_{i} : g\rightarrow q_i \cdot d(g(\bar{s}_i ), \bar{s}_i ), \mbox{ where } \bar{s}_i\subset \mathbb{QU}, \mbox{ and } q_i \in \mathbb{Q},     
$$ 
has a grey subgroup from $\mathcal{R}_0$. 
Indeed, let $\bar{s}$ consist of all elements appearing in $\bar{s}_i$, 
$i\le k$, and $q = max ( q_1 ,...,q_k )$. 
Then  
$$ 
H_{q, \bar{s}} : g\rightarrow q \cdot d(g(\bar{s}), \bar{s})    
$$ 
is a grey subgroup of $max(H_1 ,...,H_k )$.   
}
\end{remark} 

\bigskip 
 
Let $L$ be a relational language of a continuous signature with 
inverse continuity moduli $\le n\cdot id$ for $n$-ary relations. 
Let $\mathcal{L}$ be a countable fragment of $L_{\omega_1 \omega}$, 
in particular $\mathcal{L}$ be closed under first-order connectives. 
Note that inverse continuity moduli of first-order continuous formulas 
(with connectives as in Introduction) can be taken linear 
(of the form $k\cdot id(x)$). 
Thus it is easy to see that for every formula of $\mathcal{L}$ 
has linear inverse continuity moduli. 

Let $\mathcal{B}$ be the family of all grey subsets defined 
by continuous $\mathcal{L}$-sentences (with parameters) as follows  
$$
\phi (\bar{s}) : M\rightarrow \phi^{M} (\bar{s}), 
\mbox{ where }\bar{s}\in \mathbb{QU} \mbox{ and } \phi (\bar{x}) \in \mathcal{L} . 
$$ 
By linearity of inverse continuity moduli it is easy to see 
that for any continuous sentence $\phi (\bar{s})$ there is a number 
$q\in \mathbb{Q}$ 
(depending on the continuity modulus of $\phi$) such that 
the grey subset as above is $H_{q, \bar{s}}$-invariant. 
We will prove below (in Theorem \ref{mainUrysohn})  
that $\mathcal{B}$ is an $\mathcal{R}^U$-nice basis. 
This can be considered as a version of 
Theorem 1.10 of \cite{becker2} which states that 
in the discrete case of the $S_{\infty}$-space of 
countable $L$-structures all formulas as above 
already form a nice basis. 
%
In Subsection \ref{LopEsc}
we will show that nice topologies induced 
by countable continuous fragments  
of $L_{\omega_1 \omega}$ are ubiquitous.  

\begin{thm} \label{mainUrysohn} 
Let $\mathcal{B}$ be a family of grey subsets corresponding 
to a countable continuous fragment $\mathcal{L}$ of $L_{\omega_1 \omega}$. 
Then the family $\mathcal{B}$ is an $\mathcal{R}^U$-nice basis. 
\end{thm}

We will use Example 56 of \cite{ppetersen}, p.295. 
The main statement of this example can be formulated as follows. 
\begin{quote} 
Let $A=\{ a_1 ,...,a_n \}$ and $B= \{ b_1 ,...,b_n \}$ 
be finite metric spaces so that 
$$ 
|d(a_i ,a_j ) - d(b_i ,b_j )| \le \varepsilon 
\mbox{ , } 1\le i \le j \le n . 
$$ 
Then the disjoint union $A\cup B$ has a metric 
extending the metrics of $A$ and $B$ with the property 
that $d(a_i ,b_i )=\varepsilon$ for every $i\le n$. 
\end{quote}
Section 6 (Appendix) of our paper contains 
a theorem extending this statement. 
It has interesting consequences concerning definability 
of orbits of grey subgroups 
(see Section 3.4).  
\bigskip

{\em Proof of Theorem \ref{mainUrysohn}.} 
To get this statement we use the strategy of 
Theorem 1.10 of \cite{becker2}. 
Conditions (i) - (iii) and (v) of Definition \ref{NB} 
are easily seen in this case. 
For example to see condition (i) it suffices to note 
that $\tau$ is generated by open sets of the form 
$\{ M: |\psi^{M} (\bar{s}) - q_1 |<q_2\}$ where 
$\psi$ is an atomic formula, $\bar{s}\in \mathbb{QU}$ 
and $q_1 ,q_2 \in \mathbb{Q}\cap (0,1)$.  
As in \cite{becker2} we concentrate 
on condition (iv) of the definition of a nice basis.  
The following claim is the main point of the proof. 

\begin{quote} 
Let $B\in \mathcal{B}$ be a grey subset of $\mathbb{U}_L$ 
and $H\in \mathcal{R}^{U}$ be a grey coset of $Iso(\mathbb{U})$ 
corresponding to a grey subgroup of the form $H_{q, \bar{s}}$. 
Then $B^{*H}$ belongs to $\mathcal{B}$. 
\end{quote} 

We fix $B \in \mathcal{B}$ and $H\in \mathcal{R}^U$, and 
find a continuous formula $\phi$, pairwise distinct 
$a_0 ,...,a_{l-1}$, $b_0 ,...,b_{m-1}$, 
$c_0 ,...,c_{n-1}\in \mathbb{QU}$ and pairwise distinct 
$b'_0 ,...,b'_{m-1}$, $c'_0 ,...,c'_{n-1}\in \mathbb{QU}$ 
so that the following three conditions are satisfied: 

(1) the type of $b_0 ,...,b_{m-1}, c_0 ,...,c_{n-1}$ in 
the pure structue $(\mathbb{U},d)$ coincides with 

the type of $b'_0 ,...,b'_{m-1}$, $c'_0 ,...,c'_{n-1}$;  

(2) $H(g) =q\cdot max\{ d(g(b'_0 ),b_0 ),...,d(g(b'_{m-1}),b_{m-1})$, 
$d(g(c'_0 ),c_0 ),...$, $d( g(c'_{n-1}),c_{n-1})\}$ 

for an appropriate $q\in \mathbb{Q}$; 

(3) $\phi= \phi (b_0 ,...,b_{m-1}, a_0 ,...,a_{l-1})$ is 
a continuous $L$-formula with parameters so 

that the grey subset $B$ of $\mathbb{U}_L$  is defined by 
$$ 
M\rightarrow B(M) = \phi^{M}(b_0 ,...,b_{m-1}, a_0 ,...,a_{l-1}). 
$$ 
Let $\rho : [0,1] \rightarrow [0,1]$ be an inverse continuiuty modulus 
for $\phi$. 
By our assumptions on the continuous signature $L$, 
we may assume that $\rho$ is a linear function 
(see Lemma \ref{ContMod}).  
We want to show that $B^{*H}$ belongs to $\mathcal{B}$. 

%
By Example 56 of \cite{ppetersen} for any subspace 
$\{ b''_0 ,...,b''_{m-1}, c''_0,...,c''_{n-1},a''_0 ,...,a''_{l-1}\}\subset \mathbb{U}$ 
and any $\varepsilon$ so that for each pair 
$s_1 ,s_2 \in \{ b_0 ,...,b_{m-1}, c_0 ,...,c_{n-1}, a_0 ,...,a_{l-1}\}$ 
and the corresponding 
$s''_1 ,s''_2 \in \{ b''_0 ,...,b''_{m-1}, c''_0 ,...,c''_{n-1}, a''_0 ,...,a''_{l-1}\}$ 
we have 
$$
|d(s_1 ,s_2 )-d(s''_1 ,s''_2 )|\le \varepsilon ,  
$$ 
there is an embedding $\alpha$ of the  space 
$\{ b_0 ,...,b_{m-1}, c_0 ,...,c_{n-1}, a_0 ,...,a_{l-1}\}$ 
into $\mathbb{U}$ so that 
$$ 
\mbox{ for each }  s_i \in \{ b_0 ,...,b_{m-1}, c_0 ,...,c_{n-1}, a_0 ,...,a_{l-1}\} 
\mbox{ , } d(\alpha(s_i ),s''_i )\le \varepsilon . 
$$  
For such $\varepsilon$ and $\alpha$ viewing 
$$
\phi^{M}(b''_0 ,...,b''_{m-1}, a''_0 ,...,a''_{l-1})\dot- 
q\cdot max (d(b_0 ,b''_0 ),...,d(b_{m-1},b''_{m-1}), d(c_0 ,c''_0 ),..., d( c_{n-1},c''_{n-1}))
$$ 
as 
$$
[\phi^{M}(\bar{b}'',\bar{a}'')- \phi^{M}(\alpha(\bar{b)}, \alpha(\bar{a})) 
+\phi^{M}(\alpha(\bar{b}), \alpha(\bar{a}))] \dot- 
$$ 
$$ 
[q\cdot max (d(\bar{b} ,\bar{b}''),d(\bar{c},\bar{c}'')) - 
q\cdot max(d(\bar{b},\alpha(\bar{b})),d(\bar{c},\alpha(\bar{c}))) + 
$$ 
$$ 
q\cdot max (d(\bar{b},\alpha(\bar{b})),d(\bar{c},\alpha(\bar{c})))], 
$$ 
we see (by the triangle condition) that it is not greater than 
$$
\phi^{M}(\alpha(\bar{b}), \alpha(\bar{a}))\dot- q\cdot max (d(\bar{b},\alpha(\bar{b})),d(\bar{c},\alpha(\bar{c})))
+\rho (\varepsilon )+ q\varepsilon. 
$$ 

Let $\psi (u'_0 ,...,u'_{m-1},v'_0 ,...,v'_{n-1})$ 
be the following formula: 
$$
sup_{\bar{w},\bar{u},\bar{v}} [[  \phi(\bar{u}, \bar{w}) \dot- 
q\cdot max (d(u'_0 ,u_0 ),...,d(u'_{m-1},u_{m-1}), d(v'_0 ,v_0 ),..., d( v'_{n-1},v_{n-1}))] \dot-
$$ 
$$
(\rho + q\cdot id)( max (|d(z_1 ,z_2 ) - d(s_1 ,s_2 )|  :  \mbox{ the place of }
$$ 
$$ 
z_1 ,z_2 \mbox{ in } \{ u_0 ,...,u_{m-1}, v_0 ,...,v_{n-1}, w_0 ,...,w_{l-1} \}  \mbox{ corresponds to } 
$$
$$
\mbox{ the place of } s_1 ,s_2 \mbox{ in } \{ b_0 ,...,b_{m-1}, c_0 ,...,c_{n-1}, a_0 ,...,a_{l-1}\} ))]. 
$$ 
To see that $B^{*H}$ is determined by 
$\psi (b'_0 ,...,b'_{m-1},c'_0 ,...,c'_{n-1})$  
note that for a tuple 
$a''_0 ,...,a''_{l-1}, b''_0 ,...,b''_{m-1},c''_0 ,...,c''_{n-1}$ 
isomorphic to $a_0 ,...,a_{l-1}, b_0 ,...,b_{m-1},c_0 ,...,c_{n-1}$ 
and for any $L$-expansion $M$ of $(\mathbb{U},d)$ we have 
$$
\phi^{M} (\bar{b}'', \bar{a}'') \dot- 
q\cdot max (d(b'_0 ,b''_0 ),...,d(b'_{m-1},b''_{m-1}), d(c'_0 ,c''_0 ),..., d( c'_{n-1},c''_{n-1}))
\le \psi^{M} (\bar{b}',\bar{c}') .
$$  
In particular for any $r\in [0,1]$ greater than  
$\psi^{M} (b'_0 ,...,b'_{m-1},c'_0 ,...,c'_{n-1})$ 
if $g$ maps each $b''_i$ to $b_i$ and each $c''_i$ to $c_i$ 
(i.e. $H(g)=q\cdot d(\bar{b}'\bar{c}', \bar{b}''\bar{c}'')$), 
then $g$ maps $M$ to $B_{<r\dot+H(g)}$  (by 
$\phi^{g(M)}(\bar{b},\bar{a})=\phi^{M}(\bar{b}'' ,\bar{a}'' )$).  
We see that for any decomposition $r=r'- r''$  
all isometries from $H_{<r''}$ take $M$ to $B_{<r'}$. 

On the other hand if the expansion $M$ 
does not satisfy 
$$
\psi (b'_0 ,...,b'_{m-1},c'_0 ,...,c'_{n-1})< r, 
$$ 
then for any $\delta >0$ there is a tuple $\bar{a}''\bar{b}''\bar{c}''$ 
such that $r-\delta$ is less than 
$$ 
[\phi^{M}(\bar{b}'', \bar{a}'') \dot- 
q\cdot max (d(b'_0 ,b''_0 ),...,d(b'_{m-1},b''_{m-1}), d(c'_0 ,c''_0 ),..., d( c'_{n-1},c''_{n-1}))]
\dot- (\rho (q_1 )+q\cdot q_1)
$$
where 
$$
q_1 = max (|d(z_1 ,z_2 ) - d(s_1 ,s_2 )|  :  
\mbox{ the place of } z_1 ,z_2 \mbox{ in } 
\{ b''_0 ,...,b''_{m-1}, c''_0 ,...,c''_{n-1}, a''_0 ,...,a''_{l-1} \}  
$$
$$
\mbox{ corresponds to the place of } s_1 ,s_2 \mbox{ in } 
\{ b_0 ,...,b_{m-1}, c_0 ,...,c_{n-1}, a_0 ,...,a_{l-1}\} ) . 
$$  
Then as we already know by Example 56 of \cite{ppetersen} 
there is a map $\alpha$ taking  $\bar{b}'\bar{c}'\bar{a}'$ 
to a tuple in $\mathbb{U}$, say 
$\alpha (\bar{b}'\bar{c}'\bar{a}' )$, 
which is at distance 
$\le q_1$ from $\bar{b}''\bar{c}''\bar{a}''$. 
In particular, 
$$
\phi (\bar{b}'',\bar{a}'')\dot- 
q\cdot max (d(b'_0 ,b''_0 ),...,d(b'_{m-1},b''_{m-1}), d(c'_0 ,c''_0 ),..., d( c'_{n-1},c''_{n-1}))
$$ 
does not differ from  
$$ 
\phi (\alpha(\bar{b}'),\alpha (\bar{a}')) \dot- q\cdot max 
(d(b'_0 ,\alpha (b'_0) ),...,d(b'_{m-1},\alpha(b'_{m-1})), d(c'_0 ,\alpha(c'_0 )),..., d( c'_{n-1},\alpha (c'_{n-1})))
$$ 
more than $\rho (q_1 ) + q\cdot q_1 $. 
In particular we may replace $\bar{b}''\bar{c}''\bar{a}''$ 
by $\alpha(\bar{b}'\bar{c}'\bar{a}')$ 
keeping the value of the formula above greater than 
$r - \delta$ (substituting $\alpha(\bar{b}'\bar{c}'\bar{a}')$ 
we see that the final part of the formula disappears). 
As a result we may assume that 
$\bar{b}''\bar{c}''\bar{a}''$ and $\bar{b}\bar{c}\bar{a}$ 
are in the same orbit under $Iso (\mathbb{U})$. 

This argument shows that the value of 
$\psi (b'_0 ,...,b'_{m-1},c'_0 ,...,c'_{n-1})$ 
becomes $sup$ of the corresponding 
subformula  with respect to subspaces  
$a''_0 ,...,a''_{l-1}, b''_0 ,...,b''_{m-1},c''_0$, $...,c''_{n-1}$ 
isomorphic to $a_0 ,...,a_{l-1}$, 
$ b_0 ,...,b_{m-1},c_0 ,...,c_{n-1}$. 
Since any isometry $h$ taking 
$a''_0 ,...$, $a''_{l-1}, b''_0 ,...,b''_{m-1},c''_0 ,...,c''_{n-1}$ 
as above to $a_0 ,...,a_{l-1}, b_0 ,...,b_{m-1},c_0 ,...,c_{n-1}$ 
maps $M$ to $B_{> r + H(h) -\delta}$,     
for sufficiently small $\varepsilon$
there is a presentation $r -\delta = r' -r''$ 
(for example $r'' = H(h) +\varepsilon$) 
such that  the basic open set of all isometries of 
$\mathbb{U}$ taking  
$b''_0 ,...,b''_{m-1}, c''_0 ,...,c''_{n-1}, a''_0 ,...,a''_{l-1}$  
to the $\varepsilon$-neighborhood of 
$b_0 ,...,b_{m-1}, c_0 ,...,c_{n-1}, a_0$ , $...,a_{l-1}$ 
is contained in $H_{<r''}$ but  
does not contain an element taking $M$ to $B_{<r' -\varepsilon}$.   
In particular, $M\not\in B^{*H}_{<r-\delta}$.  
Thus $M\not\in B^{*H}_{<r}$.

To finish the proof of the theorem note 
that the proof given above can be easily 
generalisied to the case of a grey 
coset of the form $max (H_1 ,...,H_k )$ 
where each $H_i$ is defined as $H$ above. 
The modification concerns the form of the formula 
$\psi$: after the $\phi$-part we should apply  
$max$ to an appropriate linear functions.  
$\Box$

\bigskip 

Can Theorem \ref{mainUrysohn} be generalised 
to the statement that the family $\mathcal{B}$   
is an $\mathcal{R}$-nice basis for any  
family $\mathcal{R}$ of $G_0$-cosets of 
grey subgroups of $Iso( \mathbb{U})$ 
satisfying our standard assumptions? 
We finish this subsection by an example which shows 
that the argument of the proof of Theorem \ref{mainUrysohn} 
does not work for this generalisation. 

\bigskip 

{\bf Family $\mathcal{R}^{\sqrt{}}$.} 
We now define an extension of $\mathcal{R}^U$.  
Let $\mathcal{R}^{\sqrt{}}_{0}$ be the extension of 
$\mathcal{R}_{0}$ by grey subsets of the form 
$$ 
H^{\sqrt{}}_{q,\bar{s}} : g\rightarrow q \cdot \sqrt{d(g(\bar{s}), \bar{s})}, 
\mbox{ where } \bar{s}\subset \mathbb{QU}, \mbox{ and } q\in \mathbb{Q}^{+}.   
$$ 
To see that grey sets of this form are grey 
subgroups take $g_1 ,g_2 ,g_3 \in Iso(\mathbb{U})$ 
with $g_1 \cdot g_2 = g_3$. 
Since all $g_i$ are isometries,  
$d(g_1 g_2 (\bar{s}),\bar{s} ) \le d(g_1 (\bar{s}), \bar{s}) + d(g_2 (\bar{s}), \bar{s})$. 
Thus 
$\sqrt{d(g_1 g_2 (\bar{s}),\bar{s} )} \le \sqrt{d(g_1 (\bar{s}), \bar{s})} + \sqrt{d(g_2 (\bar{s}), \bar{s})}$ 
which implies the required inequality.
When we apply $max$ to a finite family of grey 
subgroups we obviously obtain a grey subgroup too. 
Let $\mathcal{R}^{\sqrt{}}$ be the family of all $G_0$-cosets 
of grey subgroups from $max(\mathcal{R}^{\sqrt{}}_0 )$. 
It is clear that it satisfies the corresponding assumptions 
of Theorem \ref{existence}. 
It is also clear that the corresponding version of Remark \ref{R_0} 
holds for subgroups from $\mathcal{R}^{\sqrt{}}_0$.

Let $\mathcal{L}_0$ be the first order 
fragment of $L_{\omega_1 \omega}$. 
Let $\mathcal{B}_0$ be the corresponding 
family of grey subsets of $\mathbb{U}_L$.  
The family $\mathcal{B}_0$ is an 
$\mathcal{R}^U$-nice basis and by Theorem \ref{existence}   
$\mathcal{B}_0$ can be also extended to an 
$\mathcal{R}^{\sqrt{}}$-nice basis of the $G$-space 
$\mathbb{U}_{L}$.  
We {\em conjecture} that $\mathcal{B}_0$ 
is not an $\mathcal{R}^{\sqrt{}}$-nice basis. 
The following remark supports this conjecture. 

\begin{remark} \label{notNice} 
{\em 
Let $L$ be the language corresponding to the continuous signature 
$\langle d, c\rangle$, where $c$ is a constant symbol. 
Let $u_0 \in \mathbb{U}$ and $\phi (u_0 )$ be the formula 
$10\cdot d(u_0 ,c)$ (with truncated product). 
Let $B_0$ be the grey subset of $\mathbb{U}_L$ defined by $\phi(u_0 )$ 
and let $H$ be the grey subgroup of $Iso (\mathbb{U} )$ defined by 
$$ 
g\rightarrow \sqrt{d(g(u_0 ), u_0 )} .   
$$ 
Since we do not include $\sqrt{}$ into the set of our connectives, 
the grey subset $B^{* H}_0$ of $\mathbb{U}_L$ 
cannot be defined as in the the proof of Theorem \ref{mainUrysohn}. 
Thus the proof of Theorem \ref{mainUrysohn} 
should be essentially modified in this case 
(assuming that the statement holds). 
We have some additional observation 
showing that it is difficult to avoid $\sqrt{}$ 
in these considerations. 
Let $\theta = B^{\Delta H}_0$. 

\bigskip 

{\bf Claim.} 
If a structure $M\in \mathbb{U}_L$ satisfies 
$\frac{1}{10} < d (u_0 ,c) < \frac{1}{2}$ 
then $\theta (M) = \sqrt{d(u_0 ,c)}$. 

Indeed let $q=d( u_0 ,c)$. 
Consider an isometric copy of $[0,1]$ in 
$\mathbb{U}$ so that $0$ is identified with 
$c$ and the number $q$ with $u_0$.  
Take any $a$ between $u_0$ and $c$.  
Let $\varepsilon = d(a,c)$. 
Consider an isometry $g$ taking $c$ to $a$ 
so that $d(g(u_0 ),u_0 )=\varepsilon$ 
(for example when $g$ acts along the segment $[0,1]$). 
Then $H (g) = \sqrt{\varepsilon}$ and the structure 
$g(M)$ satisfies $d(u_0 ,c) =  q - \varepsilon$ 
(in $g(M)$ the symbol $c$ is interpreted by 
the $g$-image of the interpretation of $c$ in $M$).  
In particular for every small $\delta$ we have a non-meagre 
part of $H_{<\sqrt{\varepsilon}}$ taking $M$ into  
$\{ M' : \phi (M' )<10 (q -\varepsilon) +\delta \}$. 
Thus $\theta (M)\le 10 (q - \varepsilon) + \sqrt{\varepsilon}$. 

Since for every $g$ with $H(g) = \sqrt{\varepsilon}$ 
the distance $d(g(c) ,u_0)$ is not less than $q-\varepsilon$, 
the value $\theta (M)$ coincides with 
$10 (q - \varepsilon) + \sqrt{\varepsilon}$ 
for an appropriate $\varepsilon$. 
Moreover, by elementary calculus we see that this 
number becomes the infinum only when $\varepsilon = q$ 
(here we use the assumption $\frac{1}{10} < d (u_0 ,c)$). 
This proves the claim. 
$\Box$ 
}
\end{remark} 

\subsection{When is a nice topology logical?}  \label{LopEsc}

Corolorary 1.13 of \cite{becker2} states that 
any nice topology of the $S_{\infty}$-space of discrete countable structures  
is defined by a countable fragment of $L_{\omega_1 \omega}$. 
The crucial point of the proof is the L\'{o}pez-Escobar theorem 
that in the $S_{\infty}$-space of countable structures any invariant 
Borel subset is defined by a formula of $L_{\omega_1 \omega}$. 
In this subsection we discuss some versions of 
the L\'{o}pez-Escobar theorem in the space of continuous structures 
on $\mathbb{U}$ and continuous versions of Corolorary 1.13 of \cite{becker2}. 

Let us fix some continuous language $L$,  
a sublanguage $L_0$ and  an ultrahomogeneous 
continuous $L_0$-structure $M_0$ on some 
Polish space ${\bf Y}$. 
Let $G=Aut(M_0 )$ and let  
${\bf Y}_L$ be the $G$-space  
of all continiuos $L$-structures on ${\bf Y}$. 

S.Coskey, M.Lupini, I.Ben Yaacov, A.Nies and T.Tsankov 
have proved in \cite{BNT} and \cite{CL} some versions 
of Theorem \ref{CLthm} given below in the case 
when $G= Iso(\mathbb{U})$. 
The argument from \cite{CL} works for 
the full statement of Theorem \ref{CLthm}.  
In fact the authors mention this 
after Corollary 1.2 in \cite{CL} when they 
consider the case where $M_0$ 
is just an ultrahomogeneous space.  

\begin{thm} \label{CLthm}
(S.Coskey and M.Lupini) 
Assume that ${\bf Y}$ contains 
a dense sequence $\{ p_i : i\in \omega \}$ so that 
for every $n$ the $G$-orbit of $(p_0 ,...,p_n )$ 
is definable in $M_0$. 
Let  $\lambda$ be a Borel grey subset of ${\bf Y}_L$ 
which is $G$-invariant. 

Then there is an $L_{\omega_1 \omega}$-formula $\phi$ 
which defines $\lambda$ by the map $M \rightarrow \phi^M$, 
where $M\in {\bf Y}_L$. 
\end{thm}    

The following statement  
appears in the proof of Theorem \ref{CLthm} 
in \cite{CL} as Theorem 3.2. 
\begin{quote} 
Under circumstances of Theorem \ref{CLthm} 
let ${\bf Y} = \mathbb{U}$ and let $\mu$ be 
a Borel grey subset of $\mathbb{U}_L$. 

Then for any natural number $n$ there is 
an $L_{\omega_1 \omega}$-formula $\phi (\bar{u})$ with 
$k$ free variables, so that for every subsequence 
$p_{i_1} ,...,p_{i_k }$ of $\{ p_i :i\in \omega\}$, 
if $H$ is the grey stabilizer 
$$ 
g \rightarrow n\cdot d(g(p_{i_1}),...,g(p_{i_k}),p_{i_1} ,...,p_{i_k }), 
$$ 
then the $*$-transform $\mu^{*H}$ is defined 
by $\phi ^M (p_{i_1} ,...,p_{i_k })$, $M\in \mathbb{U}_L$. 
\end{quote} 

\begin{cor} \label{CLcor}
Let us fix some continuous language $L$ and the corresponding 
space $\mathbb{U}_L$ of $L$-structures on $\mathbb{U}$. 
Let $H$ be a grey subgroup from the family $\mathcal{R}^U$ 
defined in Section 3.1 and let $\lambda$ be a Borel grey subset 
of $\mathbb{U}_L$ which is $H$-invariant. 

Then there is an $L_{\omega_1 \omega}$-formula $\phi$ 
over $\mathbb{U}$ which defines $\lambda$ by 
the map $M \rightarrow \phi^M$. 
\end{cor}    

{\em Proof.}  
Let us enumerate $\mathbb{QU}$ as 
a sequence $p_0 ,...,p_i ,...$ . 
It is observed in \cite{CL} 
(the discussion after Lemma 3.3) that  
the $G$-orbit of any $(p_0 ,...,p_i )$ 
is definable in $\mathbb{U}$. 
By Remark \ref{R_0} we may assume that 
$H$ is of the form of the grey stabiliser 
given above.   
Applying Theorem 3.2 of \cite{CL} in 
the form above where $\mu = \lambda$, 
we obtain that $\lambda^{*H}$ is defined 
by an appropriate $L_{\omega_1 \omega}$-formula 
with parameters from $\mathbb{QU}$. 
Since $\lambda$ is $H$-invariant by Lemma \ref{inv} 
we have $\lambda^{*H} = \lambda$, i.e  
the formula $\phi$ realises $\lambda$. 
$\Box$

\bigskip

We now give a version of  Corollary 1.13 of \cite{becker2}.

\begin{thm} \label{Cont113} 
Let $L$ be a continuous language and  $G=Iso (\mathbb{U})$. 
Consider the logic $G$-space $\mathbb{U}_L$ under the standard 
topology $\tau$. 
Let $\mathcal{F}$ 
be a countable family of Borel grey subsets of $\mathbb{U}_L$ 
generating a topology finer than $\tau$ such that each $\phi\in \mathcal{F}$  
is invariant with respect to some grey subgroup  $H\in \mathcal{R}^U$ . 
 
Then there is an $\mathcal{R}^U$-nice topology for $(\langle \mathbb{U}_L , \tau\rangle , G)$ 
generated by some countable fragment 
of $L_{\omega_1 \omega}$ 
such that $\mathcal{F}$ consists of open grey subsets.
\end{thm}

{\em Proof.} 
By Theorem \ref{existence} we find an appropriate 
$\mathcal{R}^U$-nice topology for $(\langle \mathbb{U}_L , \tau\rangle , G)$  
such that $\mathcal{F}$ consists of open grey subsets.
Let us consider the grey subsets of the nice basis of this topology. 
By Corollary \ref{CLcor} each subset of this basis 
is determined by an appropriate 
$L_{\omega_1 \omega}$-formula. 
Take the fragment generated by these formulas. 
Repeating the proof of Claim 2 of Theorem \ref{existence} 
we see that this fragment defines an $\mathcal{R}^U$-nice topology. 
$\Box$

\begin{remark} 
{\em 
We do not know if the statement 
of Theorem \ref{Cont113} holds for 
$\mathcal{R}^{\sqrt{}}$. 
For example let $\mathcal{B}^{\sqrt{}}$ 
be the $\mathcal{R}^{\sqrt{}}$-nice basis 
generated by $\mathcal{B}_0$ 
(it exists by Theorem \ref{existence}). 
We do not know if  
$\mathcal{B}_0 =\mathcal{B}^{\sqrt{}}$.   
Is $\mathcal{B}^{\sqrt{}}$ generated by a fragment of 
$L_{\omega_1 \omega}$?    
If this is not true we can view it as the best example 
of so called 'non-logical' basis.  
}
\end{remark} 

\bigskip

Let us consider possible generalisations of 
results of this section. 
Let us consider a Polish metric structure $(M_0 ,d)$. 
According to Defintion 5.3 of \cite{BYBM}  
a (classical) countable structure $N$
is called a countable {\em approximating substructure}
of $M_0$ if the universe of $N$ is a dense countable subset 
of $M_0$ and any automorphism of $N$ 
extends to a (necessarily unique) automorphism of
$(M_0 ,d)$ so that $Aut(N)$ is dense in $Aut(M_0 )$.
We endow $Aut(N)$  with the Polish group 
topology of point-wise convergence in the discrete set
$N$, which refines the topology induced as a subset
of $Aut(M_0 )$.
It is clear that $\mathbb{QU}$ is an approximating 
substructure of $\mathbb{U}$.  
Note that by Lemma 5.6 of \cite{BYBM} 
for every open subset $U\subset Aut(M_0 )$ 
the intersection $U\cap Aut(N)$ is open in $Aut(N)$. 

Repeating the construction of Section 3.1 
we define a family of clopen grey subgroups of 
$G=Aut (M_0)$,  say  $\mathcal{R}^0$,  which satisfies 
the conditions of  Theorem \ref{existence}. 

Let $G_0$ be a dense countable subgroup of $Aut(N)$.  
Let $\mathcal{R}_0$ be the family of all clopen 
grey subgroups of $G=Aut(M_0 )$ of the (truncated) form 
$$ 
H_{q, \bar{s}} : g\rightarrow q \cdot d(g(\bar{s}), \bar{s}), 
\mbox{ where } \bar{s}\subset N, \mbox{ and } q\in \mathbb{Q}^{+}.   
$$ 
Consider the closure of $\mathcal{R}_0$ under 
the function $max$ and define $\mathcal{R}^0$ 
to be the family of all $G_0$-cosets of grey 
subgroups from $max(\mathcal{R}_0 )$. 
Then it is easy to see that $G_0$ and 
$\mathcal{R}^0$ satisfy all the conditions of 
Remark \ref{Gbasis} and in particular $\mathcal{R}^0$ 
satisfies the conditions of Theorem \ref{existence}. 
The corresponding version of Remark \ref{R_0} 
also holds in this case. 
The following theorem generalises 
Theorem \ref{Cont113}.

\begin{thm} 
Assume that $M_0$ is an ultrahomogeneous continuous 
$L_0$-structure on some Polish  spaces ${\bf Y}$. 
Let a language $L$ extend $L_0$ and 
${\bf Y}_L$ be the $G$-space of all continuous 
$L$-structures on ${\bf Y}$ with the standard 
topology $\tau$. 

Assume that $N$ can be enumerated into 
a dense subsequence $\{ p_i : i\in \omega \}$ of 
${\bf Y}$  so that for every $n$ the $G$-orbit of $(p_0 ,...,p_n )$ 
is definable in $M_0$. 

Let $\mathcal{F}$ be a countable family of Borel grey subsets of ${\bf Y}_L$ 
generating a topology finer than $\tau$ such that each $\phi\in \mathcal{F}$  
is invariant with respect to some grey subgroup $H\in \mathcal{R}^0$. 
 
Then there is an $\mathcal{R}^0$-nice topology for $(\langle {\bf  Y}, \tau\rangle , G)$ 
generated by some countable fragment 
of $L_{\omega_1 \omega}$ 
such that $\mathcal{F}$ consists of open grey subsets.
\end{thm} 

The proof repeats the proof of Theorem \ref{Cont113}. 

\subsection{Impossible versions of the L\'{o}pez-Escobar theorem}

In this subsection we discuss some versions of 
the L\'{o}pez-Escobar theorem in the space of continuous structures. 
Let us fix some continuous language $L$ and the corresponding 
space $\mathbb{U}_L$ of $L$-structures on $\mathbb{U}$. 
S.Coskey, M.Lupini, I.Ben Yaacov, A.Nies and T.Tsankov 
have proved in \cite{BNT} and \cite{CL} that when 
$\bar{a}\in \mathbb{U}$ and $\lambda$ is a Borel grey subset of 
$\mathbb{U}_L$ which is $Iso(\mathbb{U}/\bar{a})$-invariant 
there is an $L_{\omega_1 \omega}$-formula $\phi$ realising $\lambda$, 
where the formula $\phi$ is of the form $\phi' (\bar{a})$, 
where $\bar{a}$ is a tuple from $\mathbb{U}$.  
Note that (contrary to the classical case) it can happen that 
the corresponding formula $\phi' (\bar{x})$ with variables $\bar{x}$ 
instead of $\bar{a}$ is not an $L_{\omega_1 \omega}$-formula. 
This suggests the following considerations. 

Assuming that $L$ has constant symbols let $L^{tame}_{\omega_1 \omega}$ 
consist of all formulas of the form  $\phi(\bar{z})= \phi' (\bar{c},\bar{z})$, 
where  $\bar{c}$ is a tuple of $L$-constants and 
$\phi' (\bar{x},\bar{z})$ is a continuous $L_{\omega_1 \omega}$-formula 
without constants. 
It is clear that $L^{tame}_{\omega_1 \omega}$ is a fragment of 
$L_{\omega_1 \omega}$ which is closed under finitary iterations 
of basic continuous connectives. 
Note that when $L$ is discrete,  
$L^{tame}_{\omega_1 \omega}= L_{\omega_1 \omega}$.  

We now give an example which shows that the 
L\'{o}pez-Escobar theorem does not hold for  $L^{tame}_{\omega_1 \omega}$. 
 
Let $L=\langle d, R^{1}, c\rangle$ be a continuous signature 
where $c$ is a constant symbol and $R$ is a symbol of 
a predicate with continuity modulus $id$. 
Consider the logic space $\mathbb{U}_L$ of all continuous 
$L$-structures on the Urysohn space $\mathbb{U}$.  
To each continuous $L$-structure $M\in \mathbb{U}_L$ 
we associate a real number $r_M$ defined as follows:     
$$ 
M \rightarrow r_M = sin(\frac{1}{R(c)^M}). 
$$ 
We assume that $r_M = 0$ for $R(c)^M =0$. 
Since $R(c)$ is a formula and $sin(\frac{1}{x})$ is continuous 
in $(0,1]$, the function $M\rightarrow r_M$ is Borel. 
It is obviously $Iso(\mathbb{U})$-invariant. 

We want to prove that this function is not defined by 
a continuous $L^{tame}_{\omega_1 \omega}$-formula.

\begin{prop} \label{nLE} 
There is no $L^{tame}_{\omega_1 \omega}$-sentence $\psi$ so that 
for any $M\in \mathbb{U}_L$ the numbers $r_M$ and $\psi^M$ are the same. 
\end{prop}

{\em Proof.}
If the constant symbol $c$ does not appear in 
a continuous $L_{\omega_1 \omega}$-formula $\psi$, 
then for any $(L\setminus \{ c\})$-structure on $\mathbb{U}$ 
of the form $R(x)=d(x,u_0 )$ with some fixed $u_0\in \mathbb{U}$ 
the value $\psi^M$ coincides with $\psi^{M'}$ where 
$M=\langle \mathbb{U}, d, R,c\rangle$ with $c=u_0$ and 
$M'=\langle \mathbb{U}, d, R,c'\rangle$ with $d(c',u_0 )= \frac{2}{\pi }$. 
This means that $\psi$ cannot define the function $M\rightarrow r_M$. 

Suppose that $\phi (x)$ is a continuous 
$L_{\omega_1 \omega}$-formula without $c$ so that 
for any $M\in \mathbb{U}_L$, $r_M = \phi (c)^M$. 
Let $u_0 \in \mathbb{U}$ and $M_0$ be the structure of the signature 
$\langle d,R \rangle$ on $\mathbb{U}$ where $R(x) = d(x,u_0 )$. 

Let $\delta_{\phi ,x}$ be a continuity modulus of $\phi$.  
Then for any two $L$-expansions of $M_0$ (say by $c$ and $c'$) 
$$ 
d(c,c' )\le \delta_{\phi ,x} (\varepsilon ) \Rightarrow |r_{(M_0,c)} - r_{(M_0 ,c')}|<\varepsilon . 
$$ 
Choosing $c,c'$ in an isometric copy of $[0,1]$ in $\mathbb{U}$ 
(identifying $u_0$ with $0$) 
we obtain a continuity modulus for $\sin (\frac{1}{x})$ in $(0,1]$. 
This is a contradiction. 
$\Box$


\subsection{The property of definability of orbits of grey subgroups}

The property of $\mathbb{U}$ that any $n$-orbit 
of $Iso (\mathbb{U})$ is definable, is necessary  
for the version of the L\'{o}pez-Escobar theorem 
from \cite{CL} (see Theorem \ref{CLthm} above). 
In this subsection we show a very similar property of 
the space $\mathbb{U}_L$. 
We formulate it in the most general form.

\begin{definicja} \label{APGS}
Let $(G,\mathcal{R})$ and $({\bf X},\mathcal{U})$ satisfy 
the assumptions of Remark \ref{Gbasis} and Definition \ref{gbasis}.  
We say that the $(G,\mathcal{R})$-space $({\bf X},\mathcal{U})$ 
has {\bf definable orbits of grey subgroups} 
if for any grey subgroup $H\in \mathcal{R}$ and 
for any $\varepsilon >0$ there is an $H$-invariant $\phi\in {\cal U}$ 
and $\delta >0$ such that for any $c$ and $c' \in {\bf X}$ 
the following property holds: \\ 
  
if $|\phi (c) - \phi (c')|\le \delta$,  
then there is $g\in G$ with $H(g)<\varepsilon$   
and $d(g(c'),c)< \varepsilon$. 
\end{definicja} 

Let us consider the $G$-space $\mathbb{ U}_L$ under 
any $\mathcal{R}^U$-nice basis ${\cal B}$ containing 
the family ${\cal B}_0$  of grey subsets defined by 
all continuous first-order $L$-sentences over parameters 
from $\mathbb{QU}$ (see the final part of Section 3.1).  
We consider $\mathbb{U}_L$ under the metric 
$\delta_{seq(\mathbb{QU})}$ defined as in Section 1.1 
with respect to some enumeration of $seq(\mathbb{QU})$. 
To see that $\mathbb{U}_{L}$ has the property of 
definable orbits of grey subgroups note that by 
Remark \ref{R_0} it suffices to verify the formulation 
for grey subgroups of the form 
$$ 
H_{q,\bar{s}} : g\rightarrow q \cdot d(g(\bar{s}), \bar{s}), \mbox{ where } \bar{s}\subset \mathbb{QU}, \mbox{ and } q\in \mathbb{Q}.     
$$ 

By the definition of the space 
$\mathbb{U}_L$ to guarantee that two $L$-structures $M$ 
and $N$ are distant $\le \varepsilon$ it suffices to 
find an appropriate tuple $\bar{a}\subset \mathbb{QU}$ 
and appropriate $\varepsilon'>0$ so that 
$\bar{a}$ contains tuples of a safficiently 
long initial segment of the enumeration of 
$seq(\mathbb{QU})$ and the values of 
the corresponding atomic $L$-formulas 
(say $\psi_i$, $i\in I$) on subtuples of $\bar{a}$ 
do not differ in $M$ and $N$ more than $\varepsilon'$. 

 
By the values of $\psi_i$'s in $M$ we build 
an $inf$-sentence $\phi$ over $\bar{s}$ 
so that if for appropriately chosen $\delta$,   
$|\phi^M - \phi^N|\le \delta$, 
then there is $\bar{a}'$ such that  
any $\psi^M_i (\bar{s}\bar{a})$ does not differ 
from $\psi^N_i (\bar{s}\bar{a}')$ 
more than $\frac{\varepsilon'}{2}$. 
To realise this for any $\psi_i (\bar{s}\bar{a})$ 
find a rational number $r_i$ with 
$|\psi^M_i (\bar{s}\bar{a})-r_i|\le \frac{\varepsilon'}{4}$.
Then the $inf$-formula $\phi$ just expresses 
the existence of $\bar{a}'$ so that for any $i$,  
$|\psi_i (\bar{s}\bar{a}')-r_i|\le \frac{\varepsilon'}{4}$.   

Now we use Theorem \ref{Udist} from Appendix. 
It in particular says that  there is  sufficiently 
small $\delta >0$ such that if distances between pairs 
in $\bar{s}\bar{a}$ and  the corresponding distances in 
$\bar{s}\bar{a}'$ do not differ more than  $\delta$ ,  
then there is an isometry $g\in Iso (\mathbb{U})$ 
(apply ultrahomogenity) such that 
$d(g(\bar{s}\bar{a}'), \bar{s}\bar{a} )\le \frac{\varepsilon'}{2}$. 
Choosing such $g$  and using the assumption that $id$ 
is a continuity modulus of any relation from $L$, 
we obtain that the  values of 
appropriate atomic $L$-formulas in $g(N)$ 
on $\bar{s}\bar{a}$ does not differ from 
the corresponding values in $M$ more than $ \varepsilon'$.

As a result we obtain the following proposition. 

\begin{prop} \label{apprU} 
The property of definability of orbits of grey subgroups 
holds in the space $\mathbb{U}_L$  under any nice 
basis ${\cal B}$ containing the family of 
grey subsets defined by all continuous first-order $L$-sentences 
over parameters from $\mathbb{QU}$. 
\end{prop}

\section{Canonical pieces which are $G$-orbits}

We now consider canonical pieces (see Section 2.2) 
which are $G$-orbits. 
We take the assumptions of Remark \ref{Gbasis}, 
Definition \ref{NB} and Definition \ref{nto}, 
i.e. in particular $G$ is a Polish group with a dense countable 
subgroup $G_0$ and a grey basis $\mathcal{R}$ 
consisting of clopen grey $G_0$-cosets   
(as in the case of Theorem \ref{mainUrysohn}). 
By $\mathcal{B}$ and ${\bf t}$ we denote a nice basis and 
the corresponding nice topology of 
$\langle {\bf X}, \tau \rangle$, which is a Polish $G$-space.

We start with a proposition which is a version of Lindstrom's model
completeness theorem (that an $\forall\exists$-axiomatizable
$\aleph_0$-categorical theory is model complete).

\begin{prop}\label{embed}
Let $X =Gx_0$ for some (any) $x_0\in X$ and $X$
be a $G_{\delta}$-subset of ${\bf X}$. 
Then both topologies $\tau$ and ${\bf t}$ are equal on $X$.
\end{prop}

{\em Proof.}
We have to check  that  every element of
${\cal B}$ is a $\tau$-open grey subset on $X$.
Let $\phi\in {\cal B}$, $x_1\in X$ and $\phi (x_1 )=r_1<r$.
Take an open grey subgroup $H$ of $G$  from $\mathcal{R}$  
such that $\phi$ is $H$-invariant.
Let $r'< r-r_1$. 
Since by the Effros' theorem on $G_{\delta}$-orbits,
the canonical map $g\rightarrow gx_1$ is an open map
$G\rightarrow X$, we see that $H_{<r'}x_1$ is a $\tau$-open
subset of $\phi_{<r}\cap X$ containing $x_1$. $\Box$ 
\bigskip

In the following definition we introduce farther counterparts 
of model theoretic notions.

\begin{definicja}  
Let ${\bf t}$ be an $\mathcal{R}$-nice topology 
for the $G$-space $\langle {\bf X},\tau \rangle$.
Let $H$ be an open grey subgroup from $\mathcal{R}$ 
and $X$ be an invariant $G_{\delta}$-subset of 
${\bf X}$ with respect to the ${\bf t}$-topology. \\
(1) A  family ${\cal F}$ of subsets of the form 
$\phi_{< r}$ with $H$-invariant $\phi\in {\cal B}$ 
is called an $H$-{\bf type} in $X$, if it is 
maximal with respect to the condition that 
$X\cap\bigcap {\cal F}\not=\emptyset$.\\
(2) An $H$-type ${\cal F}$ is called {\bf principal} 
if for every $\varepsilon >0$ there is an $H$-invariant 
grey basic set $\phi \in \mathcal{B}$ and there is 
$r$ such that $\phi_{<r}\in {\cal F}$ and for each 
$B=\psi_{<t}\in {\cal F}$, 
the set $\psi_{<t+\varepsilon}$ 
contains $\phi_{< r} \cap X$.
Then we say that $\phi_{< r}$ 
$\varepsilon$-{\bf defines} ${\cal F}$. 
\end{definicja}

Each type ${\cal F}$ is determined by 
any element from $X\cap\bigcap {\cal F}$. 
It is worth noting that when $x$ and $z$ determine 
different $H$-types, then there is an $H$-invariant 
$\psi$ and two rational numbers $r_1$ and $r_2 \le 1$ 
so that $\psi (x)<r_1$, $1-\psi (z)< r_2$ and 
$r_1 < 1- r_2$.   

\begin{lem} \label{compform} 
Assume that $c$  defines a principal 
$H$-type ${\cal F}$ and $\phi_{<r}$  
$\frac{\varepsilon}{2}$-defines  
${\cal F}$ (with an $H$-invariant grey basic 
set $\phi \in \mathcal{B}$ and $r>0$). 

Then for any $H$-invariant $\psi$ 
and any  $c'$ and $c''$ with 
$max(\phi (c'), \phi (c''))<r$, 
we have $|\psi (c')-\psi (c'')|<\varepsilon$. 
In particular $\phi_{<r}$ $\varepsilon$-defines 
the $H$-type of any element of $\phi_{<r}$. 
\end{lem} 

{\em Proof.}  
Note that 
$|\psi (c') -\psi (c)|<\frac{\varepsilon}{2}$. 
Indeed, if $\psi (c)< \psi (c')$, then applying 
the definition of principal types,  
$\psi (c')<\psi (c) +\frac{\varepsilon}{2}$. 
If $\psi (c')<\psi (c)$, then 
repeat this argument for 
$1-\psi (x)$.
The rest is clear. 
$\Box$ 

\bigskip   

Note that  the following lemma is related 
to omitting types theorems from logic.

\begin{lem} \label{->princ}
Let $G$, $\mathcal{R}$, $\mathcal{B}$ and ${\bf t}$ 
satisfy the assumptions of the section and $X$ 
be an invariant ${\bf t}$-$G_{\delta}$-subset.
Let $H\in \mathcal{R}$.
Then for any non-principal $H$-type
${\cal F}$ the set
$\bigcap\limits_{g\in G}\bigg(g(\bigcup\limits_{B\in {\cal F}}(X\setminus B))\bigg)$
is  nonempty and $G$-invariant.
In particular if $X$ is a $G$-orbit 
then any $H$-type of $X$ is principal.
\end{lem}

{\em Proof.} 
Suppose ${\cal F}$ is non-principal. 
Find $\varepsilon >0$ so that no $\phi_{<t}$ 
$\varepsilon$-defines ${\cal F}$, i.e. 
by Lemma \ref{compform} there is no 
$\phi_{<t}$ in $\mathcal{F}$ 
so that for any $H$-invariant $\psi$ 
and any  $c'$ and $c''$ with 
$max(\phi (c'), \phi (c''))<r$, 
we have $|\psi (c')-\psi (c'')|<\varepsilon$.  
We claim that 
$\bigcap \{ H_{<\frac{\varepsilon}{2}} B: B\in {\cal F} \}$ 
is meager in $(X,{\bf t})$.
Indeed, otherwise it would be comeager in
$X\cap B_{\cal F}$ for some non-empty $B_{\cal F}$ 
of the form $\phi_{<r}$ with $\phi \in {\cal B}$. 
We fix these $\phi$ and $r$. 
Since $\mathcal{B}$ contains constant functions 
we may assume that $r< \frac{\varepsilon}{2}$.  
Consider $\phi^{\Delta H}$. 
It is an $H$-invariant grey set from $\mathcal{B}$.  
Since $\phi$ is open, $\phi^{\Delta H}$ is open too.  
If $B\in {\cal F}$ with $B=\psi_{< t}$, then 
the closure of $H_{<\frac{\varepsilon}{2}}B\cap X$ 
contains $B_{{\cal F}}\cap X$.    
Thus by Lemma \ref{greyold} 
$(\phi^{\Delta H})_{<r}$ 
is contained in    
$\overline{H_{<r}(\psi^{\Delta H})_{< t +\frac{\varepsilon}{2}}}$ 
and $H$-invariantness of $\psi$ 
(i.e. in particular 
$(\psi^{\Delta H})_{< t} = \psi_{< t }$ and 
$H_{<\varepsilon} \psi_{< t }\subseteq \psi_{<t+\varepsilon}$)  
implies that  
$(\phi^{\Delta H})_{<r}$ is contained in  
$\psi_{\le t+r+\frac{\varepsilon}{2}}$, 
i.e. in $\psi_{< t +\varepsilon}$. 
This is a contradiction with the assumption that ${\cal F}$ is non-principal.
\parskip0pt

Now 
$X\cap \bigcup\limits_{B\in {\cal F}}({\bf X}\setminus H_{<\frac{\varepsilon}{2}}B)$
is comeagre in $X$. 
Since this set is contained in 
$$ 
\bigcap\{ X\cap h(\bigcup\limits_{B\in {\cal F}}({\bf X}\setminus B)): h\in H_{<\frac{\varepsilon}{2}}\} , 
$$ 
for every $g\in G$, the set 
$$ 
g(\bigcap\{ X\cap h(\bigcup\limits_{B\in {\cal F}}({\bf X}\setminus B)): h\in H_{<\frac{\varepsilon}{2}}\} ), 
$$ 
is comeager in $X$.

The grey group $H$ is of countable index in $G$ 
(by $H\sqsubseteq_o G$), i.e. for every $\varepsilon'$ 
the group $G$ is covered by countably many 
left translates of $H_{<\varepsilon'}$.
In particular we may choose a countable subgroup 
$G_{cd} \subset G$ such that for every $g\in G$  
there is $g'\in G_{cd} \cap gH_{<\frac{\varepsilon}{2}}$. 
It is clear that 
$$ 
\bigcap\limits_{g\in G_{cd}}
g(\bigcap\{ X\cap h(\bigcup\limits_{B\in {\cal F}}({\bf X}\setminus B)): h\in H_{<\frac{\varepsilon}{2}}\} )  
$$ 
is comeagre in $X$. 
Note that 
$$ 
\bigcap\limits_{g\in G_{cd}}
g(\bigcap\{ X\cap h(\bigcup\limits_{B\in {\cal F}}({\bf X}\setminus B)): h\in H_{<\frac{\varepsilon}{2}}\} ) 
 = \bigcap\limits_{g\in G}\bigg(g(\bigcup\limits_{B\in {\cal F}}(X\setminus B))\bigg) . 
$$ 
Indeed if 
$x\in \bigcap\limits_{g\in G_{cd}} g(\bigcap\{ X\cap h(\bigcup\limits_{B\in {\cal F}}({\bf X}\setminus B)): h\in H_{<\frac{\varepsilon}{2}}\} )$  
and $g\in G$ then find $g'\in G_{cd}\cap gH_{<\frac{\varepsilon}{2}}$ with $g=g'h'$, $h'\in H_{<\frac{\varepsilon}{2}}$. 
Thus  
$$ 
(g')^{-1}(x)\in \bigcap\{ X\cap h(\bigcup\limits_{B\in {\cal F}}({\bf X}\setminus B)): h\in H_{<\frac{\varepsilon}{2}}\} 
$$ 
and 
$(h')^{-1}(g')^{-1}(x) \in X\cap \bigcup\limits_{B\in {\cal F}}({\bf X}\setminus B)$, 
i.e. $x\in X\cap g(\bigcup\limits_{B\in {\cal F}}({\bf X}\setminus B))$.  
The rest is clear.  
We see that the intersection 
$\bigcap\limits_{g\in G}(g(\bigcup\limits_{B\in {\cal F}}(X\setminus B)))$ 
is $G$-invariant and comeagre.

To prove the remaining part suppose that ${\cal F}$ is
a non-principal $H$-type of $X$ and $X$ is
a $G$-orbit.
Then, by the previous statement,
$X\subseteq \bigcup\limits_{B\in {\cal F}}(X\setminus B)$.
This contadicts to the definition of a type.
$\Box$

\bigskip

The following statement is a version 
of Ryll-Nardzewski's theorem.
We remind the reader that we assume that 
grey cosets from $\mathcal{R}$ are clopen 
and are represented by elements of $G_0 <_{dense} G$. 
The group $G$ is considered under a left-invariant metric 
(for example defined as $\rho_S$ in Section 0). 

\begin{thm} \label{RN}
Let $({\bf X},\mathcal{U})$ be a Polish $(G,\mathcal{R})$-space and 
${\bf t}$ be a nice topology of ${\bf X}$ with the nice basis $\mathcal{B}$. 
A piece $X$ of the canonical partition with respect to the
topology ${\bf t}$ is a $G$-orbit if and only if 
for any basic clopen grey  subgroup $H\in \mathcal{R}$
any $H$-type of $X$ is principal.
\end{thm}

{\em Proof.} 
By Lemma \ref{->princ} we have the necessity
of the theorem. \parskip0pt

For sufficiency we use the back-and-forth argument.
Let $x,y\in X$. 
We build a set $\Gamma$ of tuples
$(H_{i},\phi_{i},r_{i} ,\varepsilon_{i},H'_{i},\phi'_{i},g_{i})$, 
$i\in \omega$, with the following properties:

(a) each $H_{i}$ (and $H'_{i}$) is a clopen grey 
subgroup from $\mathcal{R}$ and $g_i$ belongs to $G_0$; \\ 
each $\phi_{i}$ is an $H_{i}$-invariant basic 
${\bf t}$-clopen grey subset with $x \in (\phi_{i} )_{< r_{i}}$ 
and each $\phi'_{i}$ is an $H'_{i}$-invariant basic 
${\bf t}$-clopen grey subset with  $y\in (\phi'_{i})_{< r_{i}}$; 

(b)  for each even $i>0$, 
$H_{i+1}\sqsubset max(H_{i}, H_{i-1})$, 
$H'_{i+1} = g_{i+1}H_{i+1}g^{-1}_{i+1}$ 
(i.e. $= (^{g_{i+1}} H_{i+1})$
\footnote{by $H^{g}$ we denote the grey group $g^{-1}hg \rightarrow H(h)$ 
and by $^{g}H$  we denote the grey group $ghg^{-1}\rightarrow H(h)$ } 
), 
$\phi_{i+1}\sqsubseteq max(\phi_{i}, \phi_{i-1})$, 
$\phi'_{i+1} = g_{i+1} \phi_{i+1}$ 
and  the $H_{i+1}$-type of $x$ is $\frac{\varepsilon_{i+1}}{2}$-defined 
by $(\phi_{i+1} )_{< r_{i+1}}$ as a principal type;                

(c) for each odd $i>0$, 
$H'_{i+1}\sqsubset max(H'_{i},H'_{i-1})$, 
$H_{i+1} = g^{-1}_{i+1}H'_{i+1}g_{i+1}$, 
$\phi'_{i+1}\sqsubseteq max(\phi'_{i},\phi'_{i-1})$,   
 $\phi_{i+1} = g^{-1}_{i+1} \phi'_{i+1}$ 
and the $H'_{i+1}$-type of $y$ is 
$\frac{\varepsilon_{i+1}}{2}$-defined 
by $(\phi'_{i+1})_{< r_{i+1}}$ as a principal type;                 

(d) for each $i$, $\varepsilon_{i} \le 2^{-i}$, 
$$
max(diam((\phi_i )_{\le r_i}),diam((H_{i})_{\le\varepsilon_i}),
diam((\phi'_i )_{\le r_i}),diam((H'_{i})_{\le\varepsilon_i}))<2^{-i} ,
$$ 
$$ 
H_{i} (g^{-1}_{i} g_{i+1} ) \le \varepsilon_{i}                     
\mbox{ and  } 
H'_{i} (g_{i+1} g^{-1}_{i} )\le \varepsilon_{i}                      
\mbox{ for all  } i. 
$$ 

It is worth noting here that since $G$ is 
considered under a left-invariant metric,
$diam((g_{i}H_{i})_{<\varepsilon_i})= diam((H_{i})_{<\varepsilon_i})<2^{-i}$.

At Step 0 let $H_0=H'_0=H_{-1} =H'_{-1} =G$,  
$g_0 =id$, $\varepsilon_0 =1$ and $\phi_0=\phi'_0$ 
be a $G$-invariant basic grey set so that for some 
$r_0$, $X\subset (\phi_0 )_{<r_0}$. 

At step $i+1$ (assuming that $i$ is even) 
take any basic $C\subset (\phi_{i})_{<r_{i}}$ 
of the form $\psi_{<r}$ where $\psi \in \mathcal{B}$
is a  (${\bf t}$-clopen)  grey subset of 
$max(\phi_i ,\phi_{i-1})$, with $x\in C$ and $r$ so small 
that $diam(C)<2^{-(i+1)}$. 
Let $H_{i+1}$ be a basic clopen grey subgroup of 
$max (H_{i},H_{i-1})$ such that $\psi$ is $H_{i+1}$-invariant. 
We may choose 
$H_{i+1}$ so that 
for some $\varepsilon<2^{-(i+1)}$,   
$diam((H_{i+1} )_{\le\varepsilon})<2^{-(i+1)}$.  

We denote this $\varepsilon$ by $\varepsilon_{i+1}$. 
Let $(\phi_{i+1})_{<r_{i+1}}\subseteq C$ 
$\frac{\varepsilon_{i+1}}{2}$-define the 
(principal) $H_{i+1}$-type of $x$. 
We may assume that $\phi_{i+1} \sqsubseteq \psi$ 
and $r_{i+1}\le r$ for $\psi$ and $r$ above. 
\parskip0pt

Let $g_{i+1}$ be any element of $G_0$ 
which maps the clopen set $(\phi_{i+1})_{<r_{i+1}}$ 
to a set containing $y$. 
The existence of such $g_{i+1}$ follows 
from the assumption that $x$ and $y$ 
belong to the same canonical piece. 
On the other hand by the inductive assumptions,  
$g_{i}\phi_{i+1}$ is an $^{g_{i}} H_{i+1}$-invariant 
clopen grey  subset of $\phi'_{i}$. 
Let us compute $H'_i (g_i g^{-1}_{i+1})$. 
By Lemma \ref{compform} for any $H'_i$-invariant 
grey subset $\theta \sqsubseteq \phi'_i$ 
the values $\theta (g_i g^{-1}_{i+1} y)$ and $\theta(y)$ 
do not differ more than $\varepsilon_i$. 
Since $\theta$ is $H'_i$-invariant, 
$H'_{i} ( g_i g^{-1}_{i+1})<\varepsilon_i$. 
Since $H_{i} = g^{-1}_{i}H'_{i}g_{i}$, 
$H_{i} (g^{-1}_{i+1}g_i ) = H'_{i}(g_i g^{-1}_{i+1}g_i g^{-1}_i )<\varepsilon_i$,  
i.e. the corresponding part of condition (d) is satisfied. 

Let $H'_{i+1}=(^{g_{i+1}}H_{i+1})$ and $\phi'_{i+1}=g_{i+1}\phi_{i+1}$.
It is clear that properties (a)-(d) are satisfied for $i$.
\parskip0pt 

The case of odd $i$ is symmetric.\parskip0pt

As a result we have a Cauchy sequence $g_{i}$, $i\in \omega$.
Let $g\in G$ be the limit of this sequence. 
Note that for each $i$ the element $g_{i}$ maps $(\phi_{i})_{<r_i}$ 
to $(\phi'_{i})_{<r_i}$ and 
$$ 
\{ x\} =\bigcap \{ (\phi_{i})_{<r_i}:i\in \omega\} \mbox{ and } 
\{ y\} = \bigcap \{ (\phi'_{i})_{<r_i}:i\in\omega \} .
$$ 
Thus $g$ maps a Cauchy sequence converging to $x$ to a Cauchy sequence 
converging to $y$.  
Thus $g$ maps $x$ to $y$.
$\Box$

\section{Complexity of some subsets of the logic space}

Viewing the logic space ${\bf Y}_L$ as a Polish space 
one may consider Borel/algorithmic complexity of interesting subsets of ${\bf Y}_L$. 
In this section we fix a countable dense subset $S_{{\bf Y}}$ of ${\bf Y}$ 
and study subsets of  ${\bf Y}_L$ which are invariant with respect to 
isometries stabilising $S_{{\bf Y}}$ setwise. 
The best example of this situation is the logic space $\mathbb{U}_L$ over 
the  bounded Urysohn metric space $\mathbb{U}$ where distinguishing 
the countable counterpart $\mathbb{QU}$ of $\mathbb{U}$ (see Section 3) 
we study $Iso (\mathbb{QU})$-invariant subsets of $\mathbb{U}_L$. 

This approach corresponds to considering a structure on ${\bf Y}$ 
(say $M$) together with its {\em presentation over} $S_{{\bf Y}}$, 
i.e. the set 
$$ 
Diag (M,S_{{\bf Y}})=\{ (\phi ,q): M\models \phi <q, \mbox{ where } q\in [0,1]\cap \mathbb{Q} 
\mbox{ and } \phi \mbox{ is a continuous } 
$$ 
$$ 
\mbox{ sentence with parameters from }S_{{\bf Y}} \}.     
$$ 
In this section we examine separable categoricity and ultrahomogenity. 
In particular in Section 5.1 we find a Borel subset $\mathcal{SC}$ of 
${\bf Y}_L$ which is $Iso (S_{{\bf Y}})$-invariant and each separably 
categorical structure on ${\bf Y}$ is homeomorphic 
to a structure from $\mathcal{SC}$.  

Since any Polish group can be realised as 
the automorphism group of a approximately ultrahomogeneous 
structures, it makes sense in order to characterise some special 
properties of Polish groups to study the corresponding subclasses 
of approximately ultrahomogeneous structures and then to study 
the complexity of these subclasses. 
In Section 5.2 we consider two opposite subsets of approximately 
ultrahomogeneous structures from ${\bf Y}_L$ which correspond 
to natuaral properties of automorphism groups:  
separable oligomorphicity and admitting of complete left invariant metrics. 

In Section 5.3 we study complexity of the index set of computable members 
of $\mathcal{SC}$. 
We believe that these ideas can be applied to investigation of complexity 
of other topological and model theoretic notions. 
In fact our intention in this section is to demonstrate some 
new settings arising in the approach of the logic space of continuous 
structures. 
  
Our special attension to separable categoricity is motivated 
by ubiquity of it in this paper: the Urysohn space is separably categorical, 
the material of Section 4  is an abstract form of Ryll-Nardzewski's theorem,  
definability of orbits of grey subgroups for 
logic spaces over the Urysohn space (Section 3.4) 
is a consequence of categoricity.

\subsection{Separable categoricity}

We preserve all the assumptions of Section 1 on 
the space $({\bf Y},d)$. 
For simplicity we assume that all $L$-symbols are 
of continuity modulus $id$.  
We reformulate separable categoricity as follows.

\begin{prop} \label{SCHm} 
Let $M$ be a non-compact, separable, continuous, 
metric structure on $({\bf Y},d)$.
The structure $M$ is separably categorical if and only 
if  for any $n$ and $\varepsilon$ there are finitely many  
conditions $\phi_i (\bar{x})\le \delta_i$, $i\in I$,  
so that any $n$-tuple of $M$ satisfies one of these 
conditions and the following property holds: 
\begin{quote} 
for any  $i\in I$, any $a_1 ,...,a_n \in M$ realising  
$\phi_i (\bar{x})\le \delta_i$ and  
any finite set of formulas $\Delta (x_1 ,...,x_n , x_{n+1})$ 
realised in $M$ and containing  $\phi_i (\bar{x})\le \delta_i$, 
there is a tuple $b_1 ,...,b_n ,b_{n+1}$ 
realising $\Delta$ such that $max_{i\le n} d(a_i ,b_i ) <\varepsilon$.   
\end{quote} 
\end{prop} 

To prove this proposition we start with the following observation. 

\begin{lem} \label{SCH} 
Let $M$ be a non-compact, separable, continuous, 
metric structure on $({\bf Y},d)$.
The structure $M$ is separably categorical if and only if 
for any $n$ and $\varepsilon$ there are finitely many  
conditions $\phi_i (\bar{x})\le \delta_i$, $i\in I$,  
so that any $n$-tuple of $M$ satisfies one of these 
conditions and for any  $i\in I$, any $a_1 ,...,a_n \in M$ 
realising  $\phi_i (\bar{x})\le \delta_i$ and any type 
$p(x_1 ,...,x_n , x_{n+1})$ realised in $M$ and containing  
$\phi_i (\bar{x})\le \delta_i$, there is a tuple 
$b_1 ,...,b_n ,b_{n+1}$ realising $p$ such that 
$max_{i\le n} d(a_i ,b_i ) <\varepsilon$.   
\end{lem} 

{\em Proof.} 
By Theorem 12.10 of \cite{BYBHU} a complete theory $T$ is separably categorical 
if and only if for each $n>0$, every $n$-type is principal. 
An equivalent condition states that for each $n>0$, the metric space 
$(S_n (T),d)$ is compact.  
In particular for every $n$ and every $\varepsilon$ there is 
a finite family of principal $n$-types $p_1 ,...,p_m$ so that 
their $\varepsilon/2$-neighbourhoods cover $S_n(T)$. 

Thus when $M$ is separably categorical, given $n$ and $\varepsilon$, 
we find appropriate $p_i$, $i\in I$, define    
$P_i (\bar{x})=dist(\bar{x}, p_i (M))$, 
the corresponding definable predicates and  
$n$-conditions $\phi_i (\bar{x})\le \delta_i$ describing the corresponding 
$\varepsilon /2$-neighbourhoods of $p_i$.  
The rest follows by strong $\omega$-near-homogenity. 

To see the converse assume that $M$ satisfies the property from the 
formulation. 
To see that $G=Aut (M)$ is approximately oligomorphic 
take any $n$ and $\varepsilon$ and find finitely many  
conditions $\phi_i (\bar{x})\le \delta_i$, $i\in I$,  
satisfying the property from the formulation for $n$ and $\varepsilon /4$. 
Choose $\bar{a}_i$ with $\phi_i (\bar{a}_i )\le \delta_i$ and 
let $F =\{ \bar{a}_i : i\in I\}$.  
To see that $G\cdot F$ is $\varepsilon$-dense 
we only need to show that if $\bar{a}$ satisies 
$\phi_i (\bar{x})\le \delta_i$, then there is an automorphism 
which takes $\bar{a}$ to the $\varepsilon$-neighbourhood 
of $\bar{a}_i$. 
This is verified by "back-and-forth" as follows. 
Let $(\varepsilon_k )$ be an infinite sequence of positive 
real numbers whose sum is less than $\varepsilon /4$. 
At every step $l$ (assuming that $l\ge n$) we build 
a finite elementary map $\alpha_l$ and $l$-tuples 
$\bar{c}_l$ and $\bar{d}_l$ so that  
\begin{itemize} 
\item $\bar{c}_n =\bar{a}$ and $\bar{d}_n =\bar{a}_i$;  
\item for $l>n$, $\alpha_l$ takes $\bar{c}_l$ to $\bar{d}_l$
\item for $l>n+1$, the first $l-1$ coordinates of $\bar{c}_l$ 
(resp. $\bar{d}_l$) are at distance less than  $\varepsilon_l$ 
away from the corresponding coordinates of $\bar{c}_{l-1}$ 
(resp. $\bar{d}_{l-1}$); 
\item the sets $\bigcup\{ \bar{c}_l :l\in \omega\}$ and 
$\bigcup\{ \bar{d}_l: l\in\omega\}$ are dense in $M$. 
\end{itemize} 
In fact we additionally arrange that for even $l$, 
$\bar{c}_{l+1}$ extends $\bar{c}_l$ and for odd $l$ 
$\bar{d}_{l+1}$ extends $\bar{d}_l$. 
In particular the type of  $\bar{c}_{l+1}$ 
always extends the type of $\bar{c}_l$.   
At the $(n+1)$-th step we find finitely many conditions 
$\phi'_j (\bar{x})\le \delta'_j$, $j\in J$,  
so that any $(n+1)$-tuple of $M$ satisfies one of these 
conditions and for any  $j\in J$, any $a'_1 ,...,a'_{n+1} \in M$ 
realising  $\phi'_j (\bar{x})\le \delta'_j$ and any type 
$p(x_1 ,...,x_{n+1} , x_{n+2})$ realised in $M$ and 
containing  $\phi'_j (\bar{x})\le \delta_j$, there is 
a tuple $b_1 ,...,b_{n+1} ,b_{n+2}$ realising $p$ such 
that $max_{t\le n+1} d(a'_t ,b_t ) <\varepsilon_{n+1}$. 
Now by the choice of $i$ for any extension 
of $\bar{a}=\bar{c}_{n}$ to an $(n+1)$-tuple $\bar{c}_{n+1}$ 
we can find a tuple $\bar{d}_{n+1}$ realising $tp(\bar{c}_{n+1})$  
so that the first $n$ coordinates of $\bar{d}_{n+1}$ are 
at distance less than $\varepsilon /4$ away from 
the corresponding coordinates of $\bar{d}_{n}=\bar{a}_i$.  
If $n$ is even we choose such $\bar{c}_{n+1}$ and 
$\bar{d}_{n+1}$; if $n$ is odd we replace the roles 
of $\bar{c}_{n+1}$ and $\bar{d}_{n+1}$. 
For the next step we fix the condition 
$\phi'_j (\bar{x})\le \delta'_j$ satisfied 
by $\bar{c}_{n+1}$ and $\bar{d}_{n+1}$. 

The $(l+1)$-th step is as follows. 
Assume that $l$ is even (the odd case is symmetric). 
Extend $\bar{c}_l$ to an appropriate $\bar{c}_{l+1}$ 
(aiming to density of $\bigcup \{ \bar{c}_l :l\in \omega\}$). 
There are finitely many conditions $\phi''_k (\bar{x})\le \delta''_k$, 
$k\in K$, so that any $(l+1)$-tuple of $M$ satisfies one of these 
conditions and for any  $k\in K$, any $a'_1 ,...,a'_{l+1} \in M$ 
realising  $\phi''_k (\bar{x})\le \delta''_k$ and any type 
$p(x_1 ,...,x_{l+1} , x_{l+2})$ realised in $M$ and 
containing  $\phi''_k (\bar{x})\le \delta''_k$, there is 
a tuple $b_1 ,...,b_{l+1} ,b_{l+2}$ realising $p$ such that 
$max_{t\le l+1} d(a'_t ,b_t ) <\varepsilon_{l+1}$.   
We find the condition satisfied by $\bar{c}_{l+1}$ 
and a tuple $\bar{d}_{l+1}$ realising $tp(\bar{c}_{l+1})$  
so that the first $l$ coordinates of $\bar{d}_{l+1}$ are at distance 
less than $\varepsilon_{l}$ away from the corresponding 
coordinates of $\bar{d}_{l}$. 

As a result for every $k$ we obtain Cauchy sequences of 
$k$-restrictions of $\bar{c}_l$-s and $\bar{d}_l$-s. 
For $k=n$ their limits are not distant from $\bar{a}$ 
and $\bar{a}_i$ more than $\varepsilon/2$. 
Moreover the limits $lim \{\bar{c}_l \}$ and $lim \{\bar{d}_l \}$ 
are dense subsets of ${\bf Y}$ and realise the same type. 
This defines the required automorphism of $M$. 
$\Box$

\bigskip 

{\em Proof of Proposition \ref{SCHm}.} 
It suffices to show that the condition of the 
formulation implies the corresponding condition of Lemma \ref{SCH}.  
Given $n$ and $\varepsilon$ take the family 
$\phi_i (\bar{x})\le \delta_i$, $i\in I$, satisfying 
the condition of the proposition for $n$ and $\varepsilon /2$. 
Let $p(\bar{x},x_{n+1})$ be a type with  
$\phi_i (\bar{x})\le \delta_i$ and $a_1 ,...,a_n$ 
be as in the formulation. 

Let $(\varepsilon_k )$ be an infinite sequence of positive 
real numbers whose sum is less than $\varepsilon /2$. 
Now apply the condition of the formulation of the proposition 
to $n+1$ and $\varepsilon_1 /2$ and find an appropriate 
finite family of inequalities such that one of them, 
say  $\psi (\bar{x},x_{n+1}) \le \tau$, belongs to $p$ 
and for any $c_1 ,...,c_n ,c_{n+1} \in M$ 
realising  $\psi (\bar{x},x_{n+1})\le \tau$, 
and any finite subset $\Delta \subset p$ containing  
$\psi (\bar{x},x_{n+1})\le \tau$
there is a tuple $c'_1 ,...,c'_n ,c'_{n+1}$ realising  $\Delta$,  
such that $max_{i\le n+1} d(c_i ,c'_i ) <\varepsilon_1 /2$. 
Then let    $b^1_1 ,...,b^1_n ,b^1_{n+1}$ be a tuple realising   
$\phi_i (\bar{x})\le \delta_i$ and 
$\psi (\bar{x},x_{n+1})\le \tau$ such that 
$max_{i\le n} d(a_i ,b^1_i ) <\varepsilon /2$. 

For $n+1$ and $\varepsilon_2 /2$ find an appropriate condition 
$\psi' (\bar{x},x_{n+1}) \le \tau'$ from $p$ 
so that  any $c_1 ,...,c_n ,c_{n+1} \in M$ realising  
$\psi' (\bar{x},x_{n+1})\le \tau'$, 
and any finite subset $\Delta \subset p$ containing  
$\psi' (\bar{x},x_{n+1})\le \tau'$
there is a tuple $c'_1 ,...,c'_n ,c'_{n+1}$ realising  $\Delta$,  
such that $max_{i\le n+1} d(c_i ,c'_i ) <\varepsilon_2 /2$. 
Let    $b^2_1 ,...,b^2_n ,b^2_{n+1}$ be a tuple 
realising  $\phi_i (\bar{x}) <\delta_i$,  
$\psi (\bar{x},x_{n+1})\le \tau$ and 
$\psi' (\bar{x},x_{n+1})\le \tau'$ such that 
$max_{i\le n+1} d(b^1_i ,b^2_i ) <\varepsilon_1 /2$. 
Note that  $max_{i\le n} d(a_i ,b^2_i ) <\varepsilon /2 +\varepsilon_1 /2$. 

Continuing this procedure we obtain a Cauchy sequence of 
$(n+1)$-tuples so that its limit satisfies $p$ 
and is not distant from $\bar{a}$ more than $\varepsilon$. 
$\Box$ 
\bigskip

We now apply this proposition in order to prove the following theorem. 

\begin{thm} \label{SCBorel}
Let $S=S_{{\bf Y}}$ be a dense countable subset of ${\bf Y}$. 
There is an $Iso (S)$-invariant Borel subset ${\cal SC}_S \subset {\bf Y}_L$ 
consisting of separably categorical continuous structures on $({\bf Y},d)$ 
such that any separably categorical continuous $L$-structure on ${\bf Y}$ 
is homeomorphic to a structure from $\mathcal{SC}_S$.  
\end{thm} 

{\em Proof.} 
Let ${\cal SC}_S$  be the set of all $L$-structures $M$ on ${\bf Y}$ 
so that for every $n$ and rational $\varepsilon$ there is
a finite set $F$ of tuples $\bar{a}_i$ from $S$ together 
with conditions $\phi_i (\bar{x})\le \delta_i$  ($i\in I$ and all $\delta_i$ are rational) 
with $\phi^{M}_i (\bar{a}_i )\le \delta_i$, $i\in I$, and the following property    
\begin{quote} 
any $n$-tuple $\bar{a}$ from $S$ satisfies in $M$ one of these  
$\phi_i (\bar{x})\le \delta_i$ and 

when $\phi^{M}_i (\bar{a})\le \delta_i$ and $\bar{c}$ is an $(n+1)$-tuple from $S$ 
with $c_1 ,...,c_n$ satisfying $\phi_i (\bar{x})\le \delta_i$ in $M$, for any finite set 
$\Delta$ of $L$-formulas $\phi (\bar{y})$, $|\bar{y}| =n+1$ with $\phi^{M} (\bar{c})=0$
there is an $(n+1)$-tuple $\bar{b}\in S$ so that 
$max_{j\le n} (d(a_j ,b_j ))\le \varepsilon /2$ and $\phi^{M}(\bar{b})=0$ 
for all formulas $\phi \in \Delta$. 
\end{quote}  
To see that $\mathcal{SC}_S$ is a Borel subsets of ${\bf Y}_L$  
it suffices to note that given rational $\varepsilon >0$, 
finitely many formulas $\phi_i (\bar{x})$, $i\in I$,  
with $|\bar{x}|=n+1$, and an $n$-tuple  $\bar{a}$  from $S$  
the set of $L$-structures $M$ on ${\bf Y}$ with the property  that 
\begin{quote} 
there is an $(n+1)$-tuple $\bar{b}\in S$ so that 
$max_{j\le n} (d(a_j ,b_j ))\le \varepsilon$ and 
$\phi^M_i (\bar{b})=0$ for all $i\in I$,  
\end{quote}  
is a Borel subset of ${\bf Y}_L$. 
The latter follows from Lemma \ref{EsLo}, which 
in particular says that any set of $L$-structures of the form 
$$
\{ M: M \models  max (max_{j\le n}(d(a_j ,b_j)\dot{-} \varepsilon ), max_{i\in I} (\phi_i (\bar{b}) ))=0 \} 
$$
is a Borel subset of ${\bf Y}_L$. 

Note that the density of the set of all tuples from $S$ in all ${\bf Y}^n$ 
implies that any continuous structure $M$ from ${\cal SC}_S$ satisfies the following property 
\begin{quote}  
for every $n$ and rational $\varepsilon$ there is a finite set $F_{n,\varepsilon}$ of $n$-tuples 
$\bar{a}_i$ from $S$ together with conditions $\phi_i (\bar{x})\le \delta_i$,  
$i\in I$, (all $\delta_i$ are rational) with $\phi^M_i (\bar{a}_i )\le \delta_i$, $i\in I$, so that     

any $n$-tuple $\bar{a}$ from ${\bf Y}$ satisfies one of these  $\phi_i (\bar{x})\le \delta_i$ and 

when $\phi_i (\bar{a})\le \delta_i$ and $\bar{c}$ is an $(n+1)$-tuple from ${\bf Y}$ 
with $c_1 ,...,c_n$ satisfying $\phi_i (\bar{x})\le \delta_i$, for any finite set 
$\Delta$ of $L$-formulas $\phi (\bar{y})$, $|\bar{y}|=n+1$, with $\phi (\bar{c})=0$
there is an $(n+1)$-tuple $\bar{b}\in S$ so that 
$max_{j\le n} (d(a_j ,b_j ))\le \varepsilon /2$ and $\phi(\bar{b})=0$ 
for all formulas $\phi \in \Delta$. 
\end{quote} 
It is now clear that any $M\in {\cal SC}_S$ satisfies the condition 
of Proposition \ref{SCHm}, i.e. $M$ is separably categorical. 

To finish the proof note that Proposition \ref{SCHm} 
also implies that if $M$ is a separably categorical structure on ${\bf Y}$, 
there is a dense set $S'\subseteq {\bf Y}$ so that $M$ belongs to 
the corresponding Borel set of $L$-structures ${\cal SC}_{S'}$.   
To see this we just extend $S_{{\bf Y}}$ to some countable $S'$ 
which satisfies the property of Proposition \ref{SCHm} 
in which we additionally require that $a_1 ,...,a_n \in S'$. 
It is clear any homeomorphism taking $S'$ onto $S$ 
takes $M$ into $\mathcal{SC}_S$ . 
$\Box$ 

\bigskip 

The proof above demonstrates that $\mathcal{SC}_S$ is of Borel level $\omega$.

\subsection{Complexity of sets of approximately ultrahomogeneous structures} 

Since any Polish group can be realised as 
the automorphism group of a approximately ultrahomogeneous 
structures it makes sense to consider the subset of ${\bf Y}_L$ 
of  all approximately ultrahomogeneous structures.  
Then one can try to characterise properties of Polish groups 
by description of the corresponding classes of approximately 
ultrahomogeneous structures and then to study the complexity 
of these classes. 
In this subsection we consider two opposite properties: 
separable oligomorphicity and admitting of complete left invariant metrics. 
Our results are not complete. 
For example we do not know if the class of approximately 
ultrahomogeneous $L$-structures on ${\bf Y}$ is 
a Borel subset of ${\bf Y}_L$. 

We use the following characterisation of approximate 
ultrahomogeneity from Section 6.1 of \cite{scho}.  
\begin{quote} 
Let $M$ be a separable continuous relational structure. 
Then $M$ is approximately ultrahomogeneous if and only if 
for any $\varepsilon$, for any quantifier free type 
$p(x_1 ,...,x_n , x_{n+1})$ realised in $M$ and any 
$a_1 ,...,a_n \in M$ realising the restriction of $p$ to 
$x_1 ,..., x_n$, there is a tuple $b_1 ,...,b_n ,b_{n+1}$ 
realising $p$ such that $max_{i\le n} d(a_i ,b_i ) <\varepsilon$.  
\end{quote}  

\begin{thm} \label{SCUBorel} 
Let $S=S_{{\bf Y}}$ be a dense countable subset of ${\bf Y}$. 
There is an $Iso (S)$-invariant Borel subset 
${\cal SCU}_S \subset {\bf Y}_L$ consisting of separably 
categorical approximatly ultrahomogeneous $L$-structures 
on $({\bf Y},d)$ such that any separably categorical approximately 
ultrahomogeneous $L$-structure on ${\bf Y}$ 
is homeomorphic to a structure from $\mathcal{SCU}_S$.  
\end{thm} 

{\em Proof.} 
We start with the claim that when $M$ is a non-compact, separable, continuous, metric 
structure on $({\bf Y},d)$ the following properties are equivalent: 

1. The structure $M$ is separably categorical and approximately ultrahomogeneous;  

2. For any $n$ and $\varepsilon$ there are finitely many  
quantifier free conditions $\phi_i (\bar{x})\le \delta_i$, $i\in I$,  
so that any $n$-tuple of $M$ satisfies one of these 
conditions and for any  $i\in I$, any $a_1 ,...,a_n \in M$ 
realising  $\phi_i (\bar{x})\le \delta_i$ and any quantifier free type 
$p(x_1 ,...,x_n , x_{n+1})$ realised in $M$ and containing  
$\phi_i (\bar{x})\le \delta_i$, there is a tuple 
$b_1 ,...,b_n ,b_{n+1}$ realising $p$ such that 
$max_{i\le n} d(a_i ,b_i ) <\varepsilon$.   

3. For any $n$ and $\varepsilon$ there are finitely many  
quantifier free conditions $\phi_i (\bar{x})\le \delta_i$, $i\in I$,  
so that any $n$-tuple of $M$ satisfies one of these 
conditions and the following property holds: 
\begin{quote} 
for any  $i\in I$, any $a_1 ,...,a_n \in M$ realising  
$\phi_i (\bar{x})\le \delta_i$ and  
any finite set of quantifier free formulas $\Delta (x_1 ,...,x_n , x_{n+1})$ 
realised in $M$ and containing  $\phi_i (\bar{x})\le \delta_i$, 
there is a tuple $b_1 ,...,b_n ,b_{n+1}$ 
realising $\Delta$ such that $max_{i\le n} d(a_i ,b_i ) <\varepsilon$.   
\end{quote} 
Indeed, to see the implication $1\rightarrow 2$ 
we apply the proof of Lemma \ref{SCH} as follows. 
Since the theory $T=Th(M)$ is separably categorical 
for each $n>0$, every $n$-type is principal and 
the metric space $(S_n (T),d)$ is compact.  
In particular for every $n$ and every $\varepsilon$ there is 
a finite family of principal $n$-types $p_1 ,...,p_m$ so that 
their $\varepsilon/2$-neighbourhoods cover $S_n(T)$. 

Note that by approximate ultrahomogenity for every $i\le m$, 
$$ 
dist(\bar{x}, p_i (M)) = dist(\bar{x}, p^{qe}_i (M)), 
$$ 
where $p^{qe}_i$ is the quantifier free part of $p_i$.  
Thus when $M$ is separably categorical, given $n$ and $\varepsilon$, 
we find appropriate $p_i$, $i\in I$, define    
$P_i (\bar{x})=dist(\bar{x}, p^{qe}_i (M))$, 
the corresponding definable predicates and  quantifier 
free $n$-conditions $\phi_i (\bar{x})\le \delta_i$ describing 
the corresponding $\varepsilon /2$-neighbourhoods of $p_i$.  
The rest follows by strong $\omega$-near-homogenity 
and approximate ultrahomogenity. 

Since by the reformulation of approximate ultrahomogenity above   
condition 2 obviously implies approximate ultrahomogenity, 
the implication $2\rightarrow 1$ follows easily from Lemma \ref{SCH}.  

The equivalence $2\leftrightarrow 3$ follows by arguments of  
the proof of  Proposition \ref{SCHm}.  

We now repeat the proof of Theorem \ref{SCBorel}  
to show that the following property defines a required Borel subset 
$\mathcal{SCU}_C$ of ${\bf Y}_L$: 
\begin{quote} 
for every $n$ and rational $\varepsilon$ there is a finite set $F$ of $n$-tuples 
$\bar{a}_i$ from $S$ together with quantifier free conditions 
$\phi_i (\bar{x})\le \delta_i$, $i\in I$ (all $\delta_i$ are rational),  
with $\phi_i (\bar{a}_i )\le \delta_i$, $i\in I$, so that 

any $n$-tuple $\bar{a}$ from $S$ satisfies one of these  $\phi_i (\bar{x})\le \delta_i$ and 

when $\phi_i (\bar{a})\le \delta_i$ and $\bar{c}$ is an $(n+1)$-tuple from $S$ 
with $c_1 ,...,c_n$ satisfying $\phi_i (\bar{x})\le \delta_i$, for any finite set 
$\Delta$ of quantifier free $L$-formulas $\phi (\bar{y})$, $|\bar{y}|=n+1$, with 
$\phi (\bar{c})=0$ there is an $(n+1)$-tuple $\bar{b}\in S$ so that 
$max_{j\le n} (d(a_j ,b_j ))\le \varepsilon /2$ and $\phi(\bar{b})=0$ 
for all formulas $\phi \in \Delta$. 
\end{quote}  
$\Box$ 
\bigskip

The following observation shows that the automorphism group of 
a non-compact separably categorical approximately ultahomogeneous 
structure does not admit a compatible complete left-invariant metric.  
Indeed, it it is an easy exercise that a non-compact 
separably categorical structure properly embeds ito itself. 
Thus condition (c) below supports our claim. 

\begin{prop} 
Let $G$ be the automorphism group of an approximately ultrahomogeneous 
continuous $L$-structure $M$ on the space $({\bf Y}, d)$.  
The following conditions are equivalent: \\ 
(a) the group $G$ admits a compatible complete left-invariant metric; \\ 
(b) the group $G$ is closed in $In({\bf Y},d)$; \\  
(c) there is no proper embedding of $M$ into itself. 
\end{prop} 

{\em Proof.} 
Let $\rho_S$ be the standard metric of $Iso ({\bf Y})$ 
( $S=S_{{\bf Y}}$, see Introduction). 
To see that (b) implies (a) note that closedness of $G$ 
in $In({\bf Y},d)$ guarantees that any Cauchy $\rho_S$-sequence 
of elements from $G$ has a limit in $G$, i.e. $\rho_S$ 
is a compatible complete left-invariant metric of $G$. 

For the converse note that by Lemma 2.1 of \cite{gao} when 
$G$ has a compatible complete left-invariant metric, 
any compatible left-invariant metric is complete. 
Thus the metric $\rho_S$ is complete and $G$ is closed in $In ({\bf Y},d)$. 

To see that (b) implies (c) assume that there is a proper 
embedding of $M$ into itself (say $h$). 
Then for each sequence $s_1 ,...,s_n \in S$ the quantifier 
free types of this sequence coincides with the quantifier 
free type of $h(s_1 ), ...,h(s_n )$. 
Since $M$ is approximately ultrahomogeneous for any $\varepsilon$ 
there is an automorphism $g_{n,\varepsilon}$ taking every $s_i$ to 
the $\varepsilon$-ball of $h(s_i )$. 
This produces a Cauchy $\rho_S$-sequence from $G$ with the limit $h$, 
contradicting the closedness of $G$ in $In({\bf Y},d)$. 

The negation of (b) implies the negation of (c) by an obvious reason. 
$\Box$ 

\bigskip 

What is the complexity of the class of approximately 
ultrahomogeneous structures from this proposition? 
The following statement gives a partial answer. 
It somehow corresponds to the result 
of M.Malicki \cite{mal} that the set of all Polish groups admitting 
compatible complete left-invariant metrics is coanalytic 
non-Borel as a subset of a standard Borel space of Polish groups. 

\begin{cor} \label{corclim} 
The subset of ${\bf Y}_L$ consisting of approximately 
ultrahomogeneous structures $M$ such that $Aut(M)$ admits 
a compatible complete left-invariant metric, is coanalytic 
in any Borel subset of approximately ultrahomogeneous 
structures. 
\end{cor} 

{\em Proof.}  
It is enough to show that the subset of $L$-structures $M$ admitting 
proper embeddings into itself is analytic. 
To see this consider the extension of $L$ by a unary function $f$.  
All expansions of $L$-structures satisfying the property 
that $f$ is an isometry which preserves $L$-relations, form a closed 
subset of the (Polish) space of all $L\cup \{ f\}$-structures. 

If $s\in S=S_{{\bf Y}}$ and $\varepsilon \in {\bf Q}\cap [0,1]$ then 
the condition that $f(S)$ does not intersect the 
$\varepsilon$-ball of $s$ is open. 
Thus the set of $L\cup \{ f\}$-structures with a proper 
embedding $f$ into itself, is Borel. 
The rest is easy. 
$\Box$ 


\subsection{Computable presentations}

Consider the situation of Section 2.2. 
Let $G$ be a Polish group and $\mathcal{R}$ 
be a distinguished countable family of clopen 
grey cosets so that the family 
$\{ \rho_{<q} : \rho\in \mathcal{R}$ and $q\in \mathbb{Q}^{+}\cap [0,1]\}$ 
forms a basis of the topology of $G$. 
Let $G_0$ be a countable dense subgroup of $G$ so that 
$\mathcal{R}$ is closed under $G_0$-conjugacy and consists of all  
$G_0$-cosets of grey subgroups from $\mathcal{R}$.  
We may also assume that the set of grey subgroups from 
$\mathcal{R}$ is closed under $max$ and truncated multiplication 
by positive rational numbers.
If we assume that there is a 1-1-enumeration of the family 
$$ 
\{ \rho_{<q} : \rho\in \mathcal{R}\mbox{ and }q\in \mathbb{Q}^{+}\cap [0,1]\} \cup 
\{ \rho_{>q} : \rho\in \mathcal{R} \mbox{ and }q\in \mathbb{Q}^{+}\cap [0,1]\} 
$$
so that the relation of inclusion between members 
of this family is computable we arive at the case 
that $G$ is a {\em computably presented $\omega$-continuous domain},  
\cite{EH} and \cite{ES}. 

When we consider a $(G,G_0 ,\mathcal{R})$-space $({\bf X},\tau )$ 
together with a distinguished countable $G_0$-invariant 
grey basis $\mathcal{U}$ (see Definition \ref{gbasis}) 
of clopen grey subsets so that the relation of inclusion between 
sets of the form 
$$ 
\sigma_{<r}\mbox{ or }\sigma_{>r} \mbox{ for }
\sigma \in \mathcal{U}\mbox{ and }r\in \mathbb{Q}^+ \cap [0,1], 
$$ 
is computable (under an appropriate coding) we also 
obtain a computably presented $\omega$-continuous domain.
We denote it by $\mathcal{U}_{\mathbb{Q}}$. 
Note that in the discrete case these circumstances are 
standard and in particular arise when one studies 
computability in $S_{\infty}$-spaces of logic actions.  

Let 
$$ 
\mathcal{U}^+_{\mathbb{Q}} =\{ \sigma_{<r} :\sigma \in \mathcal{U}\mbox{ and }r\in \mathbb{Q}^+ \cap [0,1]\} , 
$$ 
$$ 
\mathcal{R}^+_{\mathbb{Q}}=\{ \rho_{<q} : \rho\in \mathcal{R}\mbox{ and }q\in \mathbb{Q}^{+}\cap [0,1]\} . 
$$ 
Below we will restrict ourselves by only $+$-parts of domains above. 

\begin{remark} 
{\em It is worth noting that when we have a recursively 
presented Polish space in the sense of the book of 
Moschovakis \cite{mos} (Section 3), 
then a basis of the form $\mathcal{U}_{\mathbb{Q}}$ 
as above (in fact $\mathcal{U}^+_{\mathbb{Q}}$) 
can be naturally defined. 
Indeed, let us recall that a recursive presentation of a Polish space 
$({\bf X},d)$ is any sequence $S_{{\bf X}} =\{ x_i : i\in\omega\}$ 
which is a dense subset of ${\bf X}$ satisfying 
the condition that $(i,j,m,k)$-relations 
$$
d(x_i ,x_j )\le \frac{m}{k+1} \mbox{ and } d(x_i ,x_j ) < \frac{m}{k+1} 
$$ 
are recursive. 
If in this case for all $i$ we define grey subsets 
$\phi_i (x) = d(x,x_i )$,  then all balls 
$(\phi_i )_{<r}$, $r\in \mathbb{Q}$, form a basis 
$\mathcal{U}^+_{\mathbb{Q}}= \{ B_i :i\in\omega \}$ of ${\bf X}$ 
which under appropriate enumeration 
(together with co-balls $(\phi_i )_{>r}$) 
satisfies our requirements above. 
When $G$ is a Polish group with a left-invariant 
metric $d$, then for any $q_1 ,...,q_k\in \mathbb{Q}$ 
and any tuple $h_1 ,...,h_k \in G$ the grey subset 
$\phi_{\bar{q},\bar{h}} (x)= max_{i\le k} (q_i \cdot d(h_i , xh_i ))$ 
is a grey subgroup 
\footnote{apply $d(h_i ,xyh_i ) \le d(h_i ,x h_i ) + d(xh_i ,xy h_i )= d(h_i ,xh_i )+d(h_i ,yh_i )$ 
together with the fact that $max$ applied to grey subgroups gives grey subgroups again} 
. 
If $G$ is a recursively presented space 
with respect to a dense countable subgroup 
$G_0$ and the multiplication is recursive, 
then let $\mathcal{V}$ consist of all $\phi_{\bar{q},\bar{h}}$ 
with $\bar{h} \in G_0$ and let $\mathcal{R}$ consist 
of all $G_0$-cosets of these grey subgroups. 
The structure (domain) $(\mathcal{R}_{\mathbb{Q}}, \subset )$ 
is computably presented. 

If $G$ isometrically acts on ${\bf X}$ and $x_1 ,...,x_k$ 
is a finite subset of the recursive presentation $S_{{\bf X}}$ 
then as we already know the function 
$\psi _{\bar{q},\bar{x}}(g)=max_{i\le k} (q_i \cdot d(x_i ,g(x_i )))$ 
also defines a grey subgroup. 
When $G$ has a recursive multiplication and a recursive action 
on ${\bf X}$ (see Section 3 of \cite{mos}) so that 
the recursive presentation $S_{{\bf X}}$ is $G_0$-invariant, 
then let $\mathcal{V}$ consist of these subgroups 
and $\mathcal{R}$ consist of the $G_0$-cosets.  
Then the structure $(\mathcal{R}_{\mathbb{Q}}, \subset )$ 
is computably presented. } 
\end{remark} 

We always assume that $\mathcal{U}$ and 
the subfamily of all grey subgroups 
from $\mathcal{R}$ are closed under $max$ 
and truncated multiplying by positive rational numbers.  
Families $\mathcal{U}^+_{\mathbb{Q}}$ and $\mathcal{R}^+_{\mathbb{Q}}$ 
are bases of the corresponding spaces. 

We now take some computability assumptions. 
As we will see below they are satisfied in 
the majority of interesting cases.  

{\bf A1.} We assume that the sets of indices (under our enumeration) 
of $\mathcal{U}^+_{\mathbb{Q}}$, $\mathcal{R}^+_{\mathbb{Q}}$  
and the set of rational cones $H_{<r}$, $H_{>r}$ for basic 
{\em subgroups} from $\mathcal{R}$ are distinguished 
by computable relations on $\omega$. 
We denote the latter by $\mathcal{V}_{\mathbb{Q}}$. 

{\bf A2.} We assume that under our 1-1-enumerations of the families 
$\mathcal{R}_{\mathbb{Q}}$ and $\mathcal{U}_{\mathbb{Q}}$ 
the binary relation to be in the pair $\sigma_{<r}$, $\sigma_{>r}$ 
for $\sigma \in \mathcal{R}$ or $\sigma\in \mathcal{U}$ 
is computable. 

{\bf A3.} We also assume that the following relation is computable: 
\begin{quote} 
{\it 
$Inv(V,U)\Leftrightarrow (V\in \mathcal{V}^+_{\mathbb{Q}})\wedge (U\in \mathcal{U}^+_{\mathbb{Q}})\wedge (U$ 
is $V$-invariant $)$ }  \\ 
(recall the latter means that when $\phi (x)<r$ and $H(g)<s$ we have 
$\phi (g(x)) < r+s$).
\end{quote} 

{\bf A4.}  We assume that there is an algorithm deciding 
the problem whether for a natural number $i$ and for a basic set 
of the form $\sigma_{<r}$ for $\sigma$ 
from $\mathcal{U}$ or $\mathcal{R}$ and $r\in \mathbb{Q}$,  
the diametr of $\sigma_{<r}$  is less than $2^{-i}$.  

\bigskip 

In \cite{mos} it is defined that a point $x\in {\bf X}$ 
is {\em recursive} if the set 
$\{ s :x\in B_s , B_s \in \mathcal{U}_{\mathbb{Q}}\}$ 
is computable. 
We immitate it in the following definition.

\begin{definicja} \label{Comp} 
We say that an element $x\in {\bf X}$ is {\bf computable} if the relation 
$$
Sat_x (U)\Leftrightarrow (U \in \mathcal{U}_{\mathbb{Q}})\wedge (x\in U)
$$ 
is decidable. 
\end{definicja}

In the case of the logic action of $S_{\infty}$, 
when $x$ is a structure on $\omega$ and all $H$ and $\phi$ are 
two-valued, this notion is obviously equivalent to 
the notion of a computable structure.  

We will denote by $Sat_x ({\mathcal{U}}_{\mathbb{Q}})$ the set 
$\{ C: C\in {\mathcal{U}}_{\mathbb{Q}}$ and $Sat_x (C)$ holds $\}$. \\ 

The following lemma follows  from the assumption 
that $\mathcal{U}$ is a grey basis and satisfies {\bf A4}.   

\begin{lem} \label{kappa} 
If $x\in {\bf X}$ is computable then there is a computable function 
$\kappa :\omega \rightarrow \mathcal{U}^+_{\mathbb{Q}}$ such that for all 
natural numbers $n$, $x\in \kappa (n)$ and $diam(\kappa (n))\le 2^{-n}$. 
\end{lem} 
 
We also say that an element $g\in G$ is {\bf computable} 
if the relation $(N\in \mathcal{R}_{\mathbb{Q}})\wedge (g\in N)$ 
is computable. 
Then there is a computable function realising the 
same property as $\kappa$ above but already in the case 
of the basis $\mathcal{R}_{\mathbb{Q}}$. 

In the following lemma we use standard indexations 
of the set of computable functions and of the set 
of all finite subsets of $\omega$. 

\begin{lem} \label{21}
The following relations belong to $\Pi^0_2$:\\
(1) $\{ e:$ the function $\varphi_e$ is a characteristic 
function of a subset of $\mathcal{U}_{\mathbb{Q}} \}$; \\ 
(2) $\{ (e,e'):$ there is a computable element $x\in {\bf X}$ 
such that the function $\varphi_e$ is a characteristic function 
of the set $Sat_x ({\mathcal{U}_{\mathbb{Q}}})$ and 
the function $\varphi_{e'}$ realizes the corresponding function 
$\kappa$ defined in Lemma \ref{kappa} $\}$;\\   
(3) $\{ (e,e'):$ there is an element $g\in G$ such that the 
function $\varphi_e$ is a characteristic function of the 
subset $\{ N\in\mathcal{R}_{\mathbb{Q}}:g\in N\}$ and the function 
$\varphi_{e'}$ realises the corresponding function $\kappa$ 
defined as in Lemma \ref{kappa} $\}$.    
\end{lem} 

{\em Proof.} (1) Obviuos. 
Here and below we use the fact that a function is computable if 
and only if its graph is computably enumerable. 

(2) The corresponding definition can be described as follows: 
$$
("e \mbox{ is a characteristic function of a subset of }\mathcal{U}_{\mathbb{Q}}")\wedge 
$$ 
$$
(\forall n)((\varphi_{e'}(n)\in\mathcal{U}^+_{\mathbb{Q}})\wedge (\varphi_{e'}(n)\not=\emptyset )
\wedge (\varphi_{e}(\varphi_{e'}(n))=1 )\wedge 
(diam(\varphi_{e'}(n))<2^{-n} )) \wedge  
$$
\begin{quote} 
$(\forall d)(\exists n )(($ 
"every element $U'$ of the finite subset of $\mathcal{U}_{\mathbb{Q}}$ 
with the canonical index $d$ satisfies $\varphi_e (U') =1$") 
$\leftrightarrow ($ "$\varphi_{e'}(n)$ is contained in any element $U'$ 
of the finite subset of $\mathcal{U}_{\mathbb{Q}}$ with the canonical index $d$"$))$.  
\end{quote}
The last part of the conjunction ensures that the intersection of 
any finite subfamily of $\mathcal{U}_{\mathbb{Q}}$ of cones $U'$ 
with $\varphi_e (U') =1$ contains a closed cone of the form $\phi_{\le r}$ 
of sufficiently small diameter. 
Now the existence of the corresponding $x$ follows by Cantor's theorem. 

(3) is similar to (2). 
$\Box$
\bigskip

As a result we see that the set of indices of computable 
elements of ${\bf X}$ belongs to $\Sigma^0_3$. 
If ${\bf X}$ is of the form ${\bf Y}_L$ then it makes sense 
to study complexity of sets of indices of computable structures 
of natural model-theoretic classes. 
In the case of first order structures this approach is 
traditional, see \cite{AK}, \cite{GK}. 
We now illustrate it in the case of $\mathbb{U}_L$ 
for relational $L$ (for simplicity). 
In fact we give effective versions of 
results of Sections 5.1 and 5.2. 
We start with the following theorem. 
We will not directly apply its statement, 
but the arguments of the proof will be very helpful 
below. 

\begin{thm} 
The structure $(\mathbb{U},  s)_{ s\in\mathbb{QU}}$ of 
the expansion of the Urysohn space  by constants from 
$\mathbb{QU}$ has a decidable continuous theory: 
for every continuous sentence of the form 
$\phi (\bar{s})$  where  $\bar{s} \in \mathbb{QU}$ 
the value of $\phi (\bar{s})$ in $\mathbb{U}$ 
is a computable real number. 
\end{thm} 

{\em Proof.} 
We remind the reader that a real number $r\ge 0$ is computable 
if there is an algorithm which for any natural number $n$ 
finds a natural number $k$ such that 
$$ 
 \frac{k-1}{n} \le r \le \frac{k+1}{n}.  
$$ 
To prove the theorem we use the main result of \cite{BYP}.  
Let $T_{\mathbb{QU}}$ consist of the standard axioms of $\mathbb{U}$ 
(with rational $\varepsilon$ and $\delta$, see Section 5 in \cite{U}) 
together with all quantifier free axioms 
describing distances between constants from $\mathbb{QU}$. 
We claim that $T_{\mathbb{QU}}$ is decidable. 
Since the set of all standard axioms of $\mathbb{U}$ 
is decidable (see \cite{U}), it suffices to check that 
the elementary (not continuous) theory of the structure 
$\mathbb{QU}$ in the language of binary relations 
$$ 
d(x,y) = q \mbox{ , where } q\in \mathbb{Q}\cap [0,1], 
$$ 
is decidable. 
The latter is straightforward. 

Note that $T_{\mathbb{QU}}$ is complete, i.e. it axiomatises 
the continuous theory of some continuous structure.  
Indeed, otherwise there is a separable continuous structure 
$M\models T_{\mathbb{QU}}$ such that for some  
tuple $\bar{s}\in  \mathbb{QU}$ the structures $(\mathbb{U},\bar{s})$ 
and the reduct of $M$, say $M'$, to the signature $(d,\bar{s})$, 
do not satisfy the same inequalities of  the form 
$$
\phi (\bar{s})\le (<) q \mbox{ or }  \phi (\bar{s})\ge (>) q \mbox{ where } 
q \in \mathbb{Q}\cap [0,1].  
$$  
On the other hand since $\mathbb{U}$ is separably  categorical 
and ultrahomogeneous, the structures $M'$ and  $(\mathbb{U},\bar{s})$ 
are isomorphic, contradicting the previous sentence.  

By Corollaries 9.8 and  9.11 of \cite{BYP} there is an algorithm which for every 
continuous sentence $\phi(\bar{s})$ computes its value in $\mathbb{U}$. 
$\Box$

\bigskip 

Let $L$ be a relational language. 
Let us consider the space $\mathbb{U}_L$, the family of 
grey cosets $\mathcal{R}^{U}$ and the nice basis 
$\mathcal{B}_0$ defined in Section 3. 
The subfamily $\mathcal{B}_{qf}$ of grey subsets 
from $\mathcal{B}_0$ correspoding to quantifier free 
$L$-formulas is considered as the the grey basis 
$\mathcal{U}$ above.  

To check that the $Iso (\mathbb{U})$-space 
$\mathbb{U}_L$ satisfies the computability conditions above 
(in particular {\bf A1} - {\bf A4}), 
note that $\mathbb{QU}$ under the language of binary relations 
$$ 
d(x,y) = q \mbox{ , where } q\in \mathbb{Q}\cap [0,1], 
$$ 
has a presentation on $\omega$ 
so that all relations first-order definable 
in $\mathbb{QU}$, are decidable. 
This follows from the fact that 
the elementary theory of the structure 
$\mathbb{QU}$ in the language expanded by 
all constants has quantifier elimination and 
is computably axiomatisable (i.e. 
the corresponding theory is decidable). 
We fix such a presentation. 

Then we can code a $q'$-cone of a grey coset  
$$ 
I_{q, \bar{s},\bar{s}'} : g\rightarrow q \cdot d(g(\bar{s}'), \bar{s}), 
\mbox{ where } \bar{s}\bar{s}'\subset \mathbb{QU}, \mbox{ and } q\in \mathbb{Q}^{+}.   
$$ 
(i.e. the coset $H_{q ,\bar{s}} g_0$ of the grey subgroup 
$$ 
H_{q, \bar{s}} : g\rightarrow q \cdot d(g(\bar{s}), \bar{s}), 
\mbox{ where } \bar{s}\subset \mathbb{QU}, \mbox{ and } q\in \mathbb{Q}^{+}.   
$$ 
\centerline{ with respect to $g_0$ taking $\bar{s}'$ to $\bar{s}$)} \\ 
by the number of the tuple $(q, \bar{s},\bar{s}',q',*)$, 
where $*$ corresponds to one of the symbols $<, \le ,>, \ge$.  
It is known that for any $\bar{t}\in \mathbb{U}$  
the algebraic closure of $\bar{t}$ in $\mathbb{U}$ 
coincides with $\bar{t}$ (Fact 5.3 of \cite{ealygold}). 
Now using decidability of the elementary 
diagram of $\mathbb{QU}$ we see that  
the relation of inclusion between cones 
of this form is decidable. 
Cones of grey subgroups 
(i.e. the set $\mathcal{V}_{\mathbb{Q}}$)
are distinguished by the computable 
subset of tuples as above 
with $\bar{s}=\bar{s'}$. 
\parskip0pt 

Since we interpret elements of $\mathcal{B}_{0}$ by 
$L$-formulas with parameters from $\mathbb{QU}$ 
and without free variables,
it is obvious that both $\mathcal{B}_{0}$ and $\mathcal{B}_{qf}$ 
can be coded in $\omega$ so that the operations of 
connectives are defined by computable functions. 
Moreover $\mathcal{B}_{qf}$ is a decidable subset of $\mathcal{B}_0$. 
Now all cones of the form 
$\sigma_{<q}$, $\sigma_{>q}$, $\sigma_{\le q}$, $\sigma_{\ge q}$ 
can be enumerated so that all natural relations between them 
(in particular relations from {\bf A2}) are computable. 
To satisfy {\bf A3} we define $Inv (V,U)$ as follows. \parskip0pt

$Inv(V,U)\Leftrightarrow$ 
"$U$ is of the form $\sigma_{<l}$ for 
$\sigma\in \mathcal{B}_{qf}$, $V$ 
is of the form $H_{<k}$ for $H\in \mathcal{V}$ 
and the tuple of parameters of 
$\sigma_{<l}$ is contained in the tuple 
of elements of $\mathbb{QU}$ which 
appears in the code of $H_{<k}$".  \\
This relations is obviously decidable. 
 
Let $\phi(\bar{s})$ be a quantifier-free formula 
defining an element $A\in \mathcal{B}_{qf}$.  
To compute $diam(A)$ consider the definition of 
the metric $\delta_{seq(\mathbb{QU})}$ of the space 
$\mathbb{U}_{L}$ with respect to $sec(\mathbb{QU})$ 
in the beginning of Section 1.1. 
Assuming (for simplicity) that $\phi$ is a conjunction 
of atomic inequalities find all numbers $i$ of 
tuples $(j,\bar{s}')$ such that 
$R_j (\bar{s}')$ appears in $\phi (\bar{s})$. 
We may assume that appearance of such subformulas  
forces inequalities of the form 
$q'_i \le R_j (\bar{s}') \le q_i$ for rational 
$0\le q'_i ,q_i \le 1$. 
Let $I$ be the (finite) subset of such $i$. 
Then $diam(A)$ is computed by  
$$ 
\sum_{i=1}^{\infty} \{ 2^{-i}\mbox{ : }  i \not\in I\} + 
\sum_{i\in I} 2^{-i} |q_i - q'_i| . 
$$  
The case of basic clopen sets of $\mathcal{R}^U$ is similar. 

The following theorem is an effective version 
of main results of Sections 5.1 and 5.2. 

\begin{thm}  
Let  ${\cal SC}_{\mathbb{QU}}$  be the 
$Iso (\mathbb{QU})$-invariant Borel subset of $\mathbb{U}_L$ 
defined as in Theorem \ref{SCBorel} 
(in particular, consisting of separably categorical continuous structures on $(\mathbb{U},d)$ 
such that any separably categorical continuous $L$-structure on $\mathbb{U}$ 
is homeomorphic to a structure from $\mathcal{SC}_{\mathbb{QU}}$).  

Let ${\cal SCU}_{\mathbb{QU}}$ be the 
$Iso (\mathbb{QU})$-invariant Borel subset of $\mathbb{U}_L$ 
defined as in Theorem \ref{SCUBorel} 
(in particular consisting of separably 
categorical approximatly ultrahomogeneous $L$-structures 
on $(\mathbb{U},d)$ such that any separably categorical approximately 
ultrahomogeneous $L$-structure on $\mathbb{U}$ 
is homeomorphic to a structure from $\mathcal{SCU}_{\mathbb{QU}}$).  

Then the subsets of indices of computable structures from 
${\cal SC}_{\mathbb{QU}}$ and 
${\cal SCU}_{\mathbb{QU}}$ respectively are hyperarithmetical. 
\end{thm} 

{\em Proof.} 
We use the following observation. 
\begin{quote} 
The set of all pairs $(i,j)$ where $j$ is an index of a cone 
from $(\mathcal{B}_0 )_{\mathbb{Q}}$ and $i$ is an index 
of a computable structure from this cone, is hyperarithmetical 
of level $\omega$. 
\end{quote}   
This is an effective version of Proposition \ref{EsLo}. 
It follows from Lemma \ref{21}  by standard arguments. 
Note that as we have shown above 
(using decidability of the elementary diagram 
of $(\mathbb{QU}, s)_{s\in \mathbb{QU}}$) all assumptions 
of Lemma \ref{21} are satisfied under the circumstances of 
our theorem. 

It remains to verify that definitions of sets ${\cal SC}_{\mathbb{QU}}$ 
and ${\cal SCU}_{\mathbb{QU}}$ from the proofs of 
Theorems \ref{SCBorel} and \ref{SCUBorel} define hyperarithmetic 
subsets of indices of computable structures. 
This is straghtforward. 
$\Box$

\section{Appendix: An amalgamation property of the Urysohn space}

\begin{quote}
{\bf Abstract.} 
We present some amalgamation property of the Urysohn sphere $\mathbb{U}$. 
As a consequence we evaluate the distance between types of $Th(\mathbb{U})$ 
with parameters. 
\bigskip

{\em Keywords:}  finite metric spaces, amalgamation, Urysohn space, continuous model theory

{\em 2010 Mathematics Subject Classification:} 03C35, 54A05, 54D80

\end{quote}

\bigskip

\subsection{Introduction}

The {\bf Urysohn space of diameter 1}, which is denoted by  $\mathbb{U}$, 
is the unique Polish metric space which is universal for spaces of diameter 1 and 
{\bf ultrahomogeneous}, i.e. every isometry between finite subsets of 
$\mathbb{U}$ extends to an isometry of $\mathbb{U}$. 
It is sometimes called the {\bf Urysohn sphere}. 
As a continuous structure the space $\mathbb{U}$  is considered in 
the continuous signature  $\langle d \rangle$, where $d$ is the metric, 
see \cite{BYBHU}. 
\parskip0pt 

The countable counterpart of $\mathbb{U}$ is the 
{\bf rational Urysohn space of diameter 1} $\mathbb{QU}$, 
which is both ultrahomogeneous and universal for countable 
metric spaces with rational distances and diameter $\le 1$. 
The space $\mathbb{QU}$ is usually considered as the first-order 
structure of infinitely many binary relations 
$$ 
d(x,y) \le q \mbox{ , where } q\in \mathbb{Q} \cap [0,1] . 
$$ 
It is shown in Section 5.2 of \cite{BYBM} that there is 
an embedding of $\mathbb{QU}$ into $\mathbb{U}$ so that: \\
(i) $\mathbb{QU}$ is dense in $\mathbb{U}$; \\ 
(ii) any isometry of  $\mathbb{QU}$ extends to an isometry of  
$\mathbb{U}$ and $Iso (\mathbb{QU})$ is dense in $Iso(\mathbb{U})$;  \\ 
(iii) for any $\varepsilon>0$, any partial isometry $h$ of  
$\mathbb{QU}$ with domain $\{ a_1 ,...,a_n\}$ and any isometry 
$g$ of $\mathbb{U}$ such that $d(g(a_i ),h(a_i ))<\varepsilon$ 
for all $i$, there is an isometry $\hat{h}$ of  $\mathbb{QU}$ 
that extends $h$ and is such that for all  $x\in \mathbb{U}$, 
$d(\hat{h}(x),g(x))<\varepsilon$.

We recommend \cite{melleray-Ur}, \cite{N} and \cite{U} for basic 
information concerning $\mathbb{U}$.   
We will use below the free amalgamation 
property for finite metric spaces, see Theorem 2.1 of \cite{bogaty}. 
\begin{quote} 
Assume that $(X,d_1 )$ and $(Y,d_2 )$ are finite metric spaces 
with $Z =X\cap Y$ and $d_1 =d_2$ for elements of $Z$. 
Then there is a metric $d$ on $X\cup Y$ which agrees with $d_1$ on $X$, 
with $d_2$ on $Y$ and is defined for $x\in X\setminus Y$ 
and $y\in Y\setminus X$ by 
$$ 
d(x,y ) = min_{z\in Z}(d_1 (x,z )+d_2 (y,z)) 
\mbox{ when } Z\not=\emptyset ; 
$$ 
$$
d(x,y) = d_1(x ,x^* ) + d(x^* ,y^* ) + d_2 (y^* ,y) 
\mbox{ when } Z = \emptyset \mbox{ and } 
$$ 
$$ 
x^* \in X \mbox{ , } y^* \in Y 
\mbox{ are distinguished together with the distance between them. }
$$
\end{quote} 

The main result of this appendix is some amalgamation property 
of the Urysohn space. 
It roughly states that if $A$ and $B$ are finite subspaces of $\mathbb{U}$  
which are sufficiently similar then $B$ has a  
copy $B'$ in $\mathbb{U}$ under an isometry fixing 
$A\cap B$ so that $A$ and $B'$ are sufficiently 
close in $\mathbb{U}$. 
It is presented in Section 6.2.  
In Section 6.3 we apply it in order  to evaluate the distance 
between types in $Th (\mathbb{U})$ with parameters.  
We think that some versions of the amalgamation 
property presented in Section 6.2 are already folklore. 
On the other hand the only published material 
we have found is the parameter-free case 
(i.e. when $A\cap B =\emptyset$)  
considered in \cite{CL}.  
The construction presented in \cite{CL} is based on 
Example 56 from \cite{ppetersen}, p. 295 - 296.  
Although the general case looks slightly more 
complicated than the case from \cite{CL},  
we have found that the idea of Example 56 of \cite{ppetersen} 
works in the general case too. 
This  improves the corresponding result 
from \cite{IMI} and simplifies the proof. 

Section 6.2 does not use any special material. 
All model theoretic notions used in Section 6.3 
will be defined there.

\subsection{Amalgamation}

The following theorem states the property that 
if $A$ and $B$ are finite metric spaces  
which are sufficiently similar then there is metric space $C$ 
containing $A$ and  a  copy of $B$ under an isometry fixing 
$A\cap B$ so that $A$ and the image of $B$ are sufficiently 
close in $C$. 
We do not have any assumptions on diameters of spaces. 

\begin{thm}  \label{LemUrysohn}
Let a finite metric space $A = \{ a_1 ,...,a_n \}$ and numbers $0\le q <n$ and 
$\varepsilon >0$ satisfy all inequalities of the form 
$$ 
4\varepsilon < d(a_i ,a_j ) \mbox{ for pairs } i<j\le n \mbox{ with } q<j \mbox{ and } 
$$ 
$$  
4\varepsilon < d(a_i ,a_j ) + d(a_i ,a_k ) - d(a_j ,a_k ) 
\mbox{ for triples } a_i ,a_j ,a_k \mbox{ with } |\{ i,j,k\}| =3
$$ 
$$ 
\mbox{ and }  k\le q < min (i,j) <n   . 
$$  

Let $B$ be an $n$-element metric space consisting 
of elements $b_i$ so that for each pair $i<j\le n$, 
$|d(b_i ,b_j )-d(a_i ,a_j )|\le \varepsilon$. 
We assume that $a_1 =b_1$,...,$a_q =b_q$, $A\cap B =\{ a_1 ,...,a_q \}$, 
and the metric defined on $\{ b_1 ,...,b_q\}$ in the space $B$ 
coincides with the metric defined on $\{ a_1 ,...,a_q\}$ in $A$.   

Then there is a metric on $A \cup B$ extending metrics in $A$ and $B$ 
so that for each $q<i\le n$, $d(a_i ,b_i )=\varepsilon$. 
\end{thm} 

{\em Proof.} 
The case $A\cap B=\emptyset$ is already considered in 
Example 56 \cite{ppetersen}, p. 295 - 296. 
We will assume that $A\cap B\not=\emptyset$. 
Let us build a metric space on $\{ a_1 ,...,a_n , b_{q+1} ,...,b_n \}$ 
so that $d(a_i ,b_i )=\varepsilon$ for all $q<i\le n$. 
We define $d(a_i ,b_j )$ for $i\not=j$ with $i>q$ or $j>q$ as follows: 
$$ 
d(a_i ,b_j ) = min (\{ d(a_i ,a_k )+ d(a_k ,b_j ): k\le q\} \cup 
\{ d(a_i ,a_k ) +\varepsilon +  d(b_k ,b_j ) : k>q \} ) . 
$$
Below we will use the observation that in the cases 
$i\le q<j$ or $j\le q<i$  the metric on $B$ or $A$ respectively 
which is given in the formulation of the theorem,  
satisfies the condition of this definition. 
To see this one should note that in these cases 
$$  
d(a_i ,b_j ) \le  d(a_i ,a_k ) +\varepsilon +  d(b_k ,b_j )  \mbox{ when } k>q . 
$$ 
Indeed, if  for example $i\le q<j$ then 
$$
d(a_i ,b_j )\le d(b_i ,b_k ) + d(b_k ,b_j ) \le d(a_i ,a_k ) +\varepsilon +  d(b_k ,b_j ) . 
$$ 
Let us verify the triangle inequality. 
We may restrict ourselves by triangles which intersect both 
$A\setminus B$ and $B\setminus A$.  
We will use below the following consequence of the 
assumptions of the theorem: 
$$ 
3\varepsilon < d(b_i ,b_j ) \mbox{ for pairs } i<j\le n \mbox{ with } q<j \mbox{ and } 
$$ 
$$  
\varepsilon < d(b_i ,b_j ) + d(b_i ,b_k ) - d(b_j ,b_k ) 
\mbox{ for triples } b_i ,b_j ,b_k \mbox{ with } |\{ i,j,k\}| =3
$$ 
$$ 
\mbox{ and }  k\le q < min (i,j) <n   . 
$$

Case 1. $d(a_i ,b_j ) \le d(a_i ,a_l )+ d(a_l ,b_j )$. 
By the assumptions of the theorem we may assume that $i\not= j$. 
If $d(a_l ,b_j )= d(a_l ,a_k )+ d(a_k ,b_j )$ with  $k\le q$, then  
$$ 
d(a_i ,a_l )+ d(a_l ,b_j ) = d(a_i , a_l ) + d(a_k ,a_l )+ d(a_k ,b_j ) \ge 
d(a_i ,a_k ) + d(a_k ,b_j ) \ge d(a_i , b_j ) . 
$$ 
If $d(a_l ,b_j ) = d(a_l ,a_k ) +\varepsilon +  d(b_k ,b_j )$ for some $k>q$, 
then 
$$ 
d(a_i ,a_l )+ d(a_l ,b_j ) = d(a_i , a_l ) + d(a_k ,a_l )+ \varepsilon + d(b_k ,b_j ) \ge 
d(a_i ,a_k ) + d(b_k ,b_j ) +\varepsilon \ge d(a_i , b_j ) . 
$$ 

Case 2.  $d(a_i ,b_j ) \le d(a_i ,b_l )+ d(b_l ,b_j )$. 
This case is similar to Case 1. 

Case 3.   $d(a_i ,a_j ) \le d(a_i ,b_l )+ d(a_j ,b_l )$.
We may assume that  $i>q$ or  $j>q$.   
If 
$$
d(a_i ,b_l )= d(a_i ,a_k )+ d(a_k ,b_l ) \mbox{ and } d(a_j ,b_l )= d(a_j ,a_m )+ d(a_m ,b_l ) 
\mbox{ with  }k,m \le q, 
$$ 
then
$$ 
d(a_i ,b_l ) + d(a_j ,b_l ) \ge  d(a_i ,a_k )+ d(a_k ,b_l ) + d(a_j ,a_m )+ d(a_m ,b_l ) \ge 
$$ 
$$ 
\ge d(a_i ,a_k ) + d(a_k ,a_m ) + d(a_m , a_j )\ge d(a_i , a_j ) . 
$$ 
If 
$$
d(a_i ,b_l )= d(a_i ,a_k )+ d(a_k ,b_l )\mbox{ with } k\le q \mbox{ and } 
$$ 
$$ 
d(a_j ,b_l )= d(a_j ,a_m )+ \varepsilon + d(b_m ,b_l ) \mbox{ with } m>q  , 
$$ 
then $j> q <l$ (otherwise we are in the previous situation) and 
$$ 
d(a_i ,b_l ) + d(a_j ,b_l ) \ge  d(a_i ,a_k )+ d(a_k ,b_l ) + d(a_j ,a_m )+ \varepsilon + d(b_m ,b_l ) \ge 
$$ 
$$ 
\ge d(a_i ,a_k ) + d(a_k ,a_l ) -\varepsilon  + d(a_m , a_j ) + d(a_m ,a_l ) . 
$$ 
By the assumptions of the theorem  $a_l$ is not between $a_k$ and $a_m$ and 
$$
\varepsilon  < d(a_k , a_l ) + d(a_m ,a_l ) - d (a_k ,a_m ) . 
$$ 
Thus   
$$
d(a_i ,a_k )+ d(a_k ,a_l ) -\varepsilon + d(a_m ,a_j ) +  d(a_m ,a_l ) \ge   
d(a_i ,a_k )+ d(a_k ,a_m ) + d(a_m ,a_j ) \ge d(a_i ,a_j ) . 
$$ 
Now assume that  
$$
d(a_i ,b_l )= d(a_i ,a_k )+  \varepsilon + d(b_k ,b_l )\mbox{ with } k> q \mbox{ and } 
$$ 
$$ 
d(a_j ,b_l )= d(a_j ,a_m )+ \varepsilon + d(b_m ,b_l ) \mbox{ with } m>q  .  
$$ 
We may assume that $q< min  (i,j,l)$ (otherwise we are in one of previous situations). 
Then  
$$ 
d(a_i ,b_l ) + d(a_j ,b_l ) \ge  d(a_i ,a_k )+ \varepsilon + d(b_k ,b_l ) + d(a_j ,a_m )+ \varepsilon + d(b_m ,b_l ) \ge 
$$ 
$$ 
\ge d(a_i ,a_k ) + d(a_k ,a_l ) + d(a_m , a_j ) + d(a_m ,a_l ) \ge d(a_i ,a_j ) . 
$$ 

Case 4.   $d(b_i ,b_j ) \le d(b_i ,a_l )+ d(a_l ,b_j )$.
This case is similar to Case 3. 
Note that we can use the inequality 
$\varepsilon  < d(b_k , b_l ) + d(b_m ,b_l ) - d (b_k ,b_m )$ 
for  $k\le q <min(l,m)$. 
$\Box$

\bigskip

The following corollary states the property that 
if $A$ and $B$ are finite subspaces of $\mathbb{U}$  
which are sufficiently similar then $B$ has a  
copy $B'$ in $\mathbb{U}$ under an isometry fixing 
$A\cap B$ so that $A$ and $B'$ are sufficiently 
close in $\mathbb{U}$. 
The additional technical assumptions on $A$ 
appearing in the formulation can be easily 
satisfied in a very close subspace $A'$. 
In Section 6.3 we give an example of such a construction.

\begin{cor}  \label{CorUrysohn}
Let $A = \{ a_1 ,...,a_n \} \subset \mathbb{U}$, $0\le q <n$ and 
$\varepsilon >0$ satisfy all inequalities of the form 
$$ 
4\varepsilon < d(a_i ,a_j ) \mbox{ for pairs } i<j\le n \mbox{ with } q<j \mbox{ and } 
$$ 
$$  
4\varepsilon < d(a_i ,a_j ) + d(a_i ,a_k ) - d(a_j ,a_k ) 
\mbox{ for triples } a_i ,a_j ,a_k \mbox{ with } |\{ i,j,k\}| =3  
\mbox{ , } k\le q<i,j . 
$$  

Let $B$ be an $n$-element metric space of diameter $\le 1$ 
consisting of elements $b_i$ so that for each pair $i<j\le n$, 
$|d(b_i ,b_j )-d(a_i ,a_j )|\le \varepsilon$. 
We assume that $a_1 =b_1$,...,$a_q =b_q$ and the metric defined 
on $\{ b_1 ,...,b_q\}$ in the space $B$ coincides with the metric 
defined on $\{ a_1 ,...,a_q\}$ in the space $\mathbb{U}$.   

Then the space $B$ embeds into $\mathbb{U}$ over 
$\{ a_1 ,...,a_q\}$ so that for each $q<i\le n$, 
$d(a_i ,b_i ) =\varepsilon$. 
\end{cor} 

The corollary follows from Theorem \ref{LemUrysohn} and 
the properties of $\mathbb{U}$ that $\mathbb{U}$ is universal and ultrahomogeneous. 


\subsection{Distance between types over the Urysohn space}

In this section we present a corollary of Theorem \ref{LemUrysohn} 
which evauates distances between types over the Urysohn space. 
We now give some preliminaries concerning continuous logic 
in the case of $\mathbb{U}$. 

Atomic formulas are of the form $d(t_1 ,t_2 )$ where $t_1$ and $t_2$ 
are variables. 
Below we will also consider expansions of $\mathbb{U}$ by constants. 
In this case $t_1$ and $t_2$ can also be constants. 
{\bf Statements} concerning metric structures are usually 
formulated in the form 
$$
\phi = 0, 
$$ 
where $\phi$ is a {\bf formula}, i.e. an expression built from 
0,1 and atomic formulas by applications of the following functions: 
$$ 
x/2  \mbox{ , } x\dot- y= max (x-y,0) \mbox{ , } min(x ,y )  \mbox{ , } max(x ,y )
\mbox{ , } |x-y| \mbox{ , } 
$$ 
$$ 
\neg (x) =1-x \mbox{ , } x\dot+ y= min(x+y, 1) \mbox{ , } sup_x \mbox{ and } inf_x  
$$ 
(here $sup_x$ and $inf_x$ play the role of quantifiers). 
By $Th(\mathbb{U})$ we denote the {\bf continuous theory} of $\mathbb{U}$ 
which is the set of all statements without free variables satisfied in $\mathbb{U}$.  

For every $c_1 ,...,c_n \in \mathbb{U}$ and $A\subseteq \mathbb{U}$ 
we define the $n$-type $tp(\bar{c}/A)$ of $\bar{c}$ over $A$ 
as the set of all $\bar{x}$-conditions with parameters from $A$ 
which are satisfied by $\bar{c}$ in $M$.  
Let $S_n (T_A )$ be the set of all $n$-types over $A$ 
of the expansion of the theory $Th(\mathbb{U})$ 
by constants from $A$. 

The $d$-topology is defined by the metric 
$$
d(p,q)= inf \{ max_{i\le n} d(c_i ,b_i )| \mbox{ there is a model } M \models Th(\mathbb{U}) \mbox{ with } M\models p(\bar{c})\wedge q(\bar{b})\}. 
$$ 
By Propositions 8.7 and 8.8 of \cite{BYBHU} for any theory 
the $d$-topology is finer than so called the logic topology and 
$(S_n (T_A ),d)$ is a complete space. 
The $d$-topology coincides with the logic topology when the theory 
is separably categorical. 
This is the case of $Th(\mathbb{U})$. 
Moreover $Th(\mathbb{U})$ admits elimination of quantifiers 
\cite{U}. 

\begin{thm} \label{Udist} 
Let $A_0 =\{ a_1 ,...,a_q \}$ be a subset of $\mathbb{U}$ of size $q$. 
Let $n>q$. 
Consider subspaces $A= \{ a_{1}, ...,a_q ,..., a_{n} \}$ and 
$B = \{ b_1 ,..., b_{q}, ...,b_{n} \}$, where  
$a_1 = b_1 ,...,a_q = b_q$. 
Assume 
$max \{ |d(b_i ,b_j )-d(a_i ,a_j )| : 1\le i< j\le n \} \le \varepsilon$. 

Then the distance between types $tp ( a_{q+1}, ...,a_{n}/A_0 )$ 
and  $tp ( b_{q+1}, ...,b_{n} /A_0 )$ is not greater than 
$18\varepsilon$
\end{thm} 

{\em Proof.} 
We start with a procedure which replaces $A$ 
by a minimal subset $A'\subset A$ containing $A_0$ 
so that the distances between pairs of elements 
of $A'$ with at least one of them from $A'\setminus A_0$
are $> 4\varepsilon$ and 
$A$ is contained in the neighbourhood of $A'$ 
of radius $4\varepsilon$. 
Let $B'$ be the subset of $B$  consisting of 
elements with the same numbers 
as elements of $A'$ in $A$. 
Note that the distances between pairs 
of elements of $B'$ are $> 3\varepsilon$ 
and $B$ is contained in the neighbourhood of $B'$ 
of radius $5\varepsilon$. 
We consider $A'$ and $B'$ under the enumerations induced by 
enumerations of $A$ and $B$. 

Now consider all $a_i \in A' \setminus A_0$ which appear 
in triples $a_i ,a_j ,a_k$ in $A'$ with  
$$ 
 d(a_i ,a_j ) + d(a_i ,a_k ) - d(a_j ,a_k ) \le 2\varepsilon 
\mbox{ , where }  |\{ i,j,k\}| =3 \mbox{ , } a_k \in A_0  
\mbox{ , } a_i ,a_j \not\in A_0 .  
$$ 
Firstly find indices of such elements in $A'$ and $B'$. 
Then for each index $i_j$ of this set apply free amalgamation 
of $A'$ with the two-element subspace $\{ a_{i_j}, a'_{i_j} \}$ 
where the distance is rational and satisfies   
$$
2\varepsilon \le d(a_{i_j} , a'_{i_j} ) < 4\varepsilon . 
$$
We repeat this procedure for each element $a_{i_j}$ of our list. 
We use Theorem 2.1 of \cite{bogaty} (see Introduction above) 
for every amalgamation. 
Removing $a_{i_0},...,a_{i_j},...$ from the obtained set 
we construct  $A''$ which already satisfies the 
assumptions of Corollary \ref{CorUrysohn}. 
Indeed, it is now obvious that if $a_i ,a_j ,a_k$ is a triple from $A'$  
as above, then $a'_i ,a_j , a_k$  have the property 
$$ 
 d(a'_i ,a_j ) + d(a'_i ,a_k ) - d(a_j ,a_k ) > 4\varepsilon .  
$$  
Since $d(a_j ,a_k ) \ge max (d(a_i ,a_j ) ,d(a_i ,a_k ))$ 
the permutation of $a'_i$ and $a_j$  does not change this property. 
Moreover it is easy to see that no element of $A'\cap A''$ 
plays the role of $a_i$ in any triple of $A''$ intersecting $A_0$.  
Such an assumption whould imply the same property in $A'$ 
(possibly replacing some $a'_{i_j}$ by $a_{i_j}$).  
This follows from the property of the space obtained after 
all amalgamations above, that any non-trivial 
path from $a'_{i_j}$ contains $a_{i_j}$.    

We apply the same procedure to $B'$ where 
we put $d(b_{i_j}, b'_{i_j}) =  d(a_{i_j} , a'_{i_j} )$. 
As a result we obtain the corresponding $B''$. 

We consider $A''$ under the enumeration induced by $A'$ 
where the number of every $a'_{i_j}$ is just $i_j$. 
Note that the distance between 
elements of $A'$ and $A''$ with the same numbers 
(i.e for example $d(a_{i_j} , a'_{i_j})$)  is 
not greater than $4\varepsilon$. 
In particular $A$  is contained in the neighbourhood 
of $A''$ of radius $8\varepsilon$. 
Thus $B$  is contained in the neighbourhood 
of $B''$ of radius $9\varepsilon$.

On the other hand it is easy to see that we still have 
$$ 
max \{ |d(b_i ,b_j )-d(a_i ,a_j )| : 1\le i< j\le n \} \le \varepsilon , 
$$
for the corresponding elements of $A''$ and $B''$. 
As a result we have $A''$ and $B''$ which satisfy 
the assumptions of Corollary \ref{CorUrysohn}.   
Applying Corollary \ref{CorUrysohn} to $A''$ and $B''$ 
over $A_0$ we obtain an embeddng of $B''$ into 
$\mathbb{U}$ (with the same denotation $B''$), 
so that the distances between the corresponding 
elements with numbers $>q$ is $\varepsilon$.
Using the fact that $\mathbb{U}$ is ultrahomogeneous 
we obtain an embedding of $B$ into $\mathbb{U}$,  
say $\hat{B}$, which is distant from $B''$ 
by $\le 9\varepsilon$. 
The rest is clear. 
$\Box$


\bigskip

Institute of Mathematics, University of Wroc{\l}aw, \parskip0pt

pl.Grunwaldzki 2/4, 50-384 Wroc{\l}aw, Poland \parskip0pt

E-mail: ivanov@math.uni.wroc.pl

\end{document}